\documentclass[12pt,a4paper]{article}

\usepackage[T1]{fontenc}
\usepackage[utf8]{inputenc}
\usepackage{geometry}
\usepackage{amsmath,amssymb,amsthm,mathtools}

\usepackage{hyperref}

\hypersetup{
  hidelinks,
  pdftitle={Some computed examples of the Noncommutative Foliation Invariant (NCFI)},
  pdfauthor={Ioannis P. Zois},
  pdfsubject={Noncommutative geometry, foliations, cyclic cohomology, K-theory},
  pdfkeywords={foliation C*-algebra, transverse fundamental cyclic cocycle, Kronecker foliation, noncommutative torus}
}

\newcounter{definition}
\newcounter{proposition}
\newcounter{propositionvariant}
\newcounter{remark}
\newcounter{lemma}
\newcounter{theorem}
\newcounter{corollary}
\newcounter{note}

\renewcommand{\thedefinition}{\arabic{definition}}
\renewcommand{\theproposition}{\arabic{proposition}}
\renewcommand{\thepropositionvariant}{\theproposition}
\renewcommand{\theremark}{\arabic{remark}}
\renewcommand{\thelemma}{\arabic{lemma}}
\renewcommand{\thetheorem}{\arabic{theorem}}
\renewcommand{\thecorollary}{\arabic{corollary}}

\newcommand{\DefHead}[1]{%
  \refstepcounter{definition}%
  \par\medskip\noindent\textbf{Definition \thedefinition.}\label{#1}\ %
}

\newcommand{\PropHead}[1]{%
  \refstepcounter{proposition}%
  \par\medskip\noindent\textbf{Proposition \theproposition.}\label{#1}\ %
}

\newcommand{\PropHeadVariant}[2]{%
  \begingroup
  \renewcommand{\thepropositionvariant}{\theproposition#2}%
  \refstepcounter{propositionvariant}%
  \par\medskip\noindent
  \textbf{Proposition \theproposition#2.}\label{#1}\ %
  \endgroup
}

\newcommand{\PropHeadPrime}[1]{%
  \PropHeadVariant{#1}{\ensuremath{{}^{\prime}}}%
}

\newcommand{\PropHeadDoublePrime}[1]{%
  \PropHeadVariant{#1}{\ensuremath{{}^{\prime\prime}}}%
}

\newcommand{\RemHead}[1]{%
  \refstepcounter{remark}%
  \par\medskip\noindent\textbf{Remark \theremark.}\label{#1}\ %
}

\newcommand{\LemHead}[1]{%
  \refstepcounter{lemma}%
  \par\medskip\noindent\textbf{Lemma \thelemma.}\label{#1}\ %
}

\newcommand{\ThmHead}[1]{%
  \refstepcounter{theorem}%
  \par\medskip\noindent\textbf{Theorem \thetheorem.}\label{#1}\ %
}

\newcommand{\CorHead}[1]{%
  \refstepcounter{corollary}%
  \par\medskip\noindent\textbf{Corollary \thecorollary.}\label{#1}\ %
}

\newcommand{\ProofHead}{%
  \par\medskip\noindent\textit{Proof.}\ %
}

\title{The Noncommutative Foliation Invariant (NCFI): extension to the odd codimension case and computed examples}
\author{Ioannis P. ZOIS, Exeter College, OX1 3DP, Oxford, UK\\
email: i.zois@exeter.oxon.org}

\date{\today}

\begin{document}

\maketitle

\begin{abstract}
This article extends the definition of the Noncommutative Foliation Invariant (NCFI) for foliations of odd codimension and computes certain key examples for both even and odd codimension cases: fibrational foliations, irrational Kronecker, rational Kronecker (both the vertical and the horzontal foliation using a flat connection) and weighted Hopf/orbifold cases. We also prove some more general results along the way. More concretely,
we compute the Noncommutative Foliation Invariant (NCFI) for several basic
families of foliated manifolds.  In even codimension the invariant is the
Chern--Connes pairing of Connes' transverse fundamental cyclic cocycle with the
\(K_0\)-class associated with the transverse geometric module.  In odd
codimension the transverse cocycle has odd degree, so a numerical pairing also
requires a specified odd \(K_1\)-class; the odd-codimensional value is therefore
an invariant of the foliation together with this chosen odd-favourable
structure.

For foliations defined by fibrations over even-dimensional bases, the invariant reduces to
an ordinary characteristic-number computation on the base.  The control example
\(S^2\times S^1\to S^2\) gives zero, while the codimension-four fibration
\(\mathbb{CP}^2\times S^1\to\mathbb{CP}^2\) gives
\[
Z(\mathcal F^{\mathrm v}_4)=3.
\]
For irrational Kronecker foliations on \(T^2\), reduction to a complete
transversal gives the irrational rotation algebra
\(A_\theta=C(S^1_z)\rtimes_\alpha\mathbb Z\).  The transported transverse
fundamental cyclic cocycle is \(\psi^{(2)}_1\).  The natural return-map,
Connes--Thom and longitudinal Dirac constructions select the class \([U]\), and
therefore
\[
Z(\mathcal F_\theta;[U])=0,
\]
whereas the transversal coordinate class \([V]\) would pair to \(1\).  For the
rational vertical and horizontal Kronecker cases considered here, the selected
longitudinal odd classes are zero, so both corresponding NCFI values vanish.

Finally, for the weighted Hopf foliation of \(S^5\) with weights \((1,2,3)\),
the quotient is the weighted projective orbifold
\(\mathbb P_{\mathrm{orb}}(1,2,3)\).  This gives a non-fibrational
even-codimensional example with
\[
Z(\mathcal F_{1,2,3})=\frac73.
\]
\end{abstract}

\section{Introduction}

\subsection{Review of the Basic Ideas}

In \cite{zois2000} we introduced the
\emph{Noncommutative Foliation Invariant}, abbreviated here to \emph{NCFI}.
The original motivation came from non-linear \(\sigma\)-models and string
theory in physics, while the construction is formulated in noncommutative geometry:
holonomy groupoids, foliation \(C^*\)-algebras, \(K\)-theory and cyclic
cohomology; see \cite{Connes94,ConnesSurveyFoliations}.

Let \((M,\mathcal F)\) be a compact smooth \(m\)-manifold equipped with a
foliation \(\mathcal F\) of codimension \(q\), where \(1\leq q<m\) with $q,m\in\mathbb{N}$, given by an \((m-q)\)-dimensional \emph{integrable} \textsl{subbundle}
\[
\mathcal F\subset TM 
\]
of the tangent bundle $TM$ of $M$. Its maximal connected immersed integral submanifolds are the leaves. All leaves have the same dimension $(m-q)$ but they may have different topologies. The
transverse bundle of the foliation is
\[
t:=TM/\mathcal F,
\]
and its complexification is
\[
t_{\mathbb C}:=t\otimes \mathbb C .
\]
This notation is used throughout the paper.

To \((M,\mathcal F)\) one associates the holonomy groupoid
\[
G_{\mathcal F},
\]
or simply \(G\) when no ambiguity arises. The source and range maps are
\[
s,r:G\rightrightarrows M .
\]
The groupoid \(G_{\mathcal F}\) records both the leaf-equivalence relation and
holonomy transport along leafwise paths. For the construction and basic
properties of holonomy groupoids, see
\cite{Connes94,ConnesSurveyFoliations,Renault80,Paterson1999}.

The smooth convolution algebra of compactly supported half-densities on the
holonomy groupoid is
\[
A_{\mathrm{sm}}
   :=
C_c^\infty(G_{\mathcal F},\Omega^{1/2}).
\]
When the foliation is fixed, we also write
\[
C_c^\infty(G,\Omega^{1/2})
\]
for the same algebra. Its reduced \(C^*\)-completion is
\[
A:=C_r^*(M,\mathcal F)=C_r^*(G_{\mathcal F}).
\]
Thus the basic passage is
\[
(M,\mathcal F)
\longmapsto
G_{\mathcal F}
\longmapsto
C_c^\infty(G_{\mathcal F},\Omega^{1/2})
\longmapsto
C_r^*(G_{\mathcal F}).
\]
We use the reduced foliation algebra throughout. Subsection~1.4 fixes the
conventions concerning smooth algebras, reduced completions, dense Banach
subalgebras and Morita-equivalent models.

The \(K\)-theory of \(A\) provides the \(K\)-classes used in the invariant.
Connes' transverse fundamental cyclic cocycle, abbreviated to \emph{TFCC}, is
defined on the smooth convolution algebra \(A_{\mathrm{sm}}\). If the foliation
has codimension \(q\) and is transversally oriented, the TFCC is a cyclic
\(q\)-cocycle
\[
\varphi_{\mathcal F}\in HC^q(A_{\mathrm{sm}}),
\]
representing the noncommutative transverse fundamental class of the foliation;
see \cite{ConnesTFClass,Connes94}. It generalises the ordinary fundamental
class of a smooth transverse manifold.

We recall the \emph{transverse geometric module}, abbreviated to \emph{tgm},
constructed in \cite{zois2000}. Let \(G=G_{\mathcal F}\) denote the holonomy groupoid of the foliation and fix a smooth
Euclidean metric \(g\) on the transverse bundle
\[
t=TM/\mathcal F .
\]
The metric induces a Hermitian inner product on \(t_{\mathbb C}\), taken to be
anti-linear in the first variable. Consider the pull-back bundle
\[
r^*(t_{\mathbb C})\longrightarrow G
\]
by the range map. The smooth core of the construction is
\[
E_{\mathrm{sm}}
   :=
C_c^\infty\bigl(G,\Omega^{1/2}\otimes r^*(t_{\mathbb C})\bigr).
\]
where $\Omega ^{1/2}$ denotes half-densities. Thus the tgm starts from compactly supported smooth half-density sections over
the holonomy groupoid with coefficients in the pulled-back complexified
transverse bundle.

For \(y\in M\), write
\[
G^y:=r^{-1}(y).
\]
All integrals below are over range fibres; the half-density formalism absorbs
the Haar-system notation.

The right convolution action of
\[
C_c^\infty(G,\Omega^{1/2})
\]
on \(E_{\mathrm{sm}}\) is
\[
  (\xi f)(\gamma)
    =
    \int_{G^{r(\gamma)}}
      \xi(\gamma_1)\,
      f(\gamma_1^{-1}\gamma),
\]
for
\[
\xi\in C_c^\infty\bigl(G,\Omega^{1/2}\otimes r^*(t_{\mathbb C})\bigr),
\qquad
f\in C_c^\infty(G,\Omega^{1/2}).
\]
The \(A\)-valued inner product on \(E_{\mathrm{sm}}\) is
\[
  \langle \xi,\eta\rangle(\gamma)
    =
    \int_{G^{r(\gamma)}}
      \big\langle
        \xi(\gamma_1^{-1}),
        \eta(\gamma_1^{-1}\gamma)
      \big\rangle,
\]
where the inner product inside the integral is the Hermitian inner product on
the fibre of \(t_{\mathbb C}\) induced by \(g\). Completing \(E_{\mathrm{sm}}\)
for the norm
\[
  \|\xi\|
  =
  \big\|\langle \xi,\xi\rangle\big\|^{1/2}_{C_r^*(M,\mathcal F)}
\]
gives a Hilbert \(A\)-module
\[
E.
\]
For the general Hilbert \(C^*\)-module background, see \cite{Lance95}.

\DefHead{def:tgm}
The Hilbert \(A\)-module \(E\) constructed above is called the
\emph{transverse geometric module} of the foliation, or simply the
\emph{tgm}.

We next recall how the associated even \(K\)-class is obtained. Let
\[
h(\gamma):t_{s(\gamma)}\longrightarrow t_{r(\gamma)}
\]
be the linear holonomy action on the transverse bundle, induced by the
differential of the holonomy transformation. Using \(h\), one defines a left
action on \(E_{\mathrm{sm}}\) by
\[
  (f\xi)(\gamma)
    =
    \int_{G^{r(\gamma)}}
      f(\gamma_1)\,
      h(\gamma_1)\,
      \xi(\gamma_1^{-1}\gamma),
\]
for
\[
f\in C_c^\infty(G,\Omega^{1/2}),
\qquad
\xi\in E_{\mathrm{sm}}.
\]
Let \(\lambda(f)\) denote the corresponding operator on the smooth core. Its
formal adjoint is not, in general, given by the ordinary involution on the
scalar convolution algebra. With the metric on \(t\), it has the form
\[
  (\lambda(f)^*\xi)(\gamma)
    =
    \int_{G^{r(\gamma)}}
      f^\#(\gamma_1)\,
      h(\gamma_1)\,
      \xi(\gamma_1^{-1}\gamma),
\]
where
\[
  f^\#(\gamma)
  :=
  \widetilde f(\gamma^{-1})\,\Delta(\gamma),
\]
and
\[
  \Delta(\gamma)
  =
  h(\gamma)^{-t}h(\gamma)^{-1}
  \in
  \operatorname{End}\bigl(t_{\mathbb C}(r(\gamma))\bigr).
\]
Here \(\widetilde f\) denotes complex conjugation of the scalar half-density
coefficient, and the exponent \(^{-t}\) denotes inverse transpose with respect to the
chosen transverse metric.

Unless the metric on \(t\) is holonomy-invariant, the left action is not a
\(*\)-representation of the reduced foliation \(C^*\)-algebra. The obstruction
is the modular correction
\[
\Delta(\gamma)=h(\gamma)^{-t}h(\gamma)^{-1}.
\]
Thus the assignment \(f\mapsto \lambda(f)\) need not extend to a bounded
\(*\)-homomorphism
\[
A\longrightarrow \operatorname{End}_A(E).
\]

The construction in \cite{zois2000} passes to the domain of the closure of the
left action and equips it with the graph norm
\[
\|x\|_\lambda:=\|x\|+\|\lambda(x)\|.
\]
The graph-norm completion gives a dense Banach \(*\)-subalgebra
\[
B\subset A.
\]
The smooth module \(E_{\mathrm{sm}}\) has a corresponding graph-norm
completion, denoted
\[
E_B,
\]
which is the \(B\)-module used in the \(K\)-theoretic construction of
\cite{zois2000}. Under the hypotheses of that construction, \(E_B\) determines
a stable projective-module class
\[
[e_B]\in K_0(B).
\]
Its image under the inclusion-induced map
\[
K_0(B)\longrightarrow K_0(A)
\]
is denoted by
\[
[e_{\mathcal F}]\in K_0(A).
\]
The data \(B\) and \(E_B\) are part of the construction; the Hilbert
\(A\)-module \(E\) alone does not determine this projective \(B\)-module
class. For the general \(K\)-theory background, see \cite{Blackadar98}.

\DefHead{def:even-k-class}
The class
\[
[e_{\mathcal F}]\in K_0(A)
\]
constructed above is the even \(K\)-class associated with the transverse
geometric module of the foliation.

\DefHead{def:ncfi-even}
Assume that the foliation \((M,\mathcal F)\) is transversally oriented and has
even codimension
\[
q=2n.
\]
Let
\[
\varphi_{\mathcal F}\in HC^{2n}(A_{\mathrm{sm}})
\]
be Connes' transverse fundamental cyclic cocycle, represented on the chosen
smooth convolution algebra or on a smooth Morita-equivalent model.  Let
\[
[e_{\mathcal F}]\in K_0(A),
\qquad
A=C_r^*(M,\mathcal F),
\]
be the even \(K\)-class associated with the transverse geometric module.

We assume, as part of the data of the pairing, that \([e_{\mathcal F}]\) is
represented by an idempotent or projection over the chosen smooth algebra on
which \(\varphi_{\mathcal F}\) is defined, possibly after replacing
\(A_{\mathrm{sm}}\) by a holomorphically closed smooth dense subalgebra or by a
smooth Morita-equivalent model.  Equivalently, the class used in the formula
below is the image in \(K_0(A)\) of a smooth \(K_0\)-class on which the
Chern--Connes pairing is defined.

Under this smooth-representability hypothesis, the
\emph{Noncommutative Foliation Invariant} in even codimension is
\[
Z(\mathcal F)
:=
\langle \varphi_{\mathcal F},[e_{\mathcal F}]\rangle .
\]
In all examples computed below, the required smooth representatives are
displayed explicitly or obtained through the stated smooth Morita-equivalent
models.  See
\cite{zois2000,Connes94,Connes80,Connes85,ConnesMoscovici95}
for the definition and the background on Connes' pairing.

This definition is direct only when \(q\) is even. The tgm construction exists
in both even and odd codimension and gives an even \(K_0\)-class of the reduced
foliation \(C^*\)-algebra. Connes' pairing, however, matches even cyclic
cohomology with \(K_0\) and odd cyclic cohomology with \(K_1\). Hence, when
\(q\) is odd, the tgm class and the transverse fundamental cyclic cocycle have
opposite parity. The odd-codimensional case is treated in the next subsection.

Before proceeding with the \(K\)-theoretic parity fixing, we indicate the
organization of the article. The remaining subsections of Section~1 explain
the odd-codimensional parity correction and fix the algebraic conventions used
throughout. Section~2 treats fibrational examples in even codimension, where
the foliation algebra is Morita equivalent to the algebra of functions on the
base and the NCFI reduces to a characteristic-number computation. Sections~3--6 deal with the Irrational Kronecker foliation of the 2-torus case and compute the NCFI. Sections~7 and~8 treat the rational vertical and horizontal
cases of the Kronecker foliation of the 2-torus. Section~9 computes the NCFI for the Hopf foliation on $S^5$ which is a non-trivial, non-fibrational codim 4 example (even codim, no parity fixing needed). Finally Section~10 summarizes the computed values and the remaining problems.

\subsection{Parity fixing in odd codimension: Route A}

Assume that the foliation has odd codimension
\[
q=2n+1.
\]
Then Connes' transverse fundamental cyclic cocycle has odd degree,
\[
\varphi_{\mathcal F}\in HC^{2n+1}(A_{\mathrm{sm}}),
\]
whereas the transverse geometric module recalled in Subsection~1.1 gives an even class
\[
[e_{\mathcal F}]\in K_0(A).
\]
The even-codimensional pairing is therefore unavailable: Connes' pairing matches even
cyclic cohomology with \(K_0\) and odd cyclic cohomology with \(K_1\); see
\cite{Connes80,Connes85,Connes94,GBVF01,ConnesTFClass}.

This parity obstruction is structural.  For a general \(C^*\)-algebra \(A\) there is no
canonical functorial map
\[
K_0(A)\longrightarrow K_1(A).
\]
Bott periodicity gives the parity shift after suspension,
\[
K_0(A)\cong K_1(SA),
\]
not a canonical map from \(K_0(A)\) to \(K_1(A)\) on the same algebra; see
\cite{Blackadar98,WeggeOlsen1993,RordamLarsenLaustsen2000}.  Thus the tgm class alone
does not determine a distinguished odd class
\[
[u_{\mathcal F}]\in K_1(A).
\]

Nor is there an internal cyclic-cohomological parity change on the same smooth algebra.
For a general smooth algebra \(A_{\mathrm{sm}}\), there is no natural map
\[
HC^q(A_{\mathrm{sm}})
   \longrightarrow
HC^{q\pm 1}(A_{\mathrm{sm}})
\]
sending the odd transverse fundamental cyclic cocycle to an even cyclic class on
\(A_{\mathrm{sm}}\).  Connes' periodicity operator
\[
S:HC^m(A_{\mathrm{sm}})\longrightarrow HC^{m+2}(A_{\mathrm{sm}})
\]
shifts degree by \(2\) and preserves parity.  The operator \(B\) is part of the mixed
\((b,B)\)-complex, but it does not by itself give a canonical parity-fixing operation on
cyclic cohomology classes.  See
\cite{Connes85,Connes94,Loday98,QuillenIHES1988,Cuntz1997,CuntzQuillen,CuntzQuillenInvent1997}.

The tgm construction also contains no hidden canonical odd class.  It is naturally a
right Hilbert \(A\)-module, while the holonomy-induced left action is not, in general, a
\(*\)-representation of \(A\).  With the notation of Subsection~1.1,
\[
f^\#(\gamma)
   =
\widetilde f(\gamma^{-1})\,\Delta(\gamma),
\qquad
\Delta(\gamma)
   =
h(\gamma)^{-t}h(\gamma)^{-1}
   \in
\operatorname{End}\!\bigl(t_{\mathbb C}(r(\gamma))\bigr).
\]
Unless the chosen transverse metric is holonomy-invariant, the left action differs from
a \(*\)-representation by the modular correction \(\Delta\).  One therefore passes to
the graph-norm Banach \(*\)-algebra
\[
B\subset A
\]
and to the graph-norm module used to construct
\[
[e_B]\in K_0(B),
\qquad
[e_{\mathcal F}]\in K_0(A).
\]
Consequently the tgm construction produces the even class used above, but not a
canonical class in \(K_1(A)\); see \cite{Connes94,zois2000}.

\RemHead{rem:transverse-metric-and-odd-route}
The transverse metric used in Subsection~1.1 is already part of the tgm construction.  It
enters the Hermitian structure on
\[
C_c^\infty\!\bigl(G,\Omega^{1/2}\otimes r^*(t_{\mathbb C})\bigr)
\]
and appears explicitly in
\[
\Delta(\gamma)=h(\gamma)^{-t}h(\gamma)^{-1}.
\]
Thus Dirac-type odd classes are natural candidates in favourable examples: they do not
import unrelated geometry.  This observation, however, does not produce a canonical
odd class in \(K_1(A)\).  Such a class still requires the Clifford-theoretic and analytic
hypotheses of unbounded \(KK\)-theory; see
\cite{BaajJulg1990,Lance95,KaadLesch2012,Kaad2020,Dam07}.

In odd codimension one therefore keeps the foliation, the smooth algebra, the TFCC and
the tgm class,
\[
(M,\mathcal F),
\qquad
A_{\mathrm{sm}},
\qquad
\varphi_{\mathcal F},
\qquad
[e_{\mathcal F}]\in K_0(A),
\]
and adds a distinguished odd class
\[
[u_{\mathcal F}]\in K_1(A)
\]
from geometric, dynamical, analytic or \(KK\)-theoretic structure.

\RemHead{rem:odd-favourable-overview}
We shall say that a transversally oriented foliation \((M,\mathcal F)\) of odd
codimension
\[
q=2n+1
\]
is \emph{odd-favourable} if its reduced foliation \(C^*\)-algebra
\[
A=C_r^*(M,\mathcal F)
\]
has a specified distinguished odd class
\[
[u_{\mathcal F}]\in K_1(A)
\]
obtained from one of the standard mechanisms below.  The tgm continues to provide
\([e_{\mathcal F}]\in K_0(A)\); the odd class is extra structure natural to the example.

\medskip

\noindent
\textbf{(a) Crossed-product / return-map case.}
Assume that a complete transversal \(T\subset M\) gives, up to groupoid equivalence, a
transformation groupoid
\[
G_T\simeq T\rtimes \mathbb Z,
\]
or equivalently that, after Morita equivalence, the foliation algebra has a crossed-product
model
\[
A\sim_M C_0(T)\rtimes_\alpha \mathbb Z.
\]
In the unital cases used below, the implementing unitary of the \(\mathbb Z\)-action
defines a class
\[
[U]\in K_1\bigl(C_0(T)\rtimes_\alpha \mathbb Z\bigr),
\]
and hence, after transport through Morita equivalence, a distinguished odd class in
\(K_1(A)\); see
\cite{PimsnerVoiculescu80,Blackadar98,RaeburnWilliams1998,Williams2007}.  This is
the mechanism used later for the irrational Kronecker foliation.

For
\[
A_B:=B\rtimes_\alpha\mathbb Z,
\]
the Pimsner--Voiculescu six-term exact sequence has boundary map
\[
\partial:K_1(A_B)\longrightarrow K_0(B).
\]

If \(B\) is unital, the crossed-product implementing unitary
\(U\in B\rtimes_\alpha\mathbb Z\) satisfies, up to the sign convention for
\(\partial\),
\[
\partial([U])=[1_B].
\]
This boundary equation records the return-map component of \([U]\), but it does
not in general characterize \([U]\) uniquely: classes coming from \(K_1(B)\) lie
in the kernel of \(\partial\).  The distinguished class in the present route is
therefore the actual implementing-unitary class determined by the crossed-product
structure, not merely an arbitrary class with this boundary.

Thus the return-map structure supplies \([U]\in K_1(B\rtimes_\alpha\mathbb Z)\).  It is
not produced by a universal map \(K_0(A)\to K_1(A)\), nor by the transverse geometric
module alone.  In the irrational rotation case,
\[
B=C(S^1),
\qquad
A_\theta\cong C(S^1)\rtimes_\alpha\mathbb Z,
\]
and the distinguished odd class is this implementing-unitary class; see
\cite{PimsnerVoiculescu80,Blackadar98}.

\medskip

\noindent
\textbf{(b) Flow / Connes--Thom case.}
Assume that the foliation algebra has a natural flow model
\[
A\simeq B\rtimes\mathbb R.
\]
The Connes--Thom isomorphism gives
\[
K_j(B)\cong K_{j+1}(A),
\]
so an even class on the \(B\)-side may determine an odd class on the \(A\)-side; see
\cite{ConnesThom1981,Connes94,FackSkandalis1981}.  In favourable dynamical examples
this gives a natural odd class without choosing an explicit longitudinal operator.

\medskip

\noindent
\textbf{(c) Boundary / extension / \(KK\)-theoretic case.}
Assume that the foliation algebra sits in a natural extension, or more generally carries a
distinguished odd Kasparov self-class
\[
x_{\mathcal F}\in KK^1(A,A).
\]
Then, when such a class is available, one may set
\[
[u_{\mathcal F}]
   :=
[e_{\mathcal F}]\otimes_A x_{\mathcal F}
   \in K_1(A).
\]
This is the abstract \(KK\)-theoretic form of the parity correction.  It includes boundary
classes from extensions and \(K\)-oriented correspondences of leaf spaces; see
\cite{Kasparov1981,HilsumSkandalis87,Blackadar98,PimsnerVoiculescu80}.  The class
\(x_{\mathcal F}\) is extra input.

\medskip

\noindent
\textbf{(d) Longitudinal Dirac-type case.}
Assume that the foliation carries longitudinal geometric data giving an odd unbounded
Kasparov cycle.  Concretely, one needs a Hilbert \(A\)-module \(E_{\mathrm{odd}}\), a
nondegenerate left \(*\)-representation \(\pi\), and a densely defined first-order
longitudinal operator
\[
D_\parallel
\]
such that:
\begin{itemize}
\item \(D_\parallel\) is self-adjoint and regular as a Hilbert \(C^*\)-module operator;
\item \([D_\parallel,\pi(a)]\) is bounded on a dense smooth \(*\)-subalgebra;
\item \(\pi(a)(1+D_\parallel^2)^{-1/2}\) is \(A\)-compact.
\end{itemize}
The bounded transform
\[
Y_\parallel
   =
D_\parallel(1+D_\parallel^2)^{-1/2}
\]
defines an odd Kasparov cycle.  Depending on the left representation, the resulting
class lies in \(KK^1(C(M),A)\) or in \(KK^1(A,A)\).  In the former case, and more
generally after capping with the relevant unit class by the Kasparov product, one obtains
an odd class in
\[
K_1(A).
\]
This route is standard but not automatic: the analytic hypotheses must be verified in
each example; see
\cite{BaajJulg1990,Lance95,HigsonRoe00,KaadLesch2012,Kaad2020,Dam07,Vassout,Pierrot2006}.
In settings where self-adjointness is replaced by weaker closure conditions, one may need
half-closed chains and localization techniques; see
\cite{KaadVanSuijlekom2019,VanDenDungen2023,CareyPhillipsRennie06,Rennie}.

\medskip

\noindent
\textbf{(e) Leafwise de Rham / Dirac case from the integrable subbundle
\(\mathcal F\).}
The original geometric datum is the integrable subbundle
\[
\mathcal F\subset TM.
\]
A natural source of odd classes is therefore the leafwise differential geometry: the
leafwise de Rham operator
\[
D_{\mathcal F}:=d_{\mathcal F}+d_{\mathcal F}^*,
\]
or a longitudinal Dirac-type operator associated with the leafwise tangent bundle, after
choosing the required longitudinal metric and Clifford data.  These operators are central
in the longitudinal index theory of foliations; see
\cite{Connes94,ConnesSkandalisLongitudinal,HigsonRoe00,Vassout,BH04,BH08,BenameurFack}.

The parity of such a longitudinal Dirac class is governed by the leaf dimension,
\[
\dim(\mathcal F)=\dim(M)-q,
\]
not by the codimension \(q\).  Thus a leafwise Dirac-type construction gives an odd
class in the usual way when the relevant longitudinal operator has odd parity; in
particular, this occurs for the one-dimensional leaves in the irrational Kronecker case.
To obtain a class in \(K_1(A)\), one first obtains an odd Kasparov class, typically in
\(KK^1(C(M),A)\), or in \(KK^1(A,A)\) when a left \(A\)-representation is present, and
then applies the appropriate Kasparov product.  Hence the leafwise de Rham / Dirac
construction is a natural favourable case, not a universal parity-fixing mechanism for all
odd-codimensional foliations.

One may also use transverse or basic Dirac operators for Riemannian foliations with
bundle-like metrics.  This imposes stronger hypotheses than the longitudinal route and
is not part of the general odd-codimension definition used here; see \cite{HR09,S12}.

\medskip

\noindent
\textbf{(f) Intrinsic modular corrections.}
The tgm construction carries the modular correction
\[
\Delta(\gamma)=h(\gamma)^{-t}h(\gamma)^{-1},
\]
which measures the failure of the left action to be a \(*\)-representation.  This modular
data is geometrically meaningful, but it is not used here as a parity-fixing theorem.  In
particular, we do not assume that \(\Delta\) functorially produces a canonical class
\[
[u_{\mathcal F}]\in K_1(A)
\]
for arbitrary odd-codimensional foliations.  In codimension one \(\Delta\) is scalar, but
it still need not produce the relevant odd \(K\)-class.  For the irrational Kronecker
foliation, the odd class used below comes from the crossed-product return-map structure,
not from the modular correction.

\medskip

All mechanisms above have the same logical form: the transverse geometric module
provides
\[
[e_{\mathcal F}]\in K_0(A),
\]
while the odd class
\[
[u_{\mathcal F}]\in K_1(A)
\]
comes from additional geometric, dynamical, analytic or \(KK\)-theoretic structure.  In
practice, the lowest-cost source is often longitudinal.  If that is unavailable, one may use
a Bott or Clifford parity shift, a boundary/extension class, or stronger transverse
hypotheses such as a basic Dirac operator.  This hierarchy is a guide to natural odd
classes, not a theorem producing them in all cases.

\DefHead{def:ncfi-odd}
Let \((M,\mathcal F)\) be a transversally oriented foliated manifold of odd
codimension
\[
q=2n+1.
\]
Let
\[
\varphi_{\mathcal F}\in HC^{2n+1}(A_{\mathrm{sm}})
\]
be Connes' transverse fundamental cyclic cocycle, represented on the chosen
smooth convolution algebra or on a smooth Morita-equivalent model.  Assume that
\((M,\mathcal F)\) is odd-favourable, so that a distinguished odd class
\[
[u_{\mathcal F}]\in K_1(A),
\qquad
A=C_r^*(M,\mathcal F),
\]
has been specified by one of the mechanisms above.

We assume, as part of the odd pairing, that \([u_{\mathcal F}]\) is represented
by a smooth invertible or unitary over the smooth algebra on which
\(\varphi_{\mathcal F}\) is defined, possibly after passing to a
holomorphically closed smooth dense subalgebra or to a smooth
Morita-equivalent model.  Equivalently, the \(K_1\)-class used in the formula
below is represented in a smooth model carrying the cyclic cocycle.

Under this smooth-representability hypothesis, the
\emph{Noncommutative Foliation Invariant} in odd codimension is
\[
Z(\mathcal F)
   :=
\langle \varphi_{\mathcal F},[u_{\mathcal F}]\rangle .
\]
This is the primary definition used in the present article for odd codimension.

\RemHead{rem:ncfi-odd-favourable}
In odd codimension, the NCFI is an invariant of the foliation together with the chosen
odd-favourable structure.  Thus, when no confusion is possible, we write \(Z(\mathcal F)\);
more precisely, the paired object is
\[
(\mathcal F,[u_{\mathcal F}]).
\]

\RemHead{rem:kronecker-odd-favourable}
For the irrational Kronecker foliation treated below, the relevant odd-favourable
structure is the complete-transversal crossed-product model
\[
A_\theta\cong C(S^1)\rtimes_\alpha\mathbb Z.
\]
The distinguished odd class is the crossed-product return class
\[
[u_\theta]=[U]\in K_1(A_\theta),
\]
where \(U\) is the implementing unitary of the \(\mathbb Z\)-action.

This statement is relative to the complete-transversal crossed-product model.  Different
complete transversals give Morita-equivalent reduced groupoids and canonically
identified \(K\)-theory after transport through the induced Morita equivalences.  They
should not, in general, be described as literally unitarily equivalent crossed-product
models.  The class \([U]\) is canonical relative to the return-map structure encoding the
holonomy dynamics.

The irrational Kronecker foliation is also leafwise one-dimensional, so the leafwise
de Rham route can be compared with the crossed-product route.  The definition used
below, however, selects the crossed-product return class \([U]\).

There is also a cyclic-cohomological reformulation.  Once
\([u_{\mathcal F}]\in K_1(A)\) has been chosen, one may pass to the suspension and use
the cylinder extension
\[
0\longrightarrow SA\longrightarrow A^I
   \xrightarrow{(\mathrm{ev}_0,\mathrm{ev}_1)}
A\oplus A
\longrightarrow 0.
\]
Excision in periodic cyclic cohomology gives connecting morphisms for extensions; see
\cite{CuntzQuillen,CuntzQuillenInvent1997,Brodzki2013}.  Applied to the cylinder
extension, this gives a Chern--Simons-type transgression of the odd cyclic cocycle, in
the sense of \cite{QuillenIHES1988,Connes94}.  The original odd pairing can then be
rewritten as an even--even pairing after Bott periodicity.  This is Route B and it will be described in the next Subsection.

\subsection{Route B: suspension, transgression, and equality with Route A}

Route A defines the odd-codimensional NCFI once a distinguished odd class
\[
[u_{\mathcal F}]\in K_1(A)
\]
has been specified:
\[
Z(\mathcal F)
   :=
\langle \varphi_{\mathcal F},[u_{\mathcal F}]\rangle .
\]
Route B rewrites this same odd pairing as an even--even pairing on the suspension
algebra. It does not construct the odd class, and it does not replace the Route A
definition.

Let
\[
SA:=C_0((0,1))\otimes A
\]
be the suspension of \(A\), and let
\[
A_I:=C([0,1],A)
\]
be the cylinder algebra. At the smooth level we use
\[
A_{I,\mathrm{sm}}
   :=
C^\infty([0,1],A_{\mathrm{sm}})
\]
and
\[
SA_{\mathrm{sm}}
   :=
\{f\in C^\infty([0,1],A_{\mathrm{sm}}):f(0)=f(1)=0\}.
\]
Thus the cylinder extension is
\[
0\longrightarrow SA
  \longrightarrow A_I
  \xrightarrow{(\mathrm{ev}_0,\mathrm{ev}_1)}
  A\oplus A
  \longrightarrow 0,
\]
with smooth analogue
\[
0\longrightarrow SA_{\mathrm{sm}}
  \longrightarrow A_{I,\mathrm{sm}}
  \xrightarrow{(\mathrm{ev}_0,\mathrm{ev}_1)}
  A_{\mathrm{sm}}\oplus A_{\mathrm{sm}}
  \longrightarrow 0.
\]

\RemHead{rem:periodic-convention}
In this subsection the exactness, suspension and transgression statements are
understood in periodic cyclic cohomology.  To keep notation compatible with the
rest of the paper, we continue to write \(HC^k\) for cocycle representatives of
degree \(k\), but the functorial exactness statements are statements in
\(HP^\bullet\).  The pairings below are periodic Chern--Connes pairings
computed on smooth representatives in the sense fixed in Definitions~\ref{def:ncfi-even}
and~\ref{def:ncfi-odd}.  The relevant excision and suspension results are those
of \cite{CuntzQuillen,CuntzQuillenInvent1997,Brodzki2013}.

Let
\[
\varphi_{\mathcal F}\in HC^{2n+1}(A_{\mathrm{sm}})
\]
be the odd transverse fundamental cyclic cocycle. In the endpoint algebra
\(A_{\mathrm{sm}}\oplus A_{\mathrm{sm}}\), the anti-diagonal class
\[
(\varphi_{\mathcal F},-\varphi_{\mathcal F})
   \in HC^{2n+1}(A_{\mathrm{sm}}\oplus A_{\mathrm{sm}})
\]
is the class determined by the two oppositely oriented boundary components of the
cylinder.

The periodic cyclic cohomology exact sequence of the cylinder extension gives a
boundary map
\[
\partial_{\mathrm{cyl}}:
HP^{2n}(SA_{\mathrm{sm}})
   \longrightarrow
HP^{2n+1}(A_{\mathrm{sm}}\oplus A_{\mathrm{sm}}).
\]
The Chern--Simons transgression class
\[
CS(\varphi_{\mathcal F})\in HC^{2n}(SA_{\mathrm{sm}})
\]
is characterized, with the present orientation convention, by
\[
\partial_{\mathrm{cyl}}
  \bigl(CS(\varphi_{\mathcal F})\bigr)
   =
(\varphi_{\mathcal F},-\varphi_{\mathcal F}).
\]
Equivalently, \(CS(\varphi_{\mathcal F})\) is the class corresponding to
\(\varphi_{\mathcal F}\) under the cyclic-cohomology suspension isomorphism
\[
HP^{2n}(SA_{\mathrm{sm}})
   \cong
HP^{2n+1}(A_{\mathrm{sm}}),
\]
with sign fixed by the orientation of the cylinder. This is the noncommutative
Chern--Simons transgression of the odd cyclic cocycle; see
\cite{QuillenIHES1988,QuillenCS1990,Perrot2000,Cuntz1997,Connes94}.

This parity shift is not an internal operation on \(A_{\mathrm{sm}}\). There is no
canonical map
\[
HC^{2n+1}(A_{\mathrm{sm}})
   \longrightarrow
HC^{2n}(A_{\mathrm{sm}})
\]
producing an even cyclic class on the same algebra. The even class appears only after
passing to the suspension and the cylinder extension.

On the \(K\)-theory side, Bott periodicity gives
\[
\beta:K_1(A)\stackrel{\sim}{\longrightarrow}K_0(SA);
\]
see \cite{Blackadar98,Connes94,WeggeOlsen1993,RordamLarsenLaustsen2000}. Hence
\[
[u_{\mathcal F}]\in K_1(A)
\]
determines
\[
\beta([u_{\mathcal F}])\in K_0(SA).
\]
This class is obtained from the odd class \([u_{\mathcal F}]\), not from the even tgm
class
\[
[e_{\mathcal F}]\in K_0(A).
\]
Indeed, one suspension gives
\[
K_0(A)\cong K_1(SA),
\]
not a natural map
\[
K_0(A)\longrightarrow K_0(SA).
\]
Thus Route B still requires the genuine odd class supplied by Route A.

\PropHead{prop:routeA-routeB}
Assume that \((M,\mathcal F)\) has odd codimension
\[
q=2n+1.
\]
Let
\[
A=C_r^*(M,\mathcal F),
\qquad
A_{\mathrm{sm}}\subset A
\]
be a smooth model carrying the transverse fundamental cyclic cocycle.  Let
\[
\varphi_{\mathcal F}\in HC^{2n+1}(A_{\mathrm{sm}})
\]
be the transverse fundamental cyclic cocycle, and let
\[
[u_{\mathcal F}]\in K_1(A)
\]
be the odd class specified by Route A and represented in the chosen smooth
model.  Let
\[
CS(\varphi_{\mathcal F})\in HC^{2n}(SA_{\mathrm{sm}})
\]
denote the Chern--Simons transgression class determined in periodic cyclic
cohomology by
\[
\partial_{\mathrm{cyl}}\bigl(CS(\varphi_{\mathcal F})\bigr)
=
(\varphi_{\mathcal F},-\varphi_{\mathcal F}),
\]
and let
\[
\beta:K_1(A)\stackrel{\sim}{\longrightarrow}K_0(SA)
\]
be the Bott isomorphism with the compatible suspension convention.  Then, as an
identity of periodic Chern--Connes pairings,
\[
\langle \varphi_{\mathcal F},[u_{\mathcal F}]\rangle
=
\langle CS(\varphi_{\mathcal F}),\beta([u_{\mathcal F}])\rangle .
\]
Consequently Route B is a suspension/transgression reformulation of Route A; it
does not produce a new odd class and does not replace the Route A definition.

\ProofHead
The smooth cylinder extension
\[
0\longrightarrow SA_{\mathrm{sm}}
  \longrightarrow A_{I,\mathrm{sm}}
  \xrightarrow{(\mathrm{ev}_0,\mathrm{ev}_1)}
  A_{\mathrm{sm}}\oplus A_{\mathrm{sm}}
  \longrightarrow 0
\]
gives the suspension exact sequence in periodic cyclic cohomology; see
\cite{CuntzQuillen,CuntzQuillenInvent1997,Brodzki2013}.  With the boundary
orientation convention fixed above, the anti-diagonal endpoint class
\[
(\varphi_{\mathcal F},-\varphi_{\mathcal F})
\]
is the image of the suspended class
\[
CS(\varphi_{\mathcal F})\in HP^{2n}(SA_{\mathrm{sm}}).
\]
This is the cyclic Chern--Simons transgression convention used here; see
\cite{QuillenIHES1988,QuillenCS1990,Perrot2000,Cuntz1997}.

Bott periodicity gives
\[
\beta:K_1(A)\stackrel{\sim}{\longrightarrow}K_0(SA),
\]
see \cite{Blackadar98,Connes94,WeggeOlsen1993,RordamLarsenLaustsen2000}.  The
Chern--Connes character is natural for cyclic excision, suspension and Bott
periodicity.  Therefore, for every smooth representative of
\([u]\in K_1(A)\) on which the odd pairing with
\(\varphi\in HC^{2n+1}(A_{\mathrm{sm}})\) is defined, one has
\[
\langle \varphi,[u]\rangle
=
\langle CS(\varphi),\beta([u])\rangle ,
\]
with the sign fixed by
\[
\partial_{\mathrm{cyl}}\bigl(CS(\varphi)\bigr)
=
(\varphi,-\varphi).
\]
Applying this standard compatibility to
\[
\varphi=\varphi_{\mathcal F},
\qquad
[u]=[u_{\mathcal F}],
\]
gives the claimed equality.  Thus Route B is exactly the periodic-cyclic
suspension form of the Route A pairing. \qed

\RemHead{rem:routeA-vs-routeB}
Route B is retained because it records the functorial form of the odd invariant. Route A
is the definition:
\[
Z(\mathcal F)
   =
\langle \varphi_{\mathcal F},[u_{\mathcal F}]\rangle .
\]
Route B expresses the same number as
\[
Z(\mathcal F)
   =
\langle CS(\varphi_{\mathcal F}),\beta([u_{\mathcal F}])\rangle .
\]
Thus the odd pairing is placed in the even--even framework
\[
HC^{2n}(SA_{\mathrm{sm}})\times K_0(SA).
\]

This reformulation is useful for comparison with extension, mapping-cone and index
theoretic constructions, where even cyclic classes paired with \(K_0\)-classes are often
the natural language. It also exhibits the odd NCFI as a secondary transgression-type
pairing: \(CS(\varphi_{\mathcal F})\) is the cyclic analogue of a Chern--Simons class
obtained from the cylinder.

The two formulas therefore define the same invariant:
\[
\text{Route A:}\qquad
Z(\mathcal F)
=
\langle \varphi_{\mathcal F},[u_{\mathcal F}]\rangle ,
\]
and
\[
\text{Route B:}\qquad
Z(\mathcal F)
=
\langle CS(\varphi_{\mathcal F}),\beta([u_{\mathcal F}])\rangle .
\]
Route A is the constructive definition; Route B is its suspension, transgression and
functorial form.

\subsection{Algebraic conventions}

Throughout the paper, \((M,\mathcal F)\) denotes a smooth foliated manifold with
holonomy groupoid
\[
G=G_{\mathcal F}.
\]
We fix the algebraic conventions used below: the smooth convolution algebra, the
reduced foliation \(C^*\)-algebra, the graph-norm Banach algebra in the tgm
construction, the smooth models used for cyclic cohomology, and the auxiliary
algebras used for parity correction.

\medskip

\noindent
\textbf{Groupoid and smooth convolution algebra.}
The basic groupoid attached to the foliation is its holonomy groupoid \(G\); see
\cite{Connes85,ConnesSurveyFoliations,Connes94,Renault80,ADR2000}.  We write
\[
A_{\mathrm{sm}}(M,\mathcal F)
   :=
C_c^\infty(G,\Omega^{1/2})
\]
for the convolution \(*\)-algebra of compactly supported smooth half-densities on
\(G\).  When the foliation is fixed, we also write
\[
A_{\mathrm{sm}}:=A_{\mathrm{sm}}(M,\mathcal F).
\]
Connes' transverse fundamental cyclic cocycle is represented on this smooth
convolution algebra, or on a smooth algebra transported from it through a chosen
Morita-equivalent groupoid or crossed-product model; see
\cite{Connes80,ConnesTFClass,Connes94,GBVF01}.

\medskip

\noindent
\textbf{Full and reduced foliation \(C^*\)-algebras.}
The full and reduced \(C^*\)-completions of the convolution algebra are denoted by
\[
C^*(G),
\qquad
C_r^*(G),
\]
respectively.  In general they need not coincide; they do coincide for amenable
groupoids; see \cite{ADR2000,Renault80}.  In this article the \(K\)-theoretic side is
always taken in the reduced foliation algebra:
\[
A:=C_r^*(M,\mathcal F)=C_r^*(G).
\]
Thus \(A_{\mathrm{sm}}\) carries the cyclic cocycles, while \(A\) carries the
\(K\)-theory classes.

\medskip

\noindent
\textbf{The graph-norm Banach algebra in the tgm construction.}
The transverse geometric module recalled in Subsection~1.1 starts from
\[
E_{\mathrm{sm}}
   =
C_c^\infty\bigl(G,\Omega^{1/2}\otimes r^*(t_{\mathbb C})\bigr)
\]
and completes it to a Hilbert \(A\)-module \(E\).  The holonomy-induced left action is
generally only closable and is twisted by the modular correction
\[
\Delta(\gamma)=h(\gamma)^{-t}h(\gamma)^{-1}.
\]
Hence one does not usually obtain a bounded \(*\)-representation
\[
A\longrightarrow \operatorname{End}_A(E).
\]
Instead, as in \cite{zois2000}, one passes to the graph norm of the closed left action.
This gives a dense Banach \(*\)-subalgebra
\[
B\subset A.
\]
When it is useful to avoid conflict with other algebras later denoted by \(B\), we write
this graph-norm algebra as
\[
B_{\mathcal F}\subset A.
\]

The graph-norm completion of the smooth module is denoted by
\[
E_B.
\]
This \(B\)-module object, not the Hilbert \(A\)-module \(E\) alone, is used in the
\(K\)-theoretic construction of \cite{zois2000}.  Under the hypotheses of that
construction, \(E_B\) determines a stable projective-module class
\[
[e_B]\in K_0(B).
\]
The inclusion \(B\hookrightarrow A\) induces
\[
K_0(B)\longrightarrow K_0(A),
\]
and the image of \([e_B]\) is denoted by
\[
[e_{\mathcal F}]\in K_0(A).
\]
This is the even \(K\)-class used in the even-codimensional NCFI.  The distinction
between \(B\), \(E_B\), \(E\), and \(A\) remains important even when the final pairing is
written using a class in \(K_0(A)\); see \cite{Connes94,zois2000}.

\medskip

\noindent
\textbf{Where the cyclic cocycle lives.}
The transverse fundamental cyclic cocycle
\[
\varphi_{\mathcal F}
\]
is not, in general, first defined on the \(C^*\)-completion \(A\).  It is represented on
the smooth convolution algebra
\[
A_{\mathrm{sm}}(M,\mathcal F)=C_c^\infty(G,\Omega^{1/2}),
\]
or on a smooth algebra associated with a Morita-equivalent groupoid or crossed-product
model.  Thus the Connes pairing is read in the standard way: the \(K\)-theory class lies
in the reduced foliation \(C^*\)-algebra, while the cyclic cocycle is represented on a
smooth dense algebra on which the pairing formula is defined.

In odd codimension the same convention applies.  There is no general internal
parity-changing operation on cyclic cohomology of the same algebra sending
\[
\varphi_{\mathcal F}\in HC^{2n+1}(A_{\mathrm{sm}})
\]
to an even class in
\[
HC^{2n}(A_{\mathrm{sm}}).
\]
The cyclic-side parity correction of Route B is obtained by passing to the suspension and
the cylinder extension.  At the smooth level we use
\[
A_{I,\mathrm{sm}}
   :=
C^\infty([0,1],A_{\mathrm{sm}})
\]
and
\[
SA_{\mathrm{sm}}
   :=
\{f\in C^\infty([0,1],A_{\mathrm{sm}}):f(0)=f(1)=0\}.
\]
Thus Route B produces a transgressed class
\[
CS(\varphi_{\mathcal F})\in HC^{2n}(SA_{\mathrm{sm}}),
\]
with the periodic-cyclic-cohomology convention of Remark~\ref{rem:periodic-convention}.

\medskip

\noindent
\textbf{Smooth cores and \(K\)-theory.}
In the concrete examples below, we use standard smooth cores.  For the product
fibrational examples these are of the form
\[
C^\infty(B)\widehat{\otimes}\mathcal K^\infty,
\]
where \(B\) is the base and \(\mathcal K^\infty\) denotes smoothing compact operators on
the fibre.  For a nontrivial fibration, the corresponding object is the algebra of
fibrewise smoothing kernels.  For the irrational rotation algebra, the smooth core is
\[
A_\theta^\infty,
\]
the Fréchet algebra of smooth vectors for the gauge action.  These are the smooth
algebras on which the cyclic cocycles, derivations and explicit representatives are
written.

When a \(K\)-theory computation is transported from a smooth core to a \(C^*\)-algebra,
we use the standard fact that the smooth core in question is stable under holomorphic
functional calculus in the corresponding \(C^*\)-algebra, or we choose smooth
representatives of the relevant \(K\)-classes.  This is the convention behind the repeated
use of the same symbols
\[
[e_{\mathcal F}],
\qquad
[u_{\mathcal F}],
\qquad
\varphi_{\mathcal F}
\]
at the smooth and \(C^*\)-levels.  For smooth crossed products and related smooth
cores, see \cite{PhillipsSchweitzer1994}.

This is not a general assertion that every smooth convolution algebra of every foliation
is spectrally invariant in its \(C^*\)-completion.  It is the convention for the standard
smooth models used in the computations of this article.

\medskip

\noindent
\textbf{Transversals, reduced groupoids and crossed-product models.}
If \(T\subset M\) is a complete transversal, then the reduced groupoid
\[
G_T:=G|_T
\]
is an \'{e}tale Lie groupoid.  The groupoids \(G_T\) and \(G\) are equivalent, and the
corresponding reduced groupoid \(C^*\)-algebras are strongly Morita equivalent:
\[
C_r^*(G_T)\sim_M C_r^*(G)=A.
\]
For groupoid equivalence and Morita equivalence, see
\cite{ConnesSurveyFoliations,Connes94,Green1978,BGR77,Renault80,RaeburnWilliams1998,Williams2007}.

In many examples the reduced groupoid over a complete transversal is a transformation
groupoid, giving a crossed-product model
\[
C_r^*(G_T)\cong C_0(T)\rtimes_r\Gamma
\]
for a discrete group \(\Gamma\) acting on \(T\).  At the \(C^*\)-level, Morita equivalence
identifies the \(K\)-theory groups, so the classes
\[
[e_{\mathcal F}]\in K_0(A),
\qquad
[u_{\mathcal F}]\in K_1(A)
\]
may be computed in whichever Morita-equivalent model is most convenient.

At the smooth level, one works with the corresponding smooth convolution algebra or
smooth crossed product carrying the transported cyclic cocycle.  We use the same
symbols for transported cyclic cocycles and \(K\)-classes whenever this causes no
ambiguity.

\medskip

\noindent
\textbf{Auxiliary algebras used for parity corrections.}
In odd codimension we also use auxiliary algebras such as
\[
SA=C_0((0,1))\otimes A
\]
and the smooth suspension \(SA_{\mathrm{sm}}\) defined above.  These algebras implement
Bott periodicity, Chern--Simons transgression, excision, and the even--even
reformulation of odd pairings.  They do not correspond to new foliations; they are
auxiliary algebraic settings for rewriting the same invariant.

\RemHead{rem:section1-convention}
The standing convention in the article is as follows:
\begin{itemize}
\item the foliation \(C^*\)-algebra is the reduced algebra
\[
A=C_r^*(M,\mathcal F);
\]
\item the smooth convolution algebra is
\[
A_{\mathrm{sm}}(M,\mathcal F)=C_c^\infty(G,\Omega^{1/2});
\]
\item the even class \([e_{\mathcal F}]\) is the image in \(K_0(A)\) of the
\(K_0(B)\)-class produced by the graph-norm tgm construction;
\item in odd codimension, the odd class
\[
[u_{\mathcal F}]\in K_1(A)
\]
is used only after it has been specified by one of the odd-favourable mechanisms of
Subsection~1.2;
\item the transverse fundamental cyclic cocycle is represented on
\(A_{\mathrm{sm}}(M,\mathcal F)\), or on an equivalent smooth model obtained by
complete-transversal reduction, Morita equivalence, or crossed-product identification;
\item the cyclic parity correction of Route B is implemented by the cylinder extension
and the suspended algebra \(SA_{\mathrm{sm}}\), not by an internal odd-to-even map on
\(HC^\bullet(A_{\mathrm{sm}})\);
\item when passing to complete transversals, crossed-product models, suspensions or
other auxiliary algebras, we use the same notation for transported classes whenever this
causes no ambiguity.
\end{itemize}
With these conventions fixed, all later formulas are pairings between the appropriate
smooth cyclic-cohomology representatives and the corresponding \(K\)-theory classes on
the reduced foliation \(C^*\)-algebra or on a canonically associated equivalent model.

\RemHead{rem:section1-gvclass}
The NCFI, in its direct even-codimensional form, pairs an even transverse fundamental
cyclic cocycle with a \(K_0\)-class.  This contrasts with classical secondary invariants
such as the Godbillon--Vey class, which arise in odd-degree de Rham cohomology.  The
two constructions probe different geometric data: the NCFI uses noncommutative
transverse \(K\)-theoretic information, while the Godbillon--Vey class is a classical
differential-topological secondary class; see \cite{Connes94,ConnesTFClass,CandelConlonI2000}.

\section{Fibrational foliations in even codimension}
\label{sec:fibrational-even}


Every fibre bundle can be seen as a foliation of the total space where the fibres are the leaves. The converse of course is not true. Fibre bundles constitute the simplest and most straightforward examples of foliations, hence we begin with fibrational foliations examples in even codimension where no parity fixing is needed.  These examples are
important for two reasons.  First, they provide a control case in which the
noncommutative construction reduces to an ordinary characteristic-number
calculation on a smooth base manifold.  Second, they show that the NCFI need
not detect noncommutativity alone: even in Morita-commutative situations it may
detect nontrivial transverse characteristic classes.

The basic situation is a smooth locally trivial fibration
\[
F\hookrightarrow P\xrightarrow{\pi}B,
\]
where \(P\), \(F\), and \(B\) are compact smooth manifolds, the fibres \(F\) are
connected, and \(B\) is closed and oriented.  The associated vertical foliation
\(\mathcal F^{\mathrm v}\) has leaves the fibres of \(\pi\).  The references used in
this section are
\cite{Connes94,ConnesSurveyFoliations,Renault80,ADR2000,ConnesTFClass,
ConnesSkandalisLongitudinal,GBVF01,Milnor,botttu,BGR77,Green1978,
RaeburnWilliams1998}.

\subsection{The fibrational structure theorem}

\PropHead{prop:fibrational-structure}

Let
\[
F\hookrightarrow P\xrightarrow{\pi}B
\]
be a smooth locally trivial fibration of compact manifolds with connected compact
fibres.  Assume that \(B\) is closed and oriented of even dimension
\[
\dim B=q=2n.
\]
Let \(\mathcal F^{\mathrm v}\) be the vertical foliation of \(P\), namely the leaves are the fibres.  Then
\[
\mathcal F^{\mathrm v}=\ker(d\pi)\subset TP,
\qquad
t:=TP/\mathcal F^{\mathrm v}\cong \pi^*TB.
\]
The holonomy groupoid of \(\mathcal F^{\mathrm v}\) is the fibrewise pair groupoid
\[
G_{\mathcal F^{\mathrm v}}\cong P\times_B P.
\]
Consequently
\[
A^{\mathrm v}:=C_r^*(G_{\mathcal F^{\mathrm v}})
\]
is strongly Morita equivalent to \(C(B)\).  Under this Morita equivalence, Connes'
transverse fundamental cyclic cocycle is transported to the ordinary fundamental
homology class of \(B\), and the tgm class is transported to
\[
[TB\otimes\mathbb C]\in K^0(B).
\]
Therefore
\[
Z(\mathcal F^{\mathrm v})
=
\left\langle
\operatorname{ch}(TB\otimes\mathbb C),[B]
\right\rangle
=
\int_B
\operatorname{ch}(TB\otimes\mathbb C)_{[q]},
\]
where \(\operatorname{ch}(TB\otimes\mathbb C)_{[q]}\) denotes the degree-\(q\)
component of the Chern character.

\ProofHead
Since
\[
\mathcal F^{\mathrm v}=\ker(d\pi),
\]
the surjection \(d\pi:TP\to\pi^*TB\) induces
\[
TP/\mathcal F^{\mathrm v}\cong\pi^*TB.
\]
Thus
\[
t_{\mathbb C}\cong\pi^*(TB\otimes\mathbb C).
\]

The holonomy of a fibration foliation is trivial along each fibre, so the holonomy
groupoid is
\[
P\times_B P
=
\{(p_1,p_2)\in P\times P:\pi(p_1)=\pi(p_2)\}.
\]
Its convolution algebra is the algebra of fibrewise smoothing kernels, and its reduced
\(C^*\)-algebra is
\[
C_r^*(P\times_B P)
\cong
\mathcal K_{C(B)}\bigl(L^2(P/B;\Omega^{1/2})\bigr).
\]
This algebra is strongly Morita equivalent to \(C(B)\); see
\cite{Connes94,Renault80,BGR77,Green1978,RaeburnWilliams1998}.

Under this Morita equivalence, the transverse fundamental cyclic cocycle becomes the
ordinary fundamental class of the oriented base \(B\); this is the simple-foliation case
of Connes' transverse fundamental class
\cite{ConnesTFClass,Connes94,GBVF01}.  Since
\[
t_{\mathbb C}\cong \pi^*(TB\otimes\mathbb C),
\]
the tgm class is transported to
\[
[TB\otimes\mathbb C]\in K^0(B).
\]
The Connes pairing therefore becomes the ordinary Chern-character pairing
\[
Z(\mathcal F^{\mathrm v})
=
\left\langle
\operatorname{ch}(TB\otimes\mathbb C),[B]
\right\rangle
=
\int_B
\operatorname{ch}(TB\otimes\mathbb C)_{[q]}.
\]
This proves the proposition. \qed

The following corollary records two immediate consequences.

\CorHead{cor:codim2-fibre-foliations}

(i) With the hypotheses of Proposition~\ref{prop:fibrational-structure}, the fibrational
NCFI is a characteristic number of the base.  Since \(TB\otimes\mathbb C\) is the
complexification of a real vector bundle, the components
\[
\operatorname{ch}_{2j+1}(TB\otimes\mathbb C)
\]
vanish rationally.  Hence
\[
Z(\mathcal F^{\mathrm v})=0
\]
whenever
\[
\dim B\equiv 2\pmod 4.
\]
In particular fibrational foliations (vertical foliations) over closed oriented surfaces have vanishing NCFI (since clearly $2=2(mod 4)$).

(ii) Moreover when the dimension of the base space is equal to  \(4\),
\[
Z(\mathcal F^{\mathrm v})
=
\int_B p_1(TB).
\]
In particular 
the fibrational foliation (vertical foliation) on
\[
\mathbb{CP}^2\times S^1\longrightarrow \mathbb{CP}^2
\]
has value \(3\), as it will be computed explicitly in
Subsection~\ref{subsec:cp2-nonvanishing}.

\ProofHead
For the complexification of a real vector bundle, the Chern roots occur in opposite
pairs.  Hence the odd Chern-character components
\(\operatorname{ch}_{2j+1}\) vanish.  If
\[
\dim B\equiv 2\pmod 4,
\]
then the top-degree component of
\[
\operatorname{ch}(TB\otimes\mathbb C)
\]
is one of these odd components, and its integral over \(B\) is zero.

If \(\dim B=4\), the degree-four component is
\[
\operatorname{ch}(TB\otimes\mathbb C)_{[4]}=p_1(TB),
\]
so
\[
Z(\mathcal F^{\mathrm v})
=
\int_B p_1(TB).
\]
For these characteristic-class identities, see \cite{Milnor,botttu}. \qed

\subsection{The control example \texorpdfstring{\(S^2\times S^1\)}{S2 x S1}}
\label{sec:reference-example-S2S1}
Despite the fact that we have the general results above, it is instructive to see the details of a particular example:\\
Let
\[
P:=S^2\times S^1,
\qquad
\pi:P\to S^2,
\qquad
\pi(x,s)=x.
\]
Let \(\mathcal F^{\mathrm v}\) be the vertical foliation, whose leaves are
\[
\{x\}\times S^1.
\]
This is Proposition~\ref{prop:fibrational-structure} with base \(B=S^2\).  Thus
\[
\mathcal F^{\mathrm v}=\ker(d\pi),
\qquad
t:=TP/\mathcal F^{\mathrm v}\cong \pi^*TS^2.
\]
The holonomy groupoid is
\[
G^{\mathrm v}
=
P\times_{S^2}P
=
\{((x,s_1),(x,s_2)):\ x\in S^2,\ s_1,s_2\in S^1\}.
\]
Writing an arrow as \((x,s_1,s_2)\),
\[
s(x,s_1,s_2)=(x,s_2),
\qquad
r(x,s_1,s_2)=(x,s_1),
\]
and
\[
(x,s_1,s_2)\circ(x,s_2,s_3)=(x,s_1,s_3).
\]

The smooth convolution algebra is
\[
C_c^\infty(G^{\mathrm v},\Omega^{1/2})
\cong
C^\infty(S^2)\,\widehat\otimes\,\mathcal K^\infty,
\]
where \(\mathcal K^\infty\) denotes smoothing compact operators on \(S^1\).  The
reduced foliation algebra is
\[
A^{\mathrm v}:=C_r^*(G^{\mathrm v})
\cong
C(S^2)\otimes\mathcal K(L^2(S^1)).
\]
Thus \(A^{\mathrm v}\) is strongly Morita equivalent to \(C(S^2)\).  We use the
smooth dense algebra
\[
A^{\mathrm v,\infty}
:=
C^\infty(S^2)\,\widehat\otimes\,\mathcal K^\infty.
\]

Since
\[
\operatorname{codim}(\mathcal F^{\mathrm v})=2,
\]
the transverse fundamental cyclic cocycle has degree \(2\).  With the normalization
used here, a representative is
\begin{equation}\label{eq:tfcc-S2S1}
\varphi_{\mathcal F^{\mathrm v}}(a_0,a_1,a_2)
:=
\frac{1}{2\pi i}
\int_{S^2}
\tau_{\mathcal K}\!\bigl(a_0\,da_1\wedge da_2\bigr),
\qquad
a_0,a_1,a_2\in A^{\mathrm v,\infty},
\end{equation}
where \(d\) is the de Rham differential on \(S^2\), and \(\tau_{\mathcal K}\)
is the trace on \(\mathcal K^\infty\).  This is the transported form of Connes'
TFCC; in this simple-foliation case it is the fundamental cyclic cocycle of \(S^2\).

The tgm is built from
\[
C_c^\infty\!\bigl(
G^{\mathrm v},\Omega^{1/2}\otimes r^*(t_{\mathbb C})
\bigr),
\qquad
t_{\mathbb C}
\cong
\pi^*(TS^2\otimes\mathbb C).
\]
Hence, under
\[
K_0(A^{\mathrm v})\cong K^0(S^2),
\]
the tgm class corresponds to
\[
[TS^2\otimes\mathbb C]\in K^0(S^2).
\]

It is essential to distinguish this bundle from the holomorphic tangent line bundle.
The tgm uses
\[
TS^2\otimes\mathbb C,
\]
not
\[
T^{1,0}S^2.
\]
Although
\[
\int_{S^2}c_1(T^{1,0}S^2)=2,
\]
that is not the class entering the tgm.

Since \(TS^2\) is an oriented real \(2\)-plane bundle, it is the underlying real bundle
of a complex line bundle \(L\).  Equivalently,
\[
TS^2\otimes\mathbb C
\cong
T^{1,0}S^2\oplus T^{0,1}S^2
\cong
\mathcal O(2)\oplus\mathcal O(-2).
\]
Therefore
\[
c_1(TS^2\otimes\mathbb C)=0.
\]
In \(K^0(S^2)\), the reduced parts of
\[
[\mathcal O(2)]
\quad\text{and}\quad
[\mathcal O(-2)]
\]
cancel, so
\[
[TS^2\otimes\mathbb C]=2[1].
\]
Thus the tgm class is not zero in \(K_0\); it is the trivial rank-two class.  What
vanishes is its degree-\(2\) Chern-character component.

Choose a rank-one projection
\[
e_0\in\mathcal K^\infty.
\]
A stabilized projection representing the tgm class is
\[
e_{\mathcal F^{\mathrm v}}
:=
\mathbf 1_2\otimes e_0
\in
M_2(A^{\mathrm v,\infty}).
\]
It is constant in the \(S^2\)-direction, so
\[
d(e_{\mathcal F^{\mathrm v}})=0.
\]
Using \eqref{eq:tfcc-S2S1},
\[
\varphi_{\mathcal F^{\mathrm v}}
(e_{\mathcal F^{\mathrm v}},
 e_{\mathcal F^{\mathrm v}},
 e_{\mathcal F^{\mathrm v}})
=
0.
\]
Therefore
\[
Z(\mathcal F^{\mathrm v})
=
\big\langle
\varphi_{\mathcal F^{\mathrm v}},
[e_{\mathcal F^{\mathrm v}}]
\big\rangle
=
0.
\]
This agrees with Corollary~\ref{cor:codim2-fibre-foliations}, since
\[
\dim S^2=2\equiv 2\pmod 4.
\]

\RemHead{rem:s2-classical-comparison}
The half-density formalism is part of the analytic realization of the algebra and the
module, but it does not change the transported topological class
\[
t_{\mathbb C}=t\otimes\mathbb C.
\]
In this example the tgm class is the trivial rank-two class in \(K_0\).  The TFCC
detects only its degree-\(2\) component, which is zero.

The classical number
\[
2=\chi(S^2)
\]
appears if one pairs the transverse fundamental class with the holomorphic tangent line
bundle
\[
T^{1,0}S^2\cong\mathcal O(2):
\]
then
\[
\int_{S^2}c_1(T^{1,0}S^2)=2.
\]
This is not the NCFI from the tgm construction.  The Bott connection does not enter
the vanishing; the vanishing follows from
\[
t_{\mathbb C}\cong L\oplus\overline L.
\]
Thus \(S^2\times S^1\) is a control case.  The codimension-\(4\) example
\(\mathbb{CP}^2\times S^1\) below shows that the same fibrational mechanism can
give a nonzero value.

\RemHead{rem:s2-horizontal-flat}
For a flat \(U(1)\)-principal bundle over \(S^2\), there is no new horizontal case of
the kind considered in \cite{zois2000}.  Principal \(U(1)\)-bundles over \(S^2\) are
classified by
\[
c_1\in H^2(S^2;\mathbb Z)\cong\mathbb Z.
\]
A flat connection has zero curvature, hence zero real Chern class; because
\(H^2(S^2;\mathbb Z)\) has no torsion, the bundle must be trivial.  Since \(S^2\)
is simply connected, the trivial bundle has, up to gauge equivalence, only the trivial
flat holonomy.  The associated horizontal foliation is therefore equivalent to the
product foliation by the \(S^2\)-slices in
\[
S^2\times S^1.
\]
It is a consistency check, not a new family of NCFI values.

\subsection{Broader codimension-2 vanishing criteria}
\label{subsec:codim2-vanishing}

The fibrational surface case is one instance of a broader codimension-\(2\) vanishing
mechanism.  The following criteria isolate two situations in which the canonical tgm
class has zero degree-\(2\) Chern-character content.

\PropHeadPrime{prop:codim2-broader-vanishing}
Let \((M,\mathcal F)\) be a compact transversally oriented foliation of codimension
\(2\), let \(G=G_{\mathcal F}\), and set
\[
A_{\mathrm{sm}}:=C_c^\infty(G,\Omega^{1/2}).
\]
Let
\[
E_{\mathrm{sm}}
:=
C_c^\infty\!\bigl(G,\Omega^{1/2}\otimes r^*(t_{\mathbb C})\bigr)
\]
be the smooth tgm, let
\[
[e_{\mathcal F}]\in K_0\!\bigl(C_r^*(M,\mathcal F)\bigr)
\]
be the corresponding tgm class, and let
\[
\varphi_{\mathcal F}\in HC^2(A_{\mathrm{sm}})
\]
be Connes' TFCC.

Assume that, after passage to a complete transversal \(T\) and transport to
\(G_T:=G|_T\), one of the following conditions holds.

\[
\text{\emph{(i) Freeness condition.}}
\]
The transported smooth tgm is free of rank \(2\) as a right module over
\[
A_{T,\mathrm{sm}}:=C_c^\infty(G_T,\Omega^{1/2}):
\]
\[
E_{T,\mathrm{sm}}
:=
C_c^\infty\!\bigl(G_T,\Omega^{1/2}\otimes r^*(t_{\mathbb C}|_T)\bigr)
\cong
A_{T,\mathrm{sm}}^{\oplus 2}.
\]

\[
\text{\emph{(ii) Commutative surface condition.}}
\]
The reduced foliation algebra is strongly Morita equivalent to
\[
C(\Sigma)\otimes\mathcal K,
\]
where \(\Sigma\) is a closed oriented surface.  Under this Morita equivalence,
\(\varphi_{\mathcal F}\) is transported to the fundamental cyclic \(2\)-cocycle
of \(\Sigma\), and \([e_{\mathcal F}]\) is transported to the class of
\[
\xi\otimes\mathbb C,
\]
where \(\xi\to\Sigma\) is a real oriented rank-\(2\) vector bundle.

Then
\[
Z(\mathcal F)
=
\langle \varphi_{\mathcal F},[e_{\mathcal F}]\rangle
=
0.
\]

\ProofHead
The pairing may be computed on the reduced complete-transversal model; see
\cite{Connes94,ConnesTFClass}.

If \emph{(i)} holds, the transported tgm class is the free rank-two class
\[
[e_{\mathcal F}]=2[1].
\]
It is represented by a constant rank-two projection, so its differential in the smooth
model is zero.  Hence the degree-\(2\) Connes--Chern character component vanishes,
and
\[
\langle \varphi_{\mathcal F},[e_{\mathcal F}]\rangle=0.
\]

If \emph{(ii)} holds, then \(\xi\) is the underlying real bundle of a complex line bundle
\(L\), and
\[
\xi\otimes\mathbb C
\cong
L\oplus\overline L.
\]
Thus
\[
c_1(\xi\otimes\mathbb C)
=
c_1(L)+c_1(\overline L)
=
0.
\]
Equivalently, the odd Chern classes of the complexification of a real bundle are
\(2\)-torsion; see \cite[Chs.~14--15]{Milnor}.  The degree-\(2\) component of
the Chern character of the transported tgm class therefore vanishes.  Since the
codimension-\(2\) TFCC pairing on the commutative surface model is governed by this
component \cite{Connes94,ConnesTFClass},
\[
\langle \varphi_{\mathcal F},[e_{\mathcal F}]\rangle=0.
\]
This proves the proposition. \qed

We shall also use the following surface-current version.

\PropHeadDoublePrime{prop:codim2-surface-current-vanishing}
Let \((M,\mathcal F)\) be a compact transversally oriented foliation of codimension
\(2\), let \(G=G_{\mathcal F}\), and let
\[
A_{\mathrm{sm}}:=C_c^\infty(G,\Omega^{1/2}),
\qquad
E_{\mathrm{sm}}
:=
C_c^\infty\!\bigl(G,\Omega^{1/2}\otimes r^*(t_{\mathbb C})\bigr).
\]
Let
\[
[e_{\mathcal F}]\in K_0\!\bigl(C_r^*(M,\mathcal F)\bigr)
\]
be the transported tgm class, and let
\[
\varphi_{\mathcal F}\in HC^2(A_{\mathrm{sm}})
\]
be the TFCC.

Assume that, after passage to a complete transversal \(T\) and transport to the
reduced étale model \(G_T:=G|_T\), the following hold.

\[
\text{\emph{(i)}}
\]
The transversal \(T\) is a closed oriented surface, and the transported TFCC is
represented on
\[
A_{T,\mathrm{sm}}:=C_c^\infty(G_T,\Omega^{1/2})
\]
by the group-degree-zero cyclic \(2\)-cocycle defined by the fundamental current of
\(T\).

\[
\text{\emph{(ii)}}
\]
The transported full tgm class is represented by
\[
E_{T,\mathrm{sm}}
\cong
C_c^\infty\!\bigl(G_T,\Omega^{1/2}\otimes r^*(\xi\otimes\mathbb C)\bigr),
\]
where \(\xi\to T\) is a \(G_T\)-equivariant real oriented rank-\(2\) vector bundle.

Then
\[
Z(\mathcal F)
=
\langle \varphi_{\mathcal F},[e_{\mathcal F}]\rangle
=
0.
\]

\ProofHead
The pairing may be computed on the reduced étale model \(G_T\); see
\cite{Connes94,ConnesTFClass}.  By hypothesis, the transported TFCC is the
group-degree-zero cyclic \(2\)-cocycle defined by the fundamental current of \(T\).
Thus only the ordinary degree-\(2\) component of the Connes--Chern character of the
transported tgm class contributes; see \cite{CrainicThesis}.

Since \(\xi\) is a real oriented rank-\(2\) bundle, it is the underlying real bundle of
a complex line bundle \(L\).  Hence
\[
\xi\otimes\mathbb C
\cong
L\oplus\overline L,
\qquad
c_1(\xi\otimes\mathbb C)=0.
\]
The degree-\(2\) component of the Connes--Chern character of the transported tgm
class vanishes, and therefore
\[
\langle \varphi_{\mathcal F},[e_{\mathcal F}]\rangle=0.
\]
This proves the proposition. \qed

\RemHead{rem:codim2-status}
The preceding results are not a universal vanishing theorem for all codimension-\(2\)
foliations.  They identify broad classes in which the canonical tgm class built from
\[
t_{\mathbb C}=t\otimes\mathbb C
\]
has zero pairing with the TFCC.  This includes the fibrational surface case of
Corollary~\ref{cor:codim2-fibre-foliations}, the freeness and commutative surface
cases of Proposition~\ref{prop:codim2-broader-vanishing}, and the surface-current
case of Proposition~\ref{prop:codim2-surface-current-vanishing}.

The reason is that a cyclic \(2\)-cocycle sees the degree-\(2\) component of the
Connes--Chern character, while for the complexification of a real oriented rank-\(2\)
bundle,
\[
c_1(t_{\mathbb C})=0
\]
in the complex-valued pairing; see \cite{Connes94,ConnesTFClass,Milnor}.  Hence the
standard tgm class in codimension \(2\) has a strong built-in tendency to vanish under
the TFCC pairing.  A codimension-\(2\), non-fibrational example with nonzero NCFI for
the canonical tgm class would have to evade both the freeness mechanism and the
surface-current mechanism above.

\subsection{A nonvanishing codimension-4 fibrational example}
\label{subsec:cp2-nonvanishing}

We now give the simplest nonzero fibrational example.  Codimension-\(2\) fibrations
over surfaces vanish because the degree-\(2\) Chern-character component of the
complexified real transverse bundle vanishes.  In codimension \(4\), the relevant
component is the first Pontryagin class.

Consider
\[
P_4:=\mathbb{CP}^2\times S^1,
\qquad
\pi:P_4\to\mathbb{CP}^2,
\qquad
\pi(x,s)=x,
\]
and let \(\mathcal F^{\mathrm v}_4\) be the vertical foliation by the fibres
\[
\{x\}\times S^1.
\]
The holonomy groupoid is
\[
G^{\mathrm v}_4
\cong
P_4\times_{\mathbb{CP}^2}P_4,
\]
and
\[
A^{\mathrm v}_4
:=
C_r^*(G^{\mathrm v}_4)
\cong
C(\mathbb{CP}^2)\otimes\mathcal K(L^2(S^1)),
\]
with smooth dense core
\[
A^{\mathrm v,\infty}_4
\cong
C^\infty(\mathbb{CP}^2)\,\widehat\otimes\,\mathcal K^\infty.
\]
Thus
\[
A^{\mathrm v}_4\sim_M C(\mathbb{CP}^2).
\]
For the groupoid and Morita-equivalence background, see
\cite{Connes94,Renault80,ConnesSkandalisLongitudinal,GBVF01}.

Since
\[
\operatorname{codim}(\mathcal F^{\mathrm v}_4)=4,
\]
the TFCC has degree \(4\).  With the normalization used here, a representative on
\(A^{\mathrm v,\infty}_4\) is
\[
\varphi_{\mathcal F^{\mathrm v}_4}(a_0,a_1,a_2,a_3,a_4)
:=
\frac{1}{(2\pi i)^2\,2!}
\int_{\mathbb{CP}^2}
\tau_{\mathcal K}\!
\bigl(
a_0\,da_1\wedge da_2\wedge da_3\wedge da_4
\bigr),
\]
for
\[
a_0,\ldots,a_4\in A^{\mathrm v,\infty}_4.
\]
Here \(d\) is the de Rham differential on \(\mathbb{CP}^2\), and
\(\tau_{\mathcal K}\) is the trace on \(\mathcal K^\infty\).  This is the
transported form of Connes' TFCC; in this simple-foliation case it is the ordinary
fundamental cyclic cocycle of \(\mathbb{CP}^2\)
\cite{Connes94,ConnesTFClass,GBVF01}.

The transverse bundle is
\[
t\cong\pi^*T\mathbb{CP}^2,
\qquad
t_{\mathbb C}
\cong
\pi^*(T\mathbb{CP}^2\otimes\mathbb C).
\]
Under
\[
K_0(A^{\mathrm v}_4)
\cong
K^0(\mathbb{CP}^2),
\]
the tgm class corresponds to
\[
[T\mathbb{CP}^2\otimes\mathbb C].
\]

Let
\[
h\in H^2(\mathbb{CP}^2;\mathbb Z)
\]
be the positive generator.  The standard Chern-class identities are
\[
c_1(T^{1,0}\mathbb{CP}^2)=3h,
\qquad
c_2(T^{1,0}\mathbb{CP}^2)=3h^2;
\]
see \cite{Milnor,botttu}.  Since
\[
T\mathbb{CP}^2\otimes\mathbb C
\cong
T^{1,0}\mathbb{CP}^2
\oplus
T^{0,1}\mathbb{CP}^2,
\]
we have
\[
\operatorname{ch}(T\mathbb{CP}^2\otimes\mathbb C)_{[4]}
=
\operatorname{ch}_2(T^{1,0}\mathbb{CP}^2)
+
\operatorname{ch}_2(T^{0,1}\mathbb{CP}^2).
\]
For a complex bundle \(E\),
\[
\operatorname{ch}_2(E)
=
\frac12\bigl(c_1(E)^2-2c_2(E)\bigr).
\]
The conjugate summand contributes the same degree-four term, so
\[
\operatorname{ch}(T\mathbb{CP}^2\otimes\mathbb C)_{[4]}
=
c_1(T^{1,0}\mathbb{CP}^2)^2
-
2c_2(T^{1,0}\mathbb{CP}^2).
\]
Thus
\[
\operatorname{ch}(T\mathbb{CP}^2\otimes\mathbb C)_{[4]}
=
(3h)^2-2(3h^2)
=
3h^2.
\]
Equivalently,
\[
\operatorname{ch}(T\mathbb{CP}^2\otimes\mathbb C)_{[4]}
=
p_1(T\mathbb{CP}^2).
\]
Since
\[
\int_{\mathbb{CP}^2}h^2=1,
\]
we obtain
\[
\int_{\mathbb{CP}^2}
\operatorname{ch}(T\mathbb{CP}^2\otimes\mathbb C)_{[4]}
=
\int_{\mathbb{CP}^2}
p_1(T\mathbb{CP}^2)
=
3.
\]

Therefore
\begin{equation}\label{eq:zois-cp2s1-value}
Z(\mathcal F^{\mathrm v}_4)
=
\big\langle
\varphi_{\mathcal F^{\mathrm v}_4},
[e_{\mathcal F^{\mathrm v}_4}]
\big\rangle
=
3.
\end{equation}
Under the Morita equivalence with \(C(\mathbb{CP}^2)\), the TFCC is the de Rham
fundamental class and the tgm class is
\[
[T\mathbb{CP}^2\otimes\mathbb C],
\]
so the Connes pairing is
\[
\left\langle
\operatorname{ch}(T\mathbb{CP}^2\otimes\mathbb C),
[\mathbb{CP}^2]
\right\rangle.
\]
Thus, in this example,
\[
Z(\mathcal F^{\mathrm v}_4)
=
\int_{\mathbb{CP}^2}p_1(T\mathbb{CP}^2)
=
3.
\]
For the characteristic-class identities used above, see \cite{Milnor,botttu}.

\RemHead{rem:cp2-nonvanishing}
This is the codimension-\(4\) analogue of the vertical fibration examples above, but
it is nonzero.  It shows that in higher even codimension the NCFI can detect an
ordinary transverse characteristic number; here it detects
\[
p_1[\mathbb{CP}^2]=3.
\]

The example is still fibrational and Morita-commutative, so it is not a genuinely
noncommutative nonzero example.  It is nevertheless a useful control case: the
codimension-\(2\) vanishing is not a defect of the definition, but a consequence of
the characteristic class detected in that degree.  A non-fibrational nonzero
even-codimensional example would have to produce a nonzero top-degree component of
the Connes--Chern character of the transported tgm class in a singular leaf-space
setting.

\section{The Kronecker foliation: geometry and holonomy groupoids}
\label{sec:kronecker-geometry-holonomy}

We now turn to the Kronecker foliation on the two-torus, the main
odd-codimensional family considered in the paper.  Since its codimension is one, the
transverse fundamental cyclic cocycle has odd degree, and the parity-fixing framework
of Section~1 is required.

This section fixes the geometric and groupoid data.  The \(C^*\)-algebraic models,
smooth cores, \(K\)-theory and cyclic cohomology are treated in the following sections.
For geometric and dynamical background on linear flows on tori and Kronecker
foliations, see
\cite{CandelConlonI2000,KatokHasselblatt95,GottschalkHedlund,Walters82,Weyl1916,HardyWright}.
For holonomy groupoids, complete transversals and groupoid equivalence, see
\cite{Renault80,ConnesSurveyFoliations,Connes94}.  The rotation-algebra and
\(K\)-theoretic references \cite{Rieffel81,PimsnerVoiculescu80} enter later in the
crossed-product model.

\subsection{Geometric setup}
\label{subsec:kronecker-geometric-setup}

Let
\[
T^2=\mathbb R^2/\mathbb Z^2
\]
be the two-torus with coordinates \((x,y)\) modulo integers.  Fix
\[
\theta\in\mathbb R,
\]
and consider the constant vector field
\[
X_\theta=\partial_x+\theta\,\partial_y.
\]
Its flow is
\[
\varphi_t(x,y)=(x+t,\;y+\theta t)\qquad (\mathrm{mod}\ 1).
\]
The Kronecker foliation of slope \(\theta\), denoted
\[
\mathcal F_\theta,
\]
is the one-dimensional foliation with tangent distribution
\[
\mathcal F_\theta=\mathbb R X_\theta\subset TT^2.
\]
Thus the leaves are the orbits of \(\varphi_t\).  The transverse bundle is
\[
t_\theta:=TT^2/\mathcal F_\theta.
\]
Since \(X_\theta\) is nowhere vanishing,
\[
\operatorname{codim}(\mathcal F_\theta)=1.
\]

\subsection{Leaf structure}
\label{subsec:kronecker-leaf-structure}

The leaf structure depends on whether \(\theta\) is rational or irrational.

If
\[
\theta=\frac pq\in\mathbb Q,
\qquad
\gcd(p,q)=1,\quad q>0,
\]
then all leaves are closed circles.  The map
\[
f_{p/q}:T^2\longrightarrow S^1,
\qquad
f_{p/q}(x,y)=qy-px\pmod 1
\]
is a smooth submersion, and
\[
df_{p/q}(X_{p/q})
=
q\cdot \frac pq-p
=
0.
\]
Thus the fibres of \(f_{p/q}\) are tangent to \(X_{p/q}\).  Since the vector
\((-p,q)\) is primitive in \(\mathbb Z^2\), these fibres are connected circles.
Hence the rational Kronecker foliation is the circle fibration
\[
f_{p/q}:T^2\longrightarrow S^1.
\]

If
\[
\theta\notin\mathbb Q,
\]
then every leaf is diffeomorphic to \(\mathbb R\).  Indeed, if
\[
t\longmapsto (x_0+t,\;y_0+\theta t)\pmod 1
\]
failed to be injective, then for some \(t_1\neq t_2\),
\[
(t_2-t_1,\;\theta(t_2-t_1))\in\mathbb Z^2.
\]
Writing \(n=t_2-t_1\in\mathbb Z\), with \(n\neq 0\), gives
\[
\theta=\frac{\theta n}{n}\in\mathbb Q,
\]
a contradiction.  Hence the orbit map is injective and the leaf is \(\mathbb R\).

For irrational \(\theta\), every leaf is also dense in \(T^2\).  This is the classical
minimality of irrational linear flows on the torus; see
\cite{GottschalkHedlund,Walters82,Weyl1916}.  Thus the irrational Kronecker foliation
has no compact leaves, and its classical leaf space is singular.

\subsection{The conormal form and the transversal coordinate}
\label{subsec:kronecker-conormal-transversal-coordinate}

A global \(1\)-form annihilating \(X_\theta\) is
\[
\omega_\theta
=
dy-\theta\,dx.
\]
Indeed,
\[
\omega_\theta(X_\theta)
=
(dy-\theta\,dx)(\partial_x+\theta\partial_y)
=
\theta-\theta
=
0.
\]
Hence
\[
\ker(\omega_\theta)=\mathcal F_\theta.
\]
The form \(\omega_\theta\) is closed and nowhere vanishing, so
\(\mathcal F_\theta\) is transversally oriented and the transverse line bundle
\(t_\theta\) is trivial.

On the universal cover \(\mathbb R^2\), define
\[
z:=y-\theta x.
\]
Then
\[
dz=dy-\theta dx=\omega_\theta,
\qquad
X_\theta(z)=0.
\]
Thus \(z\) is a first integral of the lifted flow.  When
\(\theta\notin\mathbb Q\), the function \(z\) does not descend to a globally defined
\(S^1\)-valued function on \(T^2\), because deck transformations change it by
integer combinations of \(1\) and \(\theta\).  Its differential
\[
dz=\omega_\theta
\]
is \(\mathbb Z^2\)-invariant and descends to the global conormal form on \(T^2\).

This distinction is used below.  The form
\[
\omega_\theta=dy-\theta dx
\]
is the global conormal form of the foliation, while \(z\) is the coordinate on a chosen
complete transversal.  On the vertical transversal
\[
T:=\{x=0\}\subset T^2,
\]
one has
\[
z=y.
\]
Thus, on the complete transversal, the conormal form is represented by the transversal
differential \(dz\).  This is the geometric source of the transversal derivation cocycle
used later in the TFCC computation.

The transverse orientation convention used in the paper is the one determined by
\[
\omega_\theta=dy-\theta dx.
\]
Reversing this convention reverses the sign of all odd transverse pairings.

\subsection{General holonomy-groupoid convention}
\label{subsec:kronecker-holonomy-general}

For a foliated manifold \((M,\mathcal F)\), the holonomy groupoid
\(G_\mathcal F\) has object space \(M\).  Its arrows are leafwise paths modulo
holonomy equivalence, with composition induced by concatenation.  Thus
\(G_\mathcal F\) records both leaf equivalence and holonomy transport; see
\cite{Renault80,ConnesSurveyFoliations,Connes94}.

For the Kronecker foliation, the form of the holonomy groupoid depends on the slope.
In the irrational case, the leaves are simply connected and have no nontrivial holonomy,
so the holonomy groupoid agrees with the transformation groupoid of the
\(\mathbb R\)-flow.  In the rational case, the flow is periodic; the transformation
groupoid remembers the time parameter and is larger than the actual holonomy groupoid.
The actual holonomy groupoid is then the fibrewise pair groupoid of the circle
fibration.

\subsection{Irrational slope: full and reduced holonomy groupoids}
\label{subsec:kronecker-irrational-groupoids}

Assume that
\[
\theta\notin\mathbb Q.
\]
Every leaf of \(\mathcal F_\theta\) is diffeomorphic to \(\mathbb R\), hence is simply
connected, and there is no nontrivial holonomy.  The holonomy groupoid is therefore
the transformation groupoid of the \(\mathbb R\)-action generated by \(X_\theta\):
\[
G_{\mathrm{full},\theta}
=
T^2\rtimes_\tau\mathbb R,
\]
where
\[
\tau_t(x,y)=(x+t,\;y+\theta t)\qquad(\mathrm{mod}\ 1).
\]
An arrow is a pair
\[
((x,y),t),
\]
with source and range
\[
s((x,y),t)=(x,y),
\qquad
r((x,y),t)=(x+t,\;y+\theta t).
\]
The inverse is
\[
((x,y),t)^{-1}
=
((x+t,\;y+\theta t),-t),
\]
and the composition is
\[
((x+t,\;y+\theta t),s)\circ((x,y),t)
=
((x,y),t+s).
\]

A convenient complete transversal is the vertical circle
\[
T=\{x=0\}\cong S^1.
\]
We use the coordinate
\[
z\in\mathbb R/\mathbb Z
\]
on \(T\).  Since \(x=0\) on \(T\), this agrees there with the coordinate \(y\).  Starting
from \((0,z)\), flowing for time \(1\) gives
\[
(1,z+\theta),
\]
which is identified in \(T^2\) with
\[
(0,z+\theta).
\]
Thus the first-return map on \(T\) is the irrational rotation
\[
\alpha:T\longrightarrow T,
\qquad
\alpha(z)=z+\theta\pmod 1.
\]
More generally,
\[
\alpha^n(z)=z+n\theta\pmod 1,
\qquad n\in\mathbb Z.
\]

The reduced holonomy groupoid over the complete transversal is
\[
G_{\mathrm{red},\theta}
=
T\rtimes_\alpha\mathbb Z
\cong
S^1\rtimes_\alpha\mathbb Z.
\]
With the transformation-groupoid convention used here, an arrow is a pair
\[
(z,n),
\]
with source and range
\[
s(z,n)=z,
\qquad
r(z,n)=z+n\theta.
\]
The inverse is
\[
(z,n)^{-1}=(z+n\theta,-n),
\]
and the composition is
\[
(z+n\theta,m)\circ(z,n)=(z,n+m).
\]

The groupoids
\[
T^2\rtimes_\tau\mathbb R
\quad\text{and}\quad
S^1\rtimes_\alpha\mathbb Z
\]
are equivalent because \(T\) is a complete transversal.  Consequently their reduced
groupoid \(C^*\)-algebras are strongly Morita equivalent; see
\cite{Renault80,ConnesSurveyFoliations,Connes94,BGR77,RaeburnWilliams1998}.
The reduced groupoid is the model used later for the irrational rotation algebra.

The coordinate convention is important.  The unitary \(V\) in the reduced
crossed-product model is the function
\[
V(z)=e^{2\pi iz}
\]
on the transversal.  Thus \(V\) is built from the complete-transversal coordinate
\(z\), not from an independent global coordinate \(y\) on the original torus.  Since
\[
dz=dy-\theta dx,
\]
the transversal differential \(dz\) is the reduced representative of the global conormal
form \(\omega_\theta\).  This point is used later to identify the transported transverse
fundamental cyclic cocycle.

\subsection{Rational slope: fibration and holonomy groupoid}
\label{subsec:kronecker-rational-groupoids}

Assume now that
\[
\theta=\frac pq\in\mathbb Q,
\qquad
\gcd(p,q)=1,\quad q>0.
\]
Then \(\mathcal F_{p/q}\) is the circle fibration
\[
f_{p/q}:T^2\longrightarrow S^1,
\qquad
f_{p/q}(x,y)=qy-px\pmod 1.
\]
The leaves are the fibres of \(f_{p/q}\), and each leaf has trivial holonomy.  Therefore
the actual holonomy groupoid is the fibrewise pair groupoid
\[
G_{\mathrm{hol}}^{p/q}
=
T^2\times_{S^1}T^2
=
\{(m_1,m_2)\in T^2\times T^2:
  f_{p/q}(m_1)=f_{p/q}(m_2)\}.
\]

This must be distinguished from the transformation groupoid
\[
T^2\rtimes_\tau\mathbb R
\]
of the periodic flow.  In the rational case the \(\mathbb R\)-flow has period
\(q\), because flowing for time \(q\) gives
\[
(x,y)\longmapsto (x+q,\;y+p)=(x,y)
\]
on \(T^2\).  Thus the transformation groupoid retains time isotropy
\(q\mathbb Z\).  The holonomy groupoid identifies leafwise paths with the same
holonomy germ; since the fibration has trivial holonomy, it is the fibrewise pair
groupoid above.

The same distinction appears after restricting to the vertical complete transversal
\[
T=\{x=0\}\cong S^1.
\]
The first-return map is
\[
\alpha(z)=z+\frac pq\pmod 1.
\]
If one records the full return-map transformation groupoid, one obtains
\[
T\rtimes_\alpha\mathbb Z.
\]
Its orbits have \(q\) points, and every point has isotropy \(q\mathbb Z\).  The subgroup
\(q\mathbb Z\) acts trivially on the transversal and represents return paths with trivial
holonomy germ.  Hence the reduced holonomy groupoid on the transversal is the finite
equivalence-relation groupoid
\[
T\rtimes_\alpha(\mathbb Z/q\mathbb Z).
\]
Since the action of \(\mathbb Z/q\mathbb Z\) on \(T\) by rotation through \(p/q\) is
free, this finite groupoid is Morita equivalent to the unit groupoid of the quotient
circle
\[
T/(\mathbb Z/q\mathbb Z)\cong S^1.
\]
This agrees with the fibrewise pair-groupoid description of the fibration
\[
T^2\longrightarrow S^1.
\]

Thus, for rational slope, the actual foliation algebra is the algebra associated with a
circle fibration and is Morita equivalent to \(C(S^1)\).  The full return-map crossed
product by \(\mathbb Z\) is a different object: it retains the redundant isotropy
\(q\mathbb Z\).  This distinction is used later when separating the actual rational
foliation algebra from the rational rotation algebra.

\subsection{Summary of the Kronecker groupoids}
\label{subsec:kronecker-groupoid-summary}

The groupoid picture is therefore as follows.

If
\[
\theta\notin\mathbb Q,
\]
then the leaves are noncompact copies of \(\mathbb R\), every leaf is dense, and the
holonomy groupoid is
\[
G_{\mathrm{full},\theta}=T^2\rtimes_\tau\mathbb R.
\]
Reduction to the complete transversal \(T=\{x=0\}\) gives the equivalent étale
groupoid
\[
G_{\mathrm{red},\theta}=S^1\rtimes_\alpha\mathbb Z,
\qquad
\alpha(z)=z+\theta.
\]
The corresponding \(C^*\)-algebraic model is the irrational rotation algebra, treated in
the next section.

If
\[
\theta=\frac pq\in\mathbb Q,
\qquad
\gcd(p,q)=1,
\]
then the foliation is the circle fibration
\[
f_{p/q}:T^2\longrightarrow S^1,
\qquad
f_{p/q}(x,y)=qy-px\pmod 1,
\]
and the actual holonomy groupoid is the fibrewise pair groupoid
\[
T^2\times_{S^1}T^2.
\]
On a complete transversal the actual reduced holonomy groupoid is
\[
S^1\rtimes_\alpha(\mathbb Z/q\mathbb Z),
\]
not the full return-map groupoid
\[
S^1\rtimes_\alpha\mathbb Z.
\]
The latter retains the trivial isotropy \(q\mathbb Z\) and is not the foliation
holonomy groupoid in the rational case.

Finally, the global conormal form of the Kronecker foliation is
\[
\omega_\theta=dy-\theta dx.
\]
On the complete-transversal model this form is represented by
\[
dz,
\]
where \(z\) is the coordinate on \(T=\{x=0\}\).  This is the key input for identifying
the transported transverse fundamental cyclic cocycle in the irrational case.

\section{\texorpdfstring{Groupoid \(C^*\)-algebras and crossed products}
{Groupoid C*-algebras and crossed products}}
\label{sec:kronecker-cstar-models}
\label{sec:groupoid-cstar-crossed-products}

We record the analytic models associated with the Kronecker foliation.  These
include the global flow groupoid, the complete-transversal groupoid, their
convolution algebras and \(C^*\)-completions, the crossed-product presentation
of the rotation algebra, the twisted group-algebra model, and the smooth core
used for cyclic cohomology.  The irrational case is treated in detail because it
is the noncommutative case used in the odd-codimensional NCFI computation.  The
rational case is separated at the end: its actual holonomy-groupoid algebra is
Morita equivalent to \(C(S^1)\), not to the rational rotation algebra.

The main references for this section are
\cite{Renault80,Paterson1999,ADR2000,RaeburnWilliams1998,Williams2007,
Pedersen1979,Davidson96,Vassout,khalkhali,Connes94,Rieffel81,
PimsnerVoiculescu80,Green1978,pr}.

\subsection{Irrational slope: the global flow groupoid}
\label{subsec:irrational-global-flow-cstar}

Assume that
\[
\theta\notin\mathbb Q.
\]
As in Section~\ref{sec:kronecker-geometry-holonomy}, the holonomy groupoid of the
irrational Kronecker foliation agrees with the transformation groupoid of the
\(\mathbb R\)-flow generated by
\[
X_\theta=\partial_x+\theta\,\partial_y.
\]
Thus
\[
G_{\mathrm{full},\theta}
=
T^2\rtimes_\tau\mathbb R,
\]
where
\[
\tau_t(x,y)=(x+t,\;y+\theta t)\qquad(\mathrm{mod}\ 1).
\]
An arrow is a pair
\[
((x,y),t),
\]
with source and range
\[
s((x,y),t)=(x,y),
\qquad
r((x,y),t)=(x+t,\;y+\theta t).
\]
The unit at \((x,y)\) is
\[
u(x,y)=((x,y),0),
\]
the inverse is
\[
((x,y),t)^{-1}=((x+t,y+\theta t),-t),
\]
and composition is
\[
((x+t,y+\theta t),s)\circ((x,y),t)=((x,y),t+s).
\]

A Haar system is given by Lebesgue measure \(dt\) on \(\mathbb R\); see
\cite{Renault80,Paterson1999}.  The convolution and involution on
\(C_c^\infty(T^2\times\mathbb R)\) are
\[
(f*g)((x,y),t)
=
\int_{\mathbb R}
f((x,y),s)\,
g((x+s,y+\theta s),t-s)\,ds,
\]
and
\[
f^*((x,y),t)
=
\overline{f((x+t,y+\theta t),-t)}.
\]

Since \(\mathbb R\) is amenable, the full and reduced crossed-product completions
coincide; see \cite{ADR2000}.  Hence
\[
C^*(T^2\rtimes_\tau\mathbb R)
\cong
C_r^*(T^2\rtimes_\tau\mathbb R)
\cong
C(T^2)\rtimes_\tau\mathbb R.
\]
The induced action on functions is
\[
(\tau_t f)(x,y)=f(x-t,\;y-\theta t),
\qquad
f\in C(T^2).
\]

\subsection{Irrational slope: reduction to a complete transversal}
\label{subsec:irrational-transversal-crossed-product}

Use the vertical complete transversal
\[
T=\{x=0\}\cong S^1.
\]
We write
\[
S^1_z
\]
when this circle is viewed with its complete-transversal coordinate \(z\).  On \(T\),
\(z\) agrees with the coordinate \(y\), but its geometric meaning comes from the lifted
first integral
\[
z=y-\theta x.
\]
Thus
\[
dz=dy-\theta dx
\]
is the global conormal form represented on the complete transversal.

The reduced groupoid over \(T\) is
\[
G_{\mathrm{red},\theta}
=
S^1_z\rtimes_\alpha\mathbb Z,
\]
where
\[
\alpha_n(z)=z+n\theta\qquad(\mathrm{mod}\ 1),
\qquad
n\in\mathbb Z.
\]
With the convention fixed in Section~\ref{sec:kronecker-geometry-holonomy}, an arrow
\((z,n)\) goes from \(z\) to \(z+n\theta\).  Hence
\[
s(z,n)=z,
\qquad
r(z,n)=z+n\theta,
\]
\[
(z,n)^{-1}=(z+n\theta,-n),
\]
and
\[
(z+n\theta,m)\circ(z,n)=(z,n+m).
\]

The Haar system is counting measure on \(\mathbb Z\); see
\cite{Renault80,Paterson1999}.  The convolution and involution on
\(C_c^\infty(S^1_z\times\mathbb Z)\) are
\[
(f*g)(z,n)
=
\sum_{m\in\mathbb Z}
f(z,m)\,
g(z+m\theta,n-m),
\]
and
\[
f^*(z,n)
=
\overline{f(z+n\theta,-n)}.
\]

Since \(\mathbb Z\) is amenable, full and reduced completions coincide; see
\cite{ADR2000}.  Therefore
\[
C^*(S^1_z\rtimes_\alpha\mathbb Z)
\cong
C_r^*(S^1_z\rtimes_\alpha\mathbb Z)
\cong
C(S^1_z)\rtimes_\alpha\mathbb Z.
\]
In the crossed product, \(\alpha\) denotes the induced action on functions:
\[
(\alpha_n f)(z)=f(z-n\theta),
\qquad
f\in C(S^1_z).
\]

\subsection{Crossed products and the irrational rotation algebra}
\label{subsec:rotation-algebra-generators}

A covariant pair for
\[
(C(S^1_z),\alpha,\mathbb Z)
\]
on a Hilbert space \(\mathcal H\) consists of a \(*\)-representation
\[
\pi:C(S^1_z)\longrightarrow \mathcal B(\mathcal H)
\]
and a unitary representation
\[
U:\mathbb Z\longrightarrow \mathcal U(\mathcal H)
\]
such that
\[
U(1)\,\pi(f)\,U(1)^*
=
\pi(\alpha_1 f),
\qquad
f\in C(S^1_z).
\]
For crossed products and covariant representations, see
\cite{Williams2007,Pedersen1979,Davidson96}.

The integrated form is
\[
(\pi\rtimes U)\left(\sum_{n\in\mathbb Z} f_n u^n\right)
=
\sum_{n\in\mathbb Z}\pi(f_n)\,U(n),
\]
first on the dense algebra of finite sums and then by completion.

A canonical covariant pair is the Koopman pair on \(L^2(S^1_z)\):
\[
(\pi(f)\xi)(z)=f(z)\xi(z),
\qquad
(U(1)\xi)(z)=\xi(z-\theta).
\]
Let
\[
V(z)=e^{2\pi iz}\in C(S^1_z),
\]
and write
\[
U:=U(1)
\]
for the implementing unitary.  Covariance gives
\[
UVU^*=\alpha_1(V).
\]
Since
\[
\alpha_1(V)(z)=V(z-\theta)=e^{-2\pi i\theta}V(z),
\]
we obtain
\[
UVU^*=e^{-2\pi i\theta}V,
\]
or equivalently
\[
VU=e^{2\pi i\theta}UV.
\]

Thus
\[
A_\theta:=C(S^1_z)\rtimes_\alpha\mathbb Z
\]
is the universal \(C^*\)-algebra generated by unitaries \(U,V\) satisfying
\[
VU=e^{2\pi i\theta}UV.
\]
This is the irrational rotation algebra, or noncommutative two-torus; see
\cite{Rieffel81,Connes94,khalkhali}.

The generator \(U\) is the crossed-product implementing unitary and represents the
one-step return map of the complete transversal.  The generator
\[
V(z)=e^{2\pi iz}
\]
is the coordinate unitary on the transversal.  Thus \(V\) is tied to the coordinate
\(z\), and
\[
dz=dy-\theta dx
\]
is the reduced representative of the global conormal form.

\subsection{The twisted group-algebra model}
\label{subsec:twisted-group-algebra}

Equivalently, \(A_\theta\) may be written as a twisted group \(C^*\)-algebra.  Use the
normalized \(2\)-cocycle
\[
\sigma_\theta((m,n),(m',n'))
=
\exp(2\pi i\,\theta\,n m')
\]
on \(\mathbb Z^2\).  Let \(W_{m,n}\) be the canonical unitaries satisfying
\[
W_{m,n}W_{m',n'}
=
\sigma_\theta((m,n),(m',n'))\,W_{m+m',n+n'}.
\]
With
\[
U=W_{1,0},
\qquad
V=W_{0,1},
\]
one obtains
\[
VU=e^{2\pi i\theta}UV.
\]
Hence
\[
A_\theta\cong C^*(\mathbb Z^2,\sigma_\theta).
\]
For twisted group \(C^*\)-algebras and their relation to crossed products, see
\cite{PackerRaeburn1989,RaeburnWilliams1998}.  This model is useful for Fourier and
Schwartz-class calculations.

\subsection{Global and complete-transversal models: Morita equivalence}
\label{subsec:global-reduced-morita}

The global-flow groupoid
\[
T^2\rtimes_\tau\mathbb R
\]
and the complete-transversal groupoid
\[
S^1_z\rtimes_\alpha\mathbb Z
\]
are equivalent groupoids because \(T=\{x=0\}\) is a complete transversal; see
\cite{Renault80,ConnesSurveyFoliations,Connes94}.  Therefore their reduced groupoid
\(C^*\)-algebras are strongly Morita equivalent:
\[
C(T^2)\rtimes_\tau\mathbb R
\sim_M
C(S^1_z)\rtimes_\alpha\mathbb Z
\cong
A_\theta.
\]
For the Morita-equivalence background, see
\cite{Green1978,BGR77,RaeburnWilliams1998,Williams2007}.

The specified groupoid equivalence gives an imprimitivity bimodule and hence identifies
the \(K\)-theory groups of the two \(C^*\)-models.  At the smooth level, cyclic cocycles
are transported through the corresponding smooth groupoid or crossed-product models.
Thus explicit computations may be carried out in
\[
A_\theta=C(S^1_z)\rtimes_\alpha\mathbb Z,
\]
and in its smooth core
\[
A_\theta^\infty.
\]

\subsection{Smooth core, derivations, and canonical trace}
\label{subsec:smooth-core-derivations-trace}

The smooth noncommutative torus
\[
A_\theta^\infty
\]
is the Fréchet \(*\)-algebra
\[
A_\theta^\infty
=
\left\{
\sum_{m,n\in\mathbb Z} c_{m,n}U^mV^n
:
(c_{m,n})\in\mathcal S(\mathbb Z^2)
\right\}.
\]
Equivalently, it is the algebra of smooth vectors for the gauge action \(\beta\) of
\(\mathbb T^2\) given by
\[
\beta_{(s,t)}(U)=e^{2\pi is}U,
\qquad
\beta_{(s,t)}(V)=e^{2\pi it}V.
\]
For the smooth noncommutative torus and its differential calculus, see
\cite{Connes94,Rieffel1990,khalkhali,Vassout}.

The associated derivations are
\[
\delta_1(U)=2\pi i\,U,
\qquad
\delta_1(V)=0,
\]
and
\[
\delta_2(U)=0,
\qquad
\delta_2(V)=2\pi i\,V.
\]
They extend to all of \(A_\theta^\infty\) by the Leibniz rule and continuity.  The
derivation \(\delta_1\) differentiates in the crossed-product return direction, while
\(\delta_2\) differentiates in the transversal coordinate direction \(z\).

The canonical trace
\[
\tau_0:A_\theta\longrightarrow\mathbb C
\]
is defined on the smooth core by Fourier-coefficient extraction:
\[
\tau_0\left(\sum_{m,n}c_{m,n}U^mV^n\right)=c_{0,0}.
\]
It is faithful on \(A_\theta\); see \cite{Connes94,Rieffel81,Rieffel1990,khalkhali}.

The algebra \(A_\theta^\infty\) is not another \(C^*\)-algebraic model Morita
equivalent to the previous ones.  It is a dense smooth Fréchet subalgebra of
\(A_\theta\), used for cyclic cohomology, derivations, smooth \(K\)-theory
representatives, and explicit Connes-pairing computations.

\subsection{Rational slope}
\label{subsec:rational-slope-cstar}

Assume now that
\[
\theta=\frac pq\in\mathbb Q,
\qquad
\gcd(p,q)=1,\quad q>0.
\]
As explained in Section~\ref{sec:kronecker-geometry-holonomy}, the Kronecker
foliation is the circle fibration
\[
S^1\hookrightarrow T^2\longrightarrow S^1.
\]
Its actual holonomy groupoid is the fibrewise pair groupoid
\[
T^2\times_{S^1}T^2.
\]
Therefore the reduced foliation \(C^*\)-algebra is
\[
C_r^*(\mathcal F_{p/q})
\cong
C(S^1)\otimes\mathcal K(L^2(S^1)).
\]
Equivalently, because this fibration has a global section, one may use a complete
transversal on which the reduced holonomy groupoid is the unit groupoid of \(S^1\).
The Morita-equivalent commutative model is then
\[
C(S^1).
\]
This is the foliation algebraic model for the rational Kronecker foliation.

If one instead restricts to the vertical transversal used in the irrational case, the
actual reduced holonomy groupoid is the finite transformation groupoid
\[
S^1\rtimes_\alpha(\mathbb Z/q\mathbb Z),
\]
not the full return-map groupoid
\[
S^1\rtimes_\alpha\mathbb Z.
\]
The finite groupoid is Morita equivalent to the quotient circle
\[
S^1/(\mathbb Z/q\mathbb Z)\cong S^1.
\]

By contrast, the rational rotation algebra
\[
A_{p/q}:=C(S^1)\rtimes_\alpha\mathbb Z
\]
retains the redundant isotropy subgroup \(q\mathbb Z\).  It is a continuous-trace
algebra strongly Morita equivalent to \(C(T^2)\); see
\cite{PimsnerVoiculescu80,RieffelStableRank,ElliottEvans1993}.  Thus \(A_{p/q}\) is
a natural algebra attached to the rational rotation action, but it is not the reduced
foliation \(C^*\)-algebra of the rational Kronecker foliation under the
holonomy-groupoid conventions used here.

\subsection{Summary of the analytic models}
\label{subsec:kronecker-cstar-summary}

For irrational slope
\[
\theta\notin\mathbb Q,
\]
the global-flow groupoid
\[
T^2\rtimes_\tau\mathbb R
\]
and the complete-transversal groupoid
\[
S^1_z\rtimes_\alpha\mathbb Z
\]
are equivalent.  Hence
\[
C(T^2)\rtimes_\tau\mathbb R
\sim_M
C(S^1_z)\rtimes_\alpha\mathbb Z.
\]
The second algebra is the irrational rotation algebra
\[
A_\theta,
\]
generated by unitaries \(U,V\) satisfying
\[
VU=e^{2\pi i\theta}UV.
\]
Its smooth core
\[
A_\theta^\infty
\]
is the Fréchet algebra used for cyclic cohomology and differential calculus.

For rational slope
\[
\theta=\frac pq,
\]
the actual foliation algebra is obtained from the circle-fibration holonomy groupoid and
is Morita equivalent to
\[
C(S^1).
\]
The rational rotation algebra
\[
A_{p/q}=C(S^1)\rtimes_\alpha\mathbb Z
\]
is a different continuous-trace algebra, Morita equivalent to \(C(T^2)\), and is not the
actual holonomy-groupoid algebra of the rational Kronecker foliation.

Finally, in the irrational complete-transversal model, the coordinate unitary
\[
V(z)=e^{2\pi iz}
\]
is tied to the transversal coordinate \(z\), with
\[
dz=dy-\theta dx.
\]
This coordinate convention is used in the identification of the transported transverse
fundamental cyclic cocycle.

\section{\texorpdfstring{\(K\)-theory, cyclic cohomology, and the TFCC for the smooth noncommutative \(2\)-torus}
{K-theory, cyclic cohomology, and the TFCC for the smooth noncommutative 2-torus}}
\label{sec:atheta-k-cyclic-tfcc}

We now describe the \(K\)-theory and cyclic cohomology of the irrational rotation
algebra
\[
A_\theta=C(S^1_z)\rtimes_\alpha\mathbb Z,
\]
and identify the transverse fundamental cyclic cocycle transported from the irrational
Kronecker foliation.  In the complete-transversal model the cocycle differentiates in
the transversal coordinate \(z\).  Thus
\[
\varphi_\theta=\psi^{(2)}_1.
\]
This identification is the input for the odd \(K_1\)-pairings used below.

The main references for this section are
\cite{Elliott84,PimsnerVoiculescu80,Loday98,CuntzQuillen,GBVF01,HigsonRoe00,
Blackadar98,Rieffel81,Rieffel88,Rieffel1990,khalkhali,Boca01,zoisktheory18,
Hadfield,Connes94,ConnesTFClass}.

\subsection{Setup and notation}
\label{subsec:atheta-setup}

Throughout this section we work in the irrational case
\[
\theta\notin\mathbb Q.
\]
For formulas involving Rieffel projections and the canonical trace on \(K_0\), we often
use the normalization
\[
0<\theta<1.
\]
For a general irrational slope, these formulas are recovered after replacing \(\theta\) by
its fractional part, since
\[
A_{\theta+n}\cong A_\theta
\]
for all
\[
n\in\mathbb Z.
\]

As in Section~\ref{sec:kronecker-cstar-models}, the complete-transversal model is
\[
A_\theta
\cong
C(S^1_z)\rtimes_\alpha\mathbb Z,
\]
where
\[
\alpha_n(z)=z+n\theta \pmod 1
\]
on the transversal, and the induced action on functions is
\[
(\alpha_n f)(z)=f(z-n\theta).
\]
Let \(U\) denote the implementing unitary of the crossed product:
\[
UfU^*=\alpha_1(f),
\qquad
f\in C(S^1_z).
\]
Let
\[
V(z)=e^{2\pi iz}.
\]
Then
\[
UVU^*
=
\alpha_1(V)
=
e^{-2\pi i\theta}V,
\]
or equivalently
\[
VU=e^{2\pi i\theta}UV.
\]
Thus \(A_\theta\) is the universal \(C^*\)-algebra generated by unitaries \(U,V\)
satisfying this relation; see \cite{Rieffel81,Connes94,khalkhali}.

The smooth noncommutative torus is
\[
A^\infty_\theta
=
\left\{
\sum_{m,n\in\mathbb Z}c_{m,n}U^mV^n:
(c_{m,n})\in\mathcal S(\mathbb Z^2)
\right\}.
\]
It is the smooth dense Fréchet \(*\)-subalgebra of \(A_\theta\) used for cyclic
cohomology, derivations and explicit Connes pairings; see
\cite{Rieffel1990,khalkhali}.

The canonical trace is
\[
\tau_0\left(\sum_{m,n}c_{m,n}U^mV^n\right)=c_{0,0}.
\]
The standard derivations are
\[
\delta_1(U)=2\pi i\,U,
\qquad
\delta_1(V)=0,
\]
and
\[
\delta_2(U)=0,
\qquad
\delta_2(V)=2\pi i\,V.
\]
Thus \(\delta_1\) differentiates in the crossed-product return direction, while
\(\delta_2\) differentiates in the complete-transversal coordinate direction \(z\).

\subsection{Full and reduced models}
\label{subsec:atheta-full-reduced}

For the irrational Kronecker foliation, the full flow groupoid is
\[
T^2\rtimes_\tau\mathbb R,
\]
and the reduced complete-transversal groupoid is
\[
S^1_z\rtimes_\alpha\mathbb Z.
\]
As explained in Sections~\ref{sec:kronecker-geometry-holonomy} and
\ref{sec:kronecker-cstar-models}, these groupoids are equivalent because
\(S^1_z=\{x=0\}\) is a complete transversal.  Therefore their reduced groupoid
\(C^*\)-algebras are strongly Morita equivalent:
\[
C(T^2)\rtimes_\tau\mathbb R
\sim_M
C(S^1_z)\rtimes_\alpha\mathbb Z
\cong
A_\theta.
\]
For the groupoid and Morita-equivalence background, see
\cite{Renault80,ConnesSurveyFoliations,Connes94,Green1978,BGR77,
RaeburnWilliams1998}.

At the \(C^*\)-level, this Morita equivalence identifies \(K\)-theory through the
associated imprimitivity bimodule.  At the smooth level, cyclic cocycles are transported
through the corresponding smooth groupoid or crossed-product models.  Hence the
explicit computation may be carried out in
\[
A^\infty_\theta.
\]

\subsection{\texorpdfstring{The Pimsner--Voiculescu sequence and \(K\)-theory}
{The Pimsner--Voiculescu sequence and K-theory}}
\label{subsec:pv-ktheory-atheta}

Let
\[
B=C(S^1_z),
\qquad
A_\theta=B\rtimes_\alpha\mathbb Z.
\]
The Pimsner--Voiculescu six-term exact sequence for the crossed product is
\[
K_0(B)
\xrightarrow{\mathrm{id}-\alpha_*}
K_0(B)
\longrightarrow
K_0(A_\theta)
\xrightarrow{\partial}
K_1(B)
\xrightarrow{\mathrm{id}-\alpha_*}
K_1(B)
\longrightarrow
K_1(A_\theta)
\xrightarrow{\partial}
K_0(B).
\]
See \cite{PimsnerVoiculescu80,Blackadar98,zoisktheory18}.

Since
\[
B=C(S^1),
\]
one has
\[
K_0(B)\cong\mathbb Z,
\qquad
K_1(B)\cong\mathbb Z.
\]
The generator of \(K_0(B)\) is \([1]\), and the generator of \(K_1(B)\) is the class of
the coordinate unitary
\[
[V].
\]
The rotation \(\alpha\) is homotopic to the identity, so
\[
\alpha_*=\mathrm{id}
\]
on both \(K_0(B)\) and \(K_1(B)\).  Hence
\[
\mathrm{id}-\alpha_*=0,
\]
and the Pimsner--Voiculescu sequence splits into short exact sequences
\[
0\longrightarrow K_0(B)
\longrightarrow K_0(A_\theta)
\longrightarrow K_1(B)
\longrightarrow 0,
\]
and
\[
0\longrightarrow K_1(B)
\longrightarrow K_1(A_\theta)
\longrightarrow K_0(B)
\longrightarrow 0.
\]
Since these are extensions of free abelian groups, they split.  Therefore
\[
K_0(A_\theta)\cong\mathbb Z^2,
\qquad
K_1(A_\theta)\cong\mathbb Z^2.
\]
This computation is independent of the irrational value of \(\theta\); see
\cite{Elliott84,PimsnerVoiculescu80,Blackadar98,Rieffel81}.

For \(K_1(A_\theta)\), we use the standard generators
\[
K_1(A_\theta)\cong \mathbb Z[U]\oplus\mathbb Z[V].
\]
For \(K_0(A_\theta)\), one generator is the unit class
\[
[1].
\]
A second generator may be represented by a Rieffel projection
\[
e_\theta\in M_N(A^\infty_\theta),
\qquad
N\geq 1,
\]
chosen so that
\[
\tau_0(e_\theta)=\theta,
\qquad
\langle \varphi_2,[e_\theta]\rangle=1
\]
for the normalized even cyclic \(2\)-cocycle \(\varphi_2\) introduced below.  Thus
\[
K_0(A_\theta)\cong \mathbb Z[1]\oplus\mathbb Z[e_\theta].
\]
For Rieffel projections and explicit projection formulas, see
\cite{Rieffel81,Rieffel88,Rieffel1990,Boca01}.

\subsection{Periodic cyclic cohomology and standard cyclic cocycles}
\label{subsec:atheta-periodic-cyclic}

For a smooth Fréchet algebra such as \(A^\infty_\theta\), cyclic cohomology is defined
by the \((b,B)\)-bicomplex, and periodic cyclic cohomology is the corresponding
\(2\)-periodic theory.  Periodic cyclic cohomology is the theory which pairs naturally
with topological \(K\)-theory.  See
\cite{Connes85,Connes94,Loday98,CuntzQuillen,GBVF01,khalkhali}.

For the smooth noncommutative \(2\)-torus, Connes computed
\[
HP^\bullet(A^\infty_\theta)
\]
and showed that it is naturally isomorphic, as a graded vector space, to the de Rham
cohomology of the ordinary \(2\)-torus:
\[
HP^\bullet(A^\infty_\theta)
\cong
H^\bullet_{\mathrm{dR}}(T^2;\mathbb C).
\]
In particular,
\[
HP^0(A^\infty_\theta)\cong\mathbb C^2,
\qquad
HP^1(A^\infty_\theta)\cong\mathbb C^2.
\]
See \cite{Connes94,Rieffel1990,GBVF01,khalkhali}.

We record explicit representatives.  The canonical trace
\[
\tau_0
\]
is a cyclic \(0\)-cocycle.  A standard cyclic \(2\)-cocycle is
\[
\varphi_2(a_0,a_1,a_2)
=
\frac{1}{2\pi i}
\tau_0\!\left(
a_0\bigl(
\delta_1(a_1)\delta_2(a_2)
-
\delta_2(a_1)\delta_1(a_2)
\bigr)
\right).
\]
The classes
\[
[\tau_0],
\qquad
[\varphi_2]
\]
form a basis of
\[
HP^0(A^\infty_\theta).
\]

The two standard odd cyclic \(1\)-cocycles are
\[
\psi^{(1)}_1(a_0,a_1)
=
\frac{1}{2\pi i}
\tau_0\bigl(a_0\delta_1(a_1)\bigr),
\]
and
\[
\psi^{(2)}_1(a_0,a_1)
=
\frac{1}{2\pi i}
\tau_0\bigl(a_0\delta_2(a_1)\bigr).
\]
Their classes form a basis of
\[
HP^1(A^\infty_\theta).
\]
This is the standard derivation-cocycle basis for the smooth noncommutative torus; see
\cite{Connes94,Rieffel1990,GBVF01,khalkhali}.

The coordinate convention used below is the complete-transversal one.  The unitary
\[
V(z)=e^{2\pi iz}
\]
is the coordinate unitary on the transversal, while \(U\) is the crossed-product return
unitary.  Therefore
\[
\delta_2
\]
is the derivation in the transversal coordinate \(z\).

\subsection{The transported TFCC for the irrational Kronecker foliation}
\label{subsec:kronecker-tfcc}

The irrational Kronecker foliation has codimension
\[
q=1.
\]
Therefore its transverse fundamental cyclic cocycle has odd degree.  The following
lemma identifies its representative on the smooth complete-transversal model.

\LemHead{lem:kronecker-tfcc-identification}
Let
\[
\theta\notin\mathbb Q,
\]
and let \(\mathcal F_\theta\) be the irrational Kronecker foliation on \(T^2\) generated
by
\[
X_\theta=\partial_x+\theta\,\partial_y.
\]
Let
\[
T=\{x=0\}\cong S^1_z
\]
be the complete transversal with coordinate \(z\), and let
\[
A^\infty_\theta=C^\infty(S^1_z)\rtimes_\alpha\mathbb Z
\]
be the smooth complete-transversal crossed-product model.  Let \(U\) be the
crossed-product implementing unitary and let
\[
V(z)=e^{2\pi iz}.
\]
With the convention
\[
VU=e^{2\pi i\theta}UV,
\]
the transverse fundamental cyclic cocycle of \(\mathcal F_\theta\), transported to
\(A^\infty_\theta\), is represented by the cyclic \(1\)-cocycle
\[
\varphi_\theta(a_0,a_1)
=
\frac{1}{2\pi i}
\tau_0\bigl(a_0\delta_2(a_1)\bigr),
\qquad
a_0,a_1\in A^\infty_\theta.
\]
Equivalently,
\[
\varphi_\theta=\psi^{(2)}_1.
\]
Under the geometric identification
\[
z=y-\theta x
\]
on the universal cover, this cocycle corresponds to the global conormal form
\[
dz=dy-\theta dx.
\]
Changing the transverse orientation changes the sign.

\ProofHead
The global conormal form of the Kronecker foliation is
\[
\omega_\theta=dy-\theta dx,
\]
since
\[
\omega_\theta(X_\theta)=0.
\]
On the universal cover \(\mathbb R^2\), define
\[
z=y-\theta x.
\]
Then
\[
dz=dy-\theta dx=\omega_\theta.
\]
Although \(z\) itself does not descend to a globally defined circle-valued function on
\(T^2\) when \(\theta\notin\mathbb Q\), its differential descends and is precisely the
conormal form.  On the complete transversal
\[
T=\{x=0\},
\]
the coordinate \(z\) agrees with the coordinate \(y\) restricted to \(T\).

By Connes' construction of the transverse fundamental cyclic cocycle for a
complete-transversal crossed-product model, one uses the crossed-product differential
graded algebra
\[
\Omega^\bullet(T)\rtimes_\alpha\mathbb Z.
\]
The differential is the ordinary de Rham differential on the transversal:
\[
d_T\left(\sum_{n\in\mathbb Z}\omega_n U^n\right)
=
\sum_{n\in\mathbb Z}d_T\omega_n\,U^n.
\]
The closed graded trace extracts the group-degree-zero component and integrates it over
\(T\).  Therefore, in codimension one,
\[
\varphi_\theta(a_0,a_1)
=
\frac{1}{2\pi i}
\int_{S^1_z}
(a_0\,d_Ta_1)_0,
\]
where \((\cdot)_0\) denotes the \(\mathbb Z\)-degree-zero component.  This is Connes'
crossed-product formula for the transverse fundamental class; see
\cite{ConnesTFClass,Connes94,CrainicThesis}.

Writing
\[
a=\sum_{n\in\mathbb Z}a_n(z)U^n,
\]
we have
\[
d_Ta
=
\sum_{n\in\mathbb Z}a_n'(z)\,dz\,U^n.
\]
Thus the cocycle differentiates only in the transversal variable \(z\).  In the standard
smooth noncommutative-torus notation, this is the derivation
\[
\delta_2(U)=0,
\qquad
\delta_2(V)=2\pi iV.
\]
Consequently
\[
\varphi_\theta(a_0,a_1)
=
\frac{1}{2\pi i}
\tau_0\bigl(a_0\delta_2(a_1)\bigr)
=
\psi^{(2)}_1(a_0,a_1).
\]
Since
\[
dz=dy-\theta dx,
\]
this complete-transversal cocycle is the transported representative of the global
conormal class of the foliation.  \qed

For rational slope
\[
\theta=\frac pq\in\mathbb Q,
\]
the actual foliation is a circle fibration over \(S^1\), and the actual foliation algebra
is Morita equivalent to \(C(S^1)\), as explained in
Sections~\ref{sec:kronecker-geometry-holonomy} and
\ref{sec:kronecker-cstar-models}.  Its TFCC is the ordinary integration \(1\)-cocycle
on the transverse circle, with the normalization determined by the chosen transverse
coordinate.  This is distinct from the rational rotation algebra \(A_{p/q}\), which is a
continuous-trace crossed product associated with the non-effective return action.

\subsection{Pairings}
\label{subsec:atheta-pairings}

We now record the standard pairings, fixing the normalization of the cyclic cocycles and
\(K\)-theory generators.

For \(K_0(A_\theta)\), with the normalization above,
\[
\langle \tau_0,[1]\rangle=1,
\qquad
\langle \tau_0,[e_\theta]\rangle=\theta,
\]
and
\[
\langle \varphi_2,[1]\rangle=0,
\qquad
\langle \varphi_2,[e_\theta]\rangle=1.
\]
These are the standard trace and Chern-number pairings for the Rieffel projection; see
\cite{Rieffel81,Rieffel88,Rieffel1990,Boca01,Connes94}.

For \(K_1(A_\theta)\), the odd pairings with the standard generators are
\[
\langle\psi^{(1)}_1,[U]\rangle=1,
\qquad
\langle\psi^{(1)}_1,[V]\rangle=0,
\]
and
\[
\langle\psi^{(2)}_1,[U]\rangle=0,
\qquad
\langle\psi^{(2)}_1,[V]\rangle=1.
\]
Indeed,
\[
\langle\psi^{(1)}_1,[U]\rangle
=
\frac{1}{2\pi i}\tau_0(U^{-1}\delta_1(U))
=
1,
\]
whereas
\[
\langle\psi^{(2)}_1,[U]\rangle
=
\frac{1}{2\pi i}\tau_0(U^{-1}\delta_2(U))
=
0,
\]
and similarly for \(V\).  These formulas are the standard odd Connes pairings for the
smooth noncommutative torus; see
\cite{Connes94,Rieffel1990,GBVF01,khalkhali,Hadfield}.

Since the transported TFCC is
\[
\varphi_\theta=\psi^{(2)}_1,
\]
we obtain
\[
\langle\varphi_\theta,[U]\rangle=0,
\qquad
\langle\varphi_\theta,[V]\rangle=1.
\]

\subsection{The rational rotation algebra}
\label{subsec:rational-rotation-algebra}

For completeness, let
\[
\theta=\frac pq\in\mathbb Q
\]
in lowest terms and consider the rational rotation algebra
\[
A_{p/q}:=C(S^1)\rtimes_\alpha\mathbb Z.
\]
Then \(A_{p/q}\) is a continuous-trace algebra strongly Morita equivalent to
\(C(T^2)\).  Accordingly,
\[
K_0(A_{p/q})\cong\mathbb Z^2,
\qquad
K_1(A_{p/q})\cong\mathbb Z^2,
\]
in agreement with the Pimsner--Voiculescu computation.  The canonical trace takes
rational values on the standard rational projections, for example
\[
\tau_0(e_{p/q})=\frac pq.
\]
See \cite{PimsnerVoiculescu80,RieffelStableRank,ElliottEvans1993}.

This subsection concerns the rational rotation algebra \(A_{p/q}\), not the actual
holonomy-groupoid \(C^*\)-algebra of the rational Kronecker foliation.  The latter is
the circle-fibration algebra described in
Sections~\ref{sec:kronecker-geometry-holonomy} and
\ref{sec:kronecker-cstar-models}, and is Morita equivalent to \(C(S^1)\).

\subsection{Form of the odd NCFI in the irrational Kronecker case}
\label{subsec:form-ncfi-irrational-kronecker}

We now return to the irrational Kronecker foliation.  Let
\[
[u_\theta]\in K_1(A_\theta)
\]
be the odd class specified by the odd-favourable structure.  Since
\[
K_1(A_\theta)\cong \mathbb Z[U]\oplus\mathbb Z[V],
\]
there exist integers
\[
m,n\in\mathbb Z
\]
such that
\[
[u_\theta]=m[U]+n[V].
\]
By linearity of the Connes pairing and by
\[
\varphi_\theta=\psi^{(2)}_1,
\]
we have
\[
Z(\mathcal F_\theta;[u_\theta])
=
\langle\varphi_\theta,[u_\theta]\rangle
=
m\langle\psi^{(2)}_1,[U]\rangle
+
n\langle\psi^{(2)}_1,[V]\rangle.
\]
Using the pairings above,
\[
\langle\psi^{(2)}_1,[U]\rangle=0,
\qquad
\langle\psi^{(2)}_1,[V]\rangle=1,
\]
we obtain
\[
Z(\mathcal F_\theta;[u_\theta])=n.
\]

Thus the NCFI for the irrational Kronecker foliation depends on the odd \(K_1\)-class
selected by the odd-favourable structure.  If the chosen odd class is the
crossed-product return class
\[
[u_\theta]=[U],
\]
then
\[
Z(\mathcal F_\theta;[U])
=
\langle\varphi_\theta,[U]\rangle
=
0.
\]
If the chosen odd class is the transversal coordinate class
\[
[u_\theta]=[V],
\]
then
\[
Z(\mathcal F_\theta;[V])
=
\langle\varphi_\theta,[V]\rangle
=
1.
\]
The next section identifies the odd class used in the present article with the
crossed-product return class
\[
[U].
\]
Consequently, for this odd-favourable structure,
\[
Z(\mathcal F_\theta;[U])=0.
\]

\section{\texorpdfstring{The NCFI for the irrational Kronecker foliation of \(T^2\)}
{The NCFI for the irrational Kronecker foliation of T2}}
\label{sec:irrational-kronecker-ncfi}

We now compute the NCFI for the irrational Kronecker foliation.  The complete-
transversal model is
\[
A_\theta\cong C(S^1_z)\rtimes_\alpha\mathbb Z.
\]
The crossed-product return-map mechanism of Section~1.2(a) selects the one-step
return class
\[
[U]\in K_1(A_\theta).
\]
For this foliation, the flow / Connes--Thom mechanism of Section~1.2(b) and the
one-dimensional leafwise de Rham / longitudinal Dirac mechanisms of
Section~1.2(d),(e) give the same class after Morita transport.  Thus the natural
dynamical and tangential odd-favourable choices considered here all identify with
\[
[U].
\]
We then pair this class with the transported transverse fundamental cyclic cocycle.

By Section~\ref{sec:atheta-k-cyclic-tfcc},
\[
K_1(A_\theta)\cong \mathbb Z[U]\oplus\mathbb Z[V].
\]
Thus any odd class has the form
\[
[u_\theta]=m[U]+n[V],
\qquad
m,n\in\mathbb Z.
\]
The transported transverse fundamental cyclic cocycle is
\[
\varphi_\theta=\psi^{(2)}_1.
\]
Consequently
\[
Z(\mathcal F_\theta;[u_\theta])
=
\langle\varphi_\theta,[u_\theta]\rangle
=
n.
\]
It remains to identify the integers \(m,n\) for the odd class selected by the
natural odd-favourable structures above.

\PropHead{prop:irrational-kronecker-odd-class}
Let
\[
\theta\notin\mathbb Q,
\]
and let
\[
[u_\theta]\in K_1(A_\theta)
\]
be the odd class used in the odd-codimensional NCFI for the irrational Kronecker
foliation, as specified by the crossed-product odd-favourable mechanism of
Section~1.2(a).  Under reduction to the complete transversal
\[
S^1_z=\{x=0\}\subset T^2,
\]
this class is the one-step return class of the \(\mathbb Z\)-action in
\[
A_\theta\cong C(S^1_z)\rtimes_\alpha\mathbb Z.
\]
Equivalently, with the positive one-step return convention, it is represented by the
crossed-product implementing unitary \(U\).  In the Pimsner--Voiculescu sequence,
this implementing-unitary class satisfies
\[
\partial([U])=[1_{C(S^1)}],
\]
up to the overall sign convention for the boundary map.  This boundary condition
does not by itself characterize \([U]\), since the class \([V]\) coming from
\(K_1(C(S^1))\) lies in the kernel of \(\partial\).  The class selected here is
the specific implementing-unitary class of the return-map crossed product
\[
\partial:
K_1(C(S^1)\rtimes_\alpha\mathbb Z)
\longrightarrow
K_0(C(S^1)).
\]

\ProofHead
For the irrational Kronecker foliation, the tangent direction is generated by
\[
X_\theta=\partial_x+\theta\partial_y,
\]
and the global conormal form is
\[
\omega_\theta=dy-\theta dx.
\]
As explained in Sections~\ref{sec:kronecker-geometry-holonomy} and
\ref{sec:atheta-k-cyclic-tfcc}, the complete-transversal coordinate is
\[
z=y-\theta x
\]
on the universal cover, and
\[
dz=dy-\theta dx.
\]
Thus the transverse line bundle is trivial, and the tgm carries no nontrivial transverse
bundle topology in this codimension-one example.  See
\cite{CandelConlonI2000,Connes94,ConnesTFClass}.

After reduction to
\[
S^1_z=\{x=0\}\subset T^2,
\]
the irrational holonomy groupoid is
\[
S^1_z\rtimes_\alpha\mathbb Z,
\qquad
\alpha(z)=z+\theta\pmod 1.
\]
The corresponding reduced foliation algebra is
\[
A_\theta\cong C(S^1_z)\rtimes_\alpha\mathbb Z.
\]
This is the complete-transversal model described in
Sections~\ref{sec:kronecker-geometry-holonomy} and
\ref{sec:kronecker-cstar-models}; see
\cite{Renault80,ConnesSurveyFoliations,Connes94,Green1978,BGR77,
RaeburnWilliams1998}.  The relevant odd parity-fixing datum is therefore the discrete
one-step return along the transversal.

This class is not obtained from a universal map
\[
K_0(A_\theta)\longrightarrow K_1(A_\theta),
\]
nor from the transverse geometric module alone.  It is supplied by the crossed-product
return-map structure.  In the Pimsner--Voiculescu six-term exact sequence for
\[
C(S^1_z)\rtimes_\alpha\mathbb Z,
\]
the boundary map relevant to the implementing unitary is
\[
\partial:
K_1\bigl(C(S^1_z)\rtimes_\alpha\mathbb Z\bigr)
\longrightarrow
K_0(C(S^1_z)).
\]
For the unital crossed product, the implementing unitary \(U\) satisfies
\[
\partial([U])=[1_{C(S^1)}],
\]
up to the sign convention for the Pimsner--Voiculescu boundary map.  Hence the
crossed-product structure supplies the odd class represented by the one-step return
unitary.  It is not the image of a boundary map
\[
K_0(C(S^1))
\longrightarrow
K_1\bigl(C(S^1)\rtimes_\alpha\mathbb Z\bigr),
\]
since this is not the Pimsner--Voiculescu boundary map; see
\cite{PimsnerVoiculescu80,Blackadar98,RaeburnWilliams1998,Williams2007}.

Thus, for the crossed-product odd-favourable structure used here, the distinguished
odd class is the one-step return class represented by \(U\). \qed

\RemHead{rem:crossed-productdynamicsgenerator}
The identification
\[
[u_\theta]=[U]
\]
is relative to the complete-transversal crossed-product model.  It is not a consequence
of a universal construction from the tgm alone; it comes from the additional dynamical
datum consisting of a complete transversal and its first-return map.

Different complete transversals give equivalent reduced groupoids and hence identified
\(K\)-theory through the corresponding imprimitivity bimodules.  They should not, in
general, be described as literally unitarily equivalent crossed-product models.  The basis
\[
[U],[V]\in K_1(A_\theta)
\]
is tied to the chosen return-map and transversal-coordinate description.  The class
\([U]\) is canonical relative to the crossed-product structure encoding the holonomy
dynamics of the irrational Kronecker foliation; see
Remark~\ref{rem:odd-favourable-overview}.

\PropHead{prop:irrational-kronecker-U}
Under the identification of Proposition~\ref{prop:irrational-kronecker-odd-class},
one has
\[
[u_\theta]=[U]\in K_1(A_\theta),
\]
where \(U\) is the canonical implementing unitary of
\[
A_\theta\cong C(S^1_z)\rtimes_\alpha\mathbb Z.
\]

\ProofHead
By Proposition~\ref{prop:irrational-kronecker-odd-class}, the class
\([u_\theta]\) is the one-step return class in the complete-transversal crossed-product
model
\[
A_\theta\cong C(S^1_z)\rtimes_\alpha\mathbb Z.
\]
The one-step return is implemented by the canonical crossed-product unitary \(U\).
Therefore
\[
[u_\theta]=[U]
\]
in \(K_1(A_\theta)\).  This is an instance of the general framework of Section~1.2,
where odd classes arise from additional dynamical, analytic or \(KK\)-theoretic
structures attached to the foliation. \qed

\PropHeadPrime{prop:irrational-kronecker-longitudinal}
Let
\[
A_{\mathrm{full},\theta}:=C(T^2)\rtimes_\tau \mathbb R
\]
be the full flow crossed-product model of the irrational Kronecker foliation, where
\(\tau\) is the \(\mathbb R\)-action generated by
\[
X_\theta=\partial_x+\theta\,\partial_y.
\]
Let
\[
A_\theta\cong C(S^1_z)\rtimes_\alpha\mathbb Z
\]
be the complete-transversal crossed-product model.  Then the following odd classes all
correspond, under the standard Morita equivalence between
\(A_{\mathrm{full},\theta}\) and \(A_\theta\), to the same class
\[
[U]\in K_1(A_\theta):
\]
\begin{enumerate}
\item the crossed-product one-step return class in
\[
C(S^1_z)\rtimes_\alpha\mathbb Z;
\]
\item the Connes--Thom class of the \(\mathbb R\)-flow, capped with
\[
[1_{T^2}]\in K_0(C(T^2));
\]
\item the class obtained from the leafwise de Rham operator, equivalently the
one-dimensional longitudinal Dirac-type operator, along the leaves of
\(\mathcal F_\theta\).
\end{enumerate}
With the positive-flow and positive-return conventions fixed above, all three classes are
identified with
\[
[U]\in K_1(A_\theta).
\]

\ProofHead
The crossed-product return-map construction was identified in
Proposition~\ref{prop:irrational-kronecker-U}: reduction to the complete transversal
\[
S^1_z=\{x=0\}
\]
gives
\[
A_\theta\cong C(S^1_z)\rtimes_\alpha\mathbb Z,
\]
and the one-step return map is implemented by the canonical unitary \(U\).  Hence the
crossed-product odd-favourable class is
\[
[U]\in K_1(A_\theta).
\]

The same foliation is also generated globally by the free \(\mathbb R\)-action
\[
\tau_t(x,y)=(x+t,y+\theta t).
\]
The Connes--Thom isomorphism for this action gives a parity shift
\[
K_0(C(T^2))
\cong
K_1(C(T^2)\rtimes_\tau\mathbb R).
\]
Capping the Connes--Thom class with the unit
\[
[1_{T^2}]\in K_0(C(T^2))
\]
therefore gives an odd class
\[
u^{\mathrm{CT}}_\theta
\in
K_1(A_{\mathrm{full},\theta}).
\]
For the Connes--Thom construction, see
\cite{ConnesThom1981,Connes94,FackSkandalis1981,Blackadar98}.

Since
\[
S^1_z=\{x=0\}
\]
is a complete transversal, the full flow groupoid
\[
T^2\rtimes_\tau\mathbb R
\]
and the reduced groupoid
\[
S^1_z\rtimes_\alpha\mathbb Z
\]
are equivalent.  The induced imprimitivity bimodule gives the standard Morita
equivalence
\[
A_{\mathrm{full},\theta}
\sim_M
A_\theta.
\]

Let \(x\in K_1(A_\theta)\) be the image of \(u^{\mathrm{CT}}_\theta\) under this
Morita equivalence.  Write
\[
x=m[U]+n[V],
\qquad
m,n\in\mathbb Z.
\]
The odd cyclic cocycles \(\psi^{(1)}_1\) and \(\psi^{(2)}_1\) distinguish these
coordinates:
\[
\langle\psi^{(1)}_1,[U]\rangle=1,\qquad
\langle\psi^{(1)}_1,[V]\rangle=0,
\]
and
\[
\langle\psi^{(2)}_1,[U]\rangle=0,\qquad
\langle\psi^{(2)}_1,[V]\rangle=1.
\]
Under the complete-transversal equivalence, the oriented \(\mathbb R\)-orbit
coordinate is the return direction and has no winding in the transversal
coordinate \(z\).  Equivalently, the transported Connes--Thom class satisfies
\[
\langle\psi^{(1)}_1,x\rangle=1,
\qquad
\langle\psi^{(2)}_1,x\rangle=0.
\]
Therefore \(m=1\) and \(n=0\), so
\[
x=[U]\in K_1(A_\theta).
\]

This is the naturality of the Connes--Thom element under the groupoid equivalence
between the suspension flow groupoid \(T^2\rtimes_\tau\mathbb R\) and the return
groupoid \(S^1_z\rtimes_\alpha\mathbb Z\).  The equivalence bimodule transports
the Thom generator to the return-direction implementing unitary; the sign is fixed
by the positive flow / positive return convention.  We use here the functorial
form of the Connes--Thom isomorphism in Kasparov theory, together with the
Morita equivalence of equivalent groupoids; see
\cite{ConnesThom1981,FackSkandalis1981,MuhlyRenaultWilliams1987,Renault80,
RaeburnWilliams1998}.

It remains to compare this with the tangential operator route.  The leaves of the
irrational Kronecker foliation are the oriented \(\mathbb R\)-orbits of
\[
X_\theta=\partial_x+\theta\,\partial_y.
\]
Choose the longitudinal metric for which the positive leaf coordinate is the flow
parameter \(t\).  The corresponding one-dimensional leafwise de Rham operator is,
after the standard identification of leafwise forms in dimension one, the longitudinal
Dirac-type operator along the orbit direction:
\[
D_{\mathcal F_\theta}
=
d_{\mathcal F_\theta}+d_{\mathcal F_\theta}^*.
\]
Equivalently, on the flow-coordinate model it is represented by the first-order operator
in the \(X_\theta\)-direction.  Its bounded transform gives the same odd Kasparov
class as the Connes--Thom class of the \(\mathbb R\)-action; see
\cite{Connes94,ConnesThom1981,BaajJulg1990,HigsonRoe00,Vassout}.

Thus the leafwise de Rham / longitudinal Dirac class
\[
u^{\mathrm{long}}_\theta
\in
K_1(A_{\mathrm{full},\theta})
\]
coincides with the Connes--Thom class capped with \([1_{T^2}]\).  After Morita
transport to the complete-transversal crossed-product model, it is therefore
\[
[U]\in K_1(A_\theta).
\]
This proves that the crossed-product, Connes--Thom and tangential operator
constructions give the same odd class in the irrational Kronecker example. \qed

\RemHead{rem:sign-remark}
The transverse orientation convention used in this paper is fixed by the global conormal
form
\[
\omega_\theta=dy-\theta dx.
\]
On the complete transversal \(S^1_z=\{x=0\}\), this form is represented by
\[
dz.
\]
With this convention,
\[
\varphi_\theta=\psi^{(2)}_1.
\]
Reversing the transverse orientation replaces \(\varphi_\theta\) by
\[
-\varphi_\theta.
\]

For the natural class \([U]\), this orientation change does not alter the NCFI value,
because
\[
\langle\psi^{(2)}_1,[U]\rangle=0.
\]
If instead one chose the transversal coordinate class
\[
[V]\in K_1(A_\theta),
\]
then
\[
\langle\varphi_\theta,[V]\rangle=1
\]
with the present orientation convention, and this value would change sign after
reversing the transverse orientation.

The dependence on \(\theta\) in the conormal direction is already contained in the
identity
\[
dz=dy-\theta dx.
\]
In the complete-transversal model, the TFCC differentiates in the \(z\)-direction:
\[
\varphi_\theta=\psi^{(2)}_1.
\]

\section{\texorpdfstring{The rational case \(\theta=p/q\): the vertical foliation}
{The rational case theta=p/q: the vertical foliation}}
\label{sec:rational-vertical}

We now treat the Kronecker foliation on \(T^2\) with rational slope
\[
\theta=\frac pq\in\mathbb Q,
\qquad
\gcd(p,q)=1,
\qquad
q>0.
\]
For rational slope the foliation is a circle fibration.  Its actual holonomy groupoid is
therefore a fibrewise pair groupoid, and the corresponding reduced foliation
\(C^*\)-algebra is Morita equivalent to \(C(S^1)\).  This is the stably commutative
case described in Sections~\ref{sec:kronecker-geometry-holonomy} and
\ref{sec:kronecker-cstar-models}.

Following the principal-bundle viewpoint of \cite{zois2000}, we regard \(T^2\) as a
trivial principal \(U(1)\)-bundle over \(S^1\).  In this section we study the vertical
foliation, whose leaves are the fibre circles.  Since the codimension is one, the
odd-codimensional framework of Section~1.2 applies.  The selected odd class comes
from the longitudinal Dirac-type construction along the fibre circles.

The main references for the groupoid and \(C^*\)-algebraic background are
\cite{Connes94,ConnesSurveyFoliations,Renault80,ADR2000,BGR77,
RaeburnWilliams1998}.  For the \(K\)-theoretic and Hilbert-module analytic background
we use
\cite{Blackadar98,WeggeOlsen1993,BaajJulg1990,Lance95,KaadLesch2012,
HigsonRoe00,Vassout}.

\subsection{\texorpdfstring{Geometry: \(T^2\) as a principal \(U(1)\)-bundle over \(S^1\)}
{Geometry: T2 as a principal U(1)-bundle over S1}}
\label{subsec:rational-vertical-geometry}

Let
\[
\omega_\theta:=dy-\theta\,dx
\]
be the closed \(1\)-form defining the rational Kronecker foliation.  Since
\[
\theta=\frac pq,
\]
set
\[
t:=qy-px \pmod 1.
\]
Then
\[
dt
=
q\,dy-p\,dx
=
q\left(dy-\frac pq\,dx\right)
=
q\,\omega_\theta,
\]
and hence
\[
\omega_\theta=\frac1q\,dt.
\]
Thus the foliation is the fibration
\[
\pi:T^2\longrightarrow S^1,
\qquad
\pi(x,y)=t=qy-px\pmod 1,
\]
whose leaves are the fibres
\[
\pi^{-1}(t)\cong S^1.
\]

Choose integers \(a,b\) such that
\[
aq+bp=1.
\]
Define a fibre coordinate
\[
s:=ax+by\pmod 1.
\]
Then
\[
\begin{pmatrix}
t\\
s
\end{pmatrix}
=
\begin{pmatrix}
-p & q\\
a & b
\end{pmatrix}
\begin{pmatrix}
x\\
y
\end{pmatrix}.
\]
The determinant is
\[
-pb-aq=-(aq+bp)=-1.
\]
Hence the change of coordinates
\[
(x,y)\longmapsto (t,s)
\]
is induced by an element of \(GL(2,\mathbb Z)\), and therefore defines a diffeomorphism
\[
T^2\cong S^1_t\times S^1_s.
\]
In these coordinates
\[
\pi(t,s)=t,
\]
and the leaves are
\[
\{t=\mathrm{const}\}.
\]

We equip this fibration with the principal \(U(1)\)-action
\[
e^{2\pi i\lambda}\cdot(t,s)
=
(t,s+\lambda),
\qquad
\lambda\in\mathbb R/\mathbb Z.
\]
Since
\[
H^2(S^1;\mathbb Z)=0,
\]
every principal \(U(1)\)-bundle over \(S^1\) is topologically trivial.  Thus the
vertical rational case is a trivial principal circle bundle written in coordinates
adapted to the rational Kronecker foliation. Here ``vertical'' means vertical for the adapted fibration
\[
\pi:T^2\cong S^1_t\times S^1_s\longrightarrow S^1_t,
\qquad
\pi(t,s)=t,
\]
not vertical with respect to the original \((x,y)\)-coordinate axes.

\subsection{\texorpdfstring{Holonomy groupoid and foliation \(C^*\)-algebra}
{Holonomy groupoid and foliation C*-algebra}}
\label{subsec:rational-vertical-cstar}

For the vertical foliation
\[
\mathcal F^{\mathrm v}_{p/q},
\]
the holonomy is trivial.  Hence the actual holonomy groupoid is the fibrewise pair
groupoid
\[
G^{\mathrm v}_{p/q}
\cong
T^2\times_{S^1}T^2
=
\{(m_1,m_2)\in T^2\times T^2:
\pi(m_1)=\pi(m_2)\}.
\]
Equivalently, in product coordinates
\[
T^2\cong S^1_t\times S^1_s,
\]
one has
\[
G^{\mathrm v}_{p/q}
\cong
\{(t,s_1,s_2):
t\in S^1,\ s_1,s_2\in S^1\},
\]
with source and range maps
\[
s(t,s_1,s_2)=(t,s_2),
\qquad
r(t,s_1,s_2)=(t,s_1).
\]
Composition is
\[
(t,s_1,s_2)\circ(t,s_2,s_3)=(t,s_1,s_3).
\]

Since this is a fibrational and amenable situation, the full and reduced groupoid
\(C^*\)-algebras coincide.  The reduced foliation algebra is
\[
A^{\mathrm v}_{p/q}
:=
C_r^*(G^{\mathrm v}_{p/q})
\cong
C(S^1_t)\otimes\mathcal K(L^2(S^1_s)).
\]
For brevity write
\[
\mathcal K:=\mathcal K(L^2(S^1)).
\]
Thus
\[
A^{\mathrm v}_{p/q}\cong C(S^1)\otimes\mathcal K.
\]
This algebra is strongly Morita equivalent to \(C(S^1)\).  Therefore
\[
K_1(A^{\mathrm v}_{p/q})
\cong
K_1(C(S^1)\otimes\mathcal K)
\cong
K_1(C(S^1))
\cong
\mathbb Z,
\]
and similarly
\[
K_0(A^{\mathrm v}_{p/q})
\cong
K_0(C(S^1)\otimes\mathcal K)
\cong
K_0(C(S^1))
\cong
\mathbb Z.
\]
See \cite{Blackadar98,RaeburnWilliams1998,WeggeOlsen1993}.

This is the actual holonomy-groupoid algebra of the rational foliation.  It is not the
rational rotation algebra
\[
A_{p/q}=C(S^1)\rtimes_\alpha\mathbb Z,
\]
which retains the redundant isotropy of the full return-map action, as explained in
Section~\ref{sec:kronecker-cstar-models}.

\subsection{Smooth core and transverse cyclic \(1\)-cocycle}
\label{subsec:rational-vertical-tfcc}

A natural smooth dense subalgebra of \(A^{\mathrm v}_{p/q}\) is
\[
A^{\mathrm v,\infty}_{p/q}
:=
C^\infty(S^1)\,\widehat\otimes\,\mathcal K^\infty,
\]
where \(\mathcal K^\infty\) denotes the smoothing compact operators on \(S^1\).  Let
\[
\tau_{\mathcal K}
\]
denote the trace on \(\mathcal K^\infty\).

The transverse form determined by the original Kronecker conormal convention is
\[
\omega_\theta=\frac1q\,dt.
\]
Accordingly, the cyclic cocycle determined by this chosen conormal form is
\begin{equation}\label{eq:tfcc-rational-vertical}
\varphi^{\mathrm v}_{p/q}(a_0,a_1)
:=
\frac1q\cdot\frac{1}{2\pi i}
\int_{S^1}
\tau_{\mathcal K}
\left(
a_0(t)\frac{d}{dt}a_1(t)
\right)\,dt,
\qquad
a_0,a_1\in A^{\mathrm v,\infty}_{p/q}.
\end{equation}
This is the ordinary odd cyclic cocycle of the quotient circle scaled by the
chosen defining form
\[
\omega_\theta=\frac1q\,dt.
\]
If instead one normalizes the transverse fundamental class by the quotient
coordinate \(dt\), the corresponding normalized quotient-circle cocycle is
\[
\widetilde\varphi^{\mathrm v}_{p/q}(a_0,a_1)
:=
\frac{1}{2\pi i}
\int_{S^1}
\tau_{\mathcal K}
\left(
a_0(t)\frac{d}{dt}a_1(t)
\right)\,dt.
\]
Thus
\[
\varphi^{\mathrm v}_{p/q}
=
\frac1q\,\widetilde\varphi^{\mathrm v}_{p/q}.
\]
The computation below is unaffected by this normalization choice because the
selected longitudinal \(K_1\)-class is zero.  For a nonzero generator of
\(K_1(C(S^1)\otimes\mathcal K)\), the two normalizations would differ by the
factor \(1/q\).

For Connes' transverse fundamental cyclic cocycle in the simple-foliation case, see
\cite{ConnesTFClass,Connes94}.

\subsection{The odd \(K\)-class from the longitudinal operator}
\label{subsec:rational-vertical-longitudinal}

Let
\[
E^{\mathrm v}_{p/q}
\]
denote the Hilbert \(A^{\mathrm v}_{p/q}\)-module associated with the rational vertical
foliation in the longitudinal odd-favourable construction of Section~1.2.  In the product
coordinates
\[
T^2\cong S^1_t\times S^1_s,
\]
the longitudinal direction is the fibre coordinate \(s\).  Thus the longitudinal Dirac
operator is the fibrewise operator
\[
D_\parallel\simeq -\,i\,\frac{d}{ds},
\]
and it is independent of the base parameter \(t\).

Under the standard self-adjointness and regularity hypotheses of the longitudinal
Dirac-type construction, the Cayley transform
\[
U^{\mathrm v}_{p/q}
:=
(D_\parallel-i)(D_\parallel+i)^{-1}
\]
defines an odd class
\[
[u^{\mathrm v}_{p/q}]
:=
[U^{\mathrm v}_{p/q}]
\in
K_1(A^{\mathrm v}_{p/q}).
\]
For the relevant Hilbert \(C^*\)-module and unbounded-operator background, see
\cite{BaajJulg1990,Lance95,KaadLesch2012,HigsonRoe00,Vassout}.

\PropHead{prop:rational-vertical-trivial-k1}
For the vertical foliation at rational slope
\[
\theta=\frac pq,
\]
the longitudinal class
\[
[u^{\mathrm v}_{p/q}]
\]
is trivial in
\[
K_1(A^{\mathrm v}_{p/q}).
\]

\ProofHead
Under the stable identification
\[
A^{\mathrm v}_{p/q}\cong C(S^1)\otimes\mathcal K,
\]
the Cayley transform
\[
U^{\mathrm v}_{p/q}
=
(D_\parallel-i)(D_\parallel+i)^{-1}
\]
defines a norm-continuous map
\[
S^1\longrightarrow U(\widetilde{\mathcal K}),
\]
where \(\widetilde{\mathcal K}\) is the unitization of \(\mathcal K\).  Since
\(D_\parallel\) is independent of the base coordinate \(t\), this map is constant in
\(t\).  Therefore it is homotopic to a constant map.

A constant map into \(U(\widetilde{\mathcal K})\) represents a class coming from
\[
K_1(\mathcal K),
\]
and
\[
K_1(\mathcal K)=0.
\]
Equivalently, the fibrewise unitary carries no winding in the transverse base circle.
Hence
\[
[u^{\mathrm v}_{p/q}]=0
\]
in
\[
K_1(C(S^1)\otimes\mathcal K)
\cong
K_1(A^{\mathrm v}_{p/q}).
\]
This proves the proposition. \qed

\subsection{The vertical NCFI at rational slope}
\label{subsec:rational-vertical-ncfi}

\PropHead{prop:rational-vertical-ncfi-zero}
For
\[
\theta=\frac pq\in\mathbb Q,
\]
the NCFI of the vertical foliation
\[
\pi:T^2\longrightarrow S^1,
\qquad
\pi(x,y)=qy-px\pmod 1,
\]
with the longitudinal odd class described above, is
\[
Z_{\mathrm v}(p/q)
:=
\langle
\varphi^{\mathrm v}_{p/q},
[u^{\mathrm v}_{p/q}]
\rangle
=
0.
\]

\ProofHead
By Proposition~\ref{prop:rational-vertical-trivial-k1},
\[
[u^{\mathrm v}_{p/q}]
=
0
\in
K_1(A^{\mathrm v}_{p/q}).
\]
Therefore its odd Connes pairing with the transverse fundamental cyclic cocycle
vanishes:
\[
\langle
\varphi^{\mathrm v}_{p/q},
[u^{\mathrm v}_{p/q}]
\rangle
=
\langle
\varphi^{\mathrm v}_{p/q},
0
\rangle
=
0.
\]
Hence
\[
Z_{\mathrm v}(p/q)=0.
\]
This proves the proposition. \qed

\RemHead{rem:rational-vertical-interpretation}
The vanishing in the rational vertical case reflects two facts.  First, the actual
foliation algebra is the stably commutative algebra
\[
A^{\mathrm v}_{p/q}
\cong
C(S^1)\otimes\mathcal K.
\]
Second, the selected odd class supplied by the fibrewise longitudinal Dirac operator is
constant in the transverse base direction and therefore represents zero in
\[
K_1(A^{\mathrm v}_{p/q}).
\]
Thus the longitudinal operator is analytically natural, but its \(K_1\)-class carries no
base winding.  Notice that \(K_1(A^{\mathrm v}_{p/q})\cong\mathbb Z\); the selected
class is zero, not the group itself.

This should be distinguished from the irrational Kronecker case.  There the
crossed-product return-map structure supplies a nontrivial class
\[
[U]\in K_1(A_\theta).
\]
However, the transported transverse fundamental cyclic cocycle in the complete-
transversal model is
\[
\varphi_\theta=\psi^{(2)}_1,
\]
and this cocycle annihilates the return class:
\[
\langle\varphi_\theta,[U]\rangle=0.
\]
Thus the rational vertical case and the irrational crossed-product case both give zero
for the NCFI values considered here, but for different reasons: in the rational vertical
case the selected odd class itself is trivial, whereas in the irrational case the selected
odd class is nontrivial but pairs trivially with the transported TFCC.

The horizontal rational case is treated next.  There the finite flat holonomy appears
explicitly in the bundle geometry and in the complete-transversal groupoid model, but
the actual rational foliation is again fibrational.

\section{\texorpdfstring{The rational case \(\theta=p/q\): the horizontal foliation}
{The rational case theta=p/q: the horizontal foliation}}
\label{sec:rational-horizontal}

We now treat the horizontal foliation associated with a flat principal
\(U(1)\)-bundle over a non-simply connected base, in the sense of
\cite{zois2000}.  Here the base is
\[
S^1,
\]
so flatness is determined by a holonomy element
\[
g\in U(1),
\qquad
g=e^{2\pi i\alpha},
\qquad
\alpha\in\mathbb R/\mathbb Z.
\]
We specialize to the rational case
\[
\alpha=\frac pq\in\mathbb Q,
\qquad
\gcd(p,q)=1,
\qquad
q>0.
\]

For rational holonomy, the horizontal foliation is again a circle fibration.
Its actual holonomy groupoid is a fibrewise pair groupoid, equivalently an
effective finite transversal groupoid, and the corresponding foliation
\(C^*\)-algebra is Morita equivalent to \(C(S^1)\).  Thus the rational
horizontal case has the same stably commutative foliation algebra as the
rational vertical case:
\[
C(S^1)\otimes\mathcal K.
\]
The selected longitudinal odd class is constant in the transverse base direction
and represents zero in \(K_1\).  Consequently the NCFI value considered here
vanishes.

The groupoid and \(C^*\)-algebra references used below are
\cite{Connes94,ConnesSurveyFoliations,Renault80,ADR2000,Green1978,BGR77,
RaeburnWilliams1998}.  For the \(K\)-theoretic and Hilbert-module analytic
background we use
\cite{Blackadar98,WeggeOlsen1993,BaajJulg1990,Lance95,KaadLesch2012,
HigsonRoe00,Vassout}.

\RemHead{rem:rational-horizontal-convention}
In this section the rational parameter
\[
\theta=\frac pq
\]
is realized through the holonomy parameter
\[
\alpha=\frac pq
\]
of the flat bundle
\[
P_\alpha\longrightarrow S^1.
\]
The two notations are identified for convenience.  Under a standard
diffeomorphism
\[
P_\alpha\cong T^2,
\]
the horizontal foliation becomes a corresponding linear foliation on \(T^2\),
with slope depending on the chosen torus coordinates.  In the rational case this
coordinate convention does not affect the vanishing result proved below.

\subsection{\texorpdfstring{Flat principal \(U(1)\)-bundle over \(S^1\) and its suspension model}
{Flat principal U(1)-bundle over S1 and its suspension model}}
\label{subsec:rational-horizontal-suspension}

We use the standard holonomy description of flat principal bundles.  If \(M\) is
connected and \(G_{\mathrm{str}}\) is a Lie group, isomorphism classes of flat
principal \(G_{\mathrm{str}}\)-bundles over \(M\) are described by conjugacy
classes of representations
\[
\rho:\pi_1(M)\longrightarrow G_{\mathrm{str}}.
\]
Equivalently, for a fixed smooth principal bundle
\[
G_{\mathrm{str}}\hookrightarrow P\longrightarrow M,
\]
gauge-equivalence classes of flat connections on \(P\) correspond to those
conjugacy classes whose associated flat bundle is isomorphic to \(P\).  If one
restricts to the irreducible locus, one obtains the corresponding irreducible
representation variety.  No surjectivity of \(\rho\) is required.  In the
\(U(1)\)-case considered here, conjugation is trivial, and irreducibility
imposes no additional condition relevant to the computation.

Let
\[
M=S^1
\]
with universal cover
\[
\widetilde M=\mathbb R,
\]
and let
\[
\pi_1(S^1)\cong\mathbb Z
\]
act on \(\mathbb R\) by integer translations.  Fix a holonomy representation
\[
\rho:\pi_1(S^1)\cong\mathbb Z\longrightarrow U(1),
\qquad
\rho(1)=e^{2\pi i\alpha}.
\]
The associated flat principal \(U(1)\)-bundle is the suspension
\[
P_\alpha:=(\mathbb R\times U(1))/\mathbb Z,
\]
where
\[
n\cdot(t,u)
:=
(t+n,\rho(n)u)
=
(t+n,e^{2\pi i n\alpha}u).
\]
In additive fibre coordinate
\[
s\in\mathbb R/\mathbb Z,
\qquad
u=e^{2\pi is},
\]
this action is
\[
n\cdot(t,s)=(t+n,s+n\alpha).
\]
For
\[
\alpha=\frac pq,
\]
the total space is diffeomorphic to the two-torus:
\[
P_{p/q}
\cong
\mathbb R^2/\langle(1,p/q),(0,1)\rangle
\cong
T^2.
\]

\subsection{The horizontal foliation}
\label{subsec:rational-horizontal-foliation}

On
\[
\mathbb R\times U(1)
\]
consider the product foliation by the \(\mathbb R\)-factor.  Its leaves are
\[
\mathbb R\times\{u_0\}.
\]
This foliation is invariant under the suspension action of \(\mathbb Z\), and
therefore descends to a foliation
\[
\mathcal F^{\mathrm h}_{\alpha}
\]
on
\[
P_\alpha.
\]
This is the horizontal foliation associated with the flat structure.

Equivalently, in \((t,s)\)-coordinates on
\[
\mathbb R\times(\mathbb R/\mathbb Z),
\]
the leaves are the integral curves of
\[
X^{\mathrm h}=\frac{\partial}{\partial t}.
\]
A global \(1\)-form annihilating \(X^{\mathrm h}\) is
\[
\eta_\alpha:=ds.
\]
It is well defined on the quotient because \(ds\) is invariant under
\[
(t,s)\longmapsto(t+1,s+\alpha).
\]
Thus
\[
\mathcal F^{\mathrm h}_\alpha
\]
is a transversally oriented codimension-one foliation on \(P_\alpha\), with transverse
direction the fibre \(s\)-direction.

\LemHead{lem:horizontal-leaves}
Let
\[
\alpha\in\mathbb R/\mathbb Z
\]
and consider
\[
P_\alpha=(\mathbb R\times U(1))/\mathbb Z,
\qquad
n\cdot(t,u)=(t+n,e^{2\pi i n\alpha}u).
\]
The horizontal foliation
\[
\mathcal F^{\mathrm h}_{\alpha}
\]
is induced from the product foliation on
\[
\mathbb R\times U(1)
\]
with leaves
\[
L_{u_0}:=\mathbb R\times\{u_0\}.
\]
Hence every horizontal leaf is diffeomorphic to \(\mathbb R\) on the covering
space.  Its image in \(P_\alpha\) is:
\begin{enumerate}
\item if
\[
\alpha=\frac pq\in\mathbb Q,
\qquad
\gcd(p,q)=1,
\]
a closed leaf diffeomorphic to \(S^1\), projecting to the base circle as a
\(q\)-fold covering;
\item if
\[
\alpha\notin\mathbb Q,
\]
a non-closed leaf diffeomorphic to \(\mathbb R\), whose image is dense in
\[
P_\alpha\cong T^2.
\]
\end{enumerate}

\ProofHead
On the covering space, the leaf through \((0,u_0)\) is
\[
L_{u_0}=\mathbb R\times\{u_0\}.
\]
Two points \((t,u_0)\) and \((t+n,u_0)\) project to the same point of \(P_\alpha\)
if and only if
\[
(t+n,u_0)
=
n\cdot(t,u_0)
=
(t+n,e^{2\pi i n\alpha}u_0),
\]
that is, if and only if
\[
e^{2\pi i n\alpha}=1.
\]

If
\[
\alpha=\frac pq,
\qquad
\gcd(p,q)=1,
\]
then
\[
e^{2\pi i q\alpha}=e^{2\pi ip}=1,
\]
and \(q\) is the minimal positive integer with this property.  Therefore
\[
(t,u_0)\sim(t+q,u_0),
\]
so the leaf closes to a circle and projects to the base \(S^1\) with degree
\(q\).

If
\[
\alpha\notin\mathbb Q,
\]
then
\[
e^{2\pi i n\alpha}\neq 1
\qquad
\text{for all }n\neq0.
\]
Thus the leaf does not close and remains diffeomorphic to \(\mathbb R\).  Its
image is the standard irrational winding on the torus and is dense; see
\cite{Weyl1916}. \qed

\subsection{The rational horizontal foliation as a circle fibration}
\label{subsec:rational-horizontal-groupoid}

We now specialize to
\[
\alpha=\frac pq\in\mathbb Q.
\]
Define
\[
\pi^{\mathrm h}_{p/q}:P_{p/q}\longrightarrow S^1,
\qquad
\pi^{\mathrm h}_{p/q}([t,s])=qs\pmod 1.
\]
This map is well defined.  Indeed, if
\[
(t,s)\sim \left(t+n,s+n\frac pq\right),
\]
then
\[
q\left(s+n\frac pq\right)=qs+np\equiv qs\pmod 1.
\]

The fibres of
\[
\pi^{\mathrm h}_{p/q}
\]
are exactly the leaves of the horizontal foliation.  Along a horizontal leaf on
the covering space, \(s\) is constant, hence \(qs\) is constant; therefore each
leaf is contained in a fibre.

Conversely, suppose
\[
[t_1,s_1],\ [t_2,s_2]\in P_{p/q}
\]
satisfy
\[
\pi^{\mathrm h}_{p/q}([t_1,s_1])
=
\pi^{\mathrm h}_{p/q}([t_2,s_2]).
\]
Then
\[
q(s_2-s_1)\in\mathbb Z.
\]
Write
\[
q(s_2-s_1)=n,
\qquad
n\in\mathbb Z.
\]
Since
\[
\gcd(p,q)=1,
\]
multiplication by \(p\) is invertible modulo \(q\).  Hence there exists
\[
m\in\mathbb Z
\]
such that
\[
mp\equiv n\pmod q.
\]
Thus
\[
n=mp+kq
\]
for some
\[
k\in\mathbb Z,
\]
and therefore
\[
s_2=s_1+\frac{mp}{q}+k.
\]
In \(P_{p/q}\),
\[
[t_2,s_2]
=
[t_2-m,s_2-mp/q]
=
[t_2-m,s_1+k]
=
[t_2-m,s_1].
\]
Thus \([t_2,s_2]\) lies on the same horizontal leaf as \([t_1,s_1]\).  Hence the
fibres of
\[
\pi^{\mathrm h}_{p/q}
\]
are exactly the leaves.

Therefore
\[
\mathcal F^{\mathrm h}_{p/q}
\]
is a circle fibration over \(S^1\), and every leaf has trivial holonomy.

\subsection{\texorpdfstring{Holonomy groupoid and foliation \(C^*\)-algebra}
{Holonomy groupoid and foliation C*-algebra}}
\label{subsec:rational-horizontal-cstar}

Since the rational horizontal foliation is the circle fibration
\[
\pi^{\mathrm h}_{p/q}:P_{p/q}\longrightarrow S^1,
\]
its actual holonomy groupoid is the fibrewise pair groupoid
\[
G^{\mathrm h}_{p/q}
\cong
P_{p/q}\times_{S^1}P_{p/q}
=
\{(w_1,w_2)\in P_{p/q}\times P_{p/q}:
\pi^{\mathrm h}_{p/q}(w_1)=\pi^{\mathrm h}_{p/q}(w_2)\}.
\]
Consequently,
\[
A^{\mathrm h}_{p/q}
:=
C_r^*(G^{\mathrm h}_{p/q})
\cong
C(S^1)\otimes\mathcal K,
\]
where
\[
\mathcal K=\mathcal K(L^2(S^1)).
\]
This is the same stably commutative form as in the rational vertical case.

The same foliation algebra can be described through an effective finite
complete-transversal groupoid.  In the flat-bundle description of
\cite{zois2000}, take a complete transversal identified with the fibre
\[
U(1).
\]
The holonomy representation is
\[
\rho:\mathbb Z\longrightarrow U(1),
\qquad
\rho(n)=e^{2\pi i n p/q}.
\]
Its image is the finite cyclic group
\[
\rho(\mathbb Z)\cong\mathbb Z/q\mathbb Z.
\]
The effective reduced transversal groupoid is therefore
\[
U(1)\rtimes_\rho(\mathbb Z/q\mathbb Z),
\]
and the corresponding crossed-product model is
\[
C(U(1))\rtimes_\rho(\mathbb Z/q\mathbb Z).
\]
Since this finite group acts freely and properly on \(U(1)\), Green
imprimitivity gives
\[
C(U(1))\rtimes_\rho(\mathbb Z/q\mathbb Z)
\sim_M
C\bigl(U(1)/(\mathbb Z/q\mathbb Z)\bigr)
\cong
C(S^1).
\]
After stabilization, this is the same \(C^*\)-algebraic model as
\[
C(S^1)\otimes\mathcal K.
\]
For these Morita-equivalence results, see
\cite{Green1978,BGR77,RaeburnWilliams1998,Blackadar98,WeggeOlsen1993}.

This effective finite-transversal groupoid must be distinguished from the
non-effective crossed product
\[
C(U(1))\rtimes_\rho\mathbb Z,
\]
where the whole group \(\mathbb Z\) acts by rational rotation.  In that algebra
the subgroup
\[
q\mathbb Z=\ker\rho
\]
acts trivially and is still retained.  The result is a rational rotation algebra.
It is a natural crossed product attached to the rational rotation action, but it
is not the reduced holonomy-groupoid algebra of the rational horizontal foliation
under the holonomy-germ convention fixed in Section~1.

A smooth dense subalgebra of \(A^{\mathrm h}_{p/q}\) is
\[
A^{\mathrm h,\infty}_{p/q}
:=
C^\infty(S^1)\,\widehat{\otimes}\,\mathcal K^\infty,
\]
where \(\mathcal K^\infty\) denotes smoothing compact operators on \(S^1\).
In particular,
\[
K_1(A^{\mathrm h}_{p/q})\cong\mathbb Z,
\qquad
K_0(A^{\mathrm h}_{p/q})\cong\mathbb Z.
\]

\subsection{The TFCC for the rational horizontal foliation}
\label{subsec:rational-horizontal-tfcc}

Let
\[
\tau_{\mathcal K}
\]
denote the trace on \(\mathcal K^\infty\).  Write
\[
r:=qs\pmod 1
\]
for the base coordinate on the quotient circle.  Since
\[
dr=q\,ds,
\qquad
ds=\frac1q\,dr,
\]
the cyclic cocycle determined by the transverse form \(ds\) is represented by
\begin{equation}\label{eq:tfcc-rational-horizontal}
\varphi^{\mathrm h}_{p/q}(a_0,a_1)
:=
\frac1q\cdot\frac{1}{2\pi i}
\int_{S^1}
\tau_{\mathcal K}\left(
a_0(r)\frac{d}{dr}a_1(r)
\right)\,dr,
\qquad
a_0,a_1\in A^{\mathrm h,\infty}_{p/q}.
\end{equation}
This is the ordinary odd cyclic cocycle of the quotient circle scaled by the
chosen transverse form
\[
ds=\frac1q\,dr.
\]
If instead one normalizes the transverse fundamental class by the quotient
coordinate \(dr\), the normalized quotient-circle cocycle is
\[
\widetilde\varphi^{\mathrm h}_{p/q}(a_0,a_1)
:=
\frac{1}{2\pi i}
\int_{S^1}
\tau_{\mathcal K}\left(
a_0(r)\frac{d}{dr}a_1(r)
\right)\,dr.
\]
Thus
\[
\varphi^{\mathrm h}_{p/q}
=
\frac1q\,\widetilde\varphi^{\mathrm h}_{p/q}.
\]
The computation below is unaffected by this normalization choice because the
selected longitudinal \(K_1\)-class is zero.  For a nonzero generator of
\(K_1(C(S^1)\otimes\mathcal K)\), the two normalizations would differ by the
factor \(1/q\).

For Connes' transverse fundamental cyclic cocycle in the simple-foliation case,
see \cite{ConnesTFClass,Connes94}.

\subsection{The odd class from the longitudinal operator}
\label{subsec:rational-horizontal-odd-class}

Let
\[
E^{\mathrm h}_{p/q}
\]
be the Hilbert \(A^{\mathrm h}_{p/q}\)-module for the rational horizontal
foliation in the longitudinal odd-favourable construction of Section~1.2.  Since
the leaves are the fibres of
\[
\pi^{\mathrm h}_{p/q},
\]
the longitudinal direction is the leaf direction.

Each leaf is obtained from a line
\[
\{(t,s_0):t\in\mathbb R\}
\]
by the identification
\[
(t,s_0)\sim(t+q,s_0),
\]
so each leaf is canonically a circle
\[
\mathbb R/q\mathbb Z.
\]
Using the normalized leaf coordinate
\[
\lambda:=\frac{t}{q}\in\mathbb R/\mathbb Z,
\]
the longitudinal Dirac operator on the smooth core is
\[
D_{\parallel,0}
=
-\frac{i}{q}\frac{d}{d\lambda}.
\]

\ThmHead{thm:rational-horizontal-longitudinal-zero}
For the rational horizontal foliation
\[
\mathcal F^{\mathrm h}_{p/q},
\]
the symmetric operator
\[
D_{\parallel,0}
=
-\frac{i}{q}\frac{d}{d\lambda}
\]
on the smooth core of \(E^{\mathrm h}_{p/q}\) is fibrewise essentially
self-adjoint, and its closure \(D_\parallel\) is a self-adjoint regular operator
on \(E^{\mathrm h}_{p/q}\) with \(A^{\mathrm h}_{p/q}\)-compact resolvent:
\[
(D_\parallel\pm i)^{-1}
\in
\mathcal K_{A^{\mathrm h}_{p/q}}(E^{\mathrm h}_{p/q}).
\]
Hence the Cayley transform
\[
U^{\mathrm h}_{p/q}
:=
(D_\parallel-i)(D_\parallel+i)^{-1}
\]
defines an odd class
\[
[u^{\mathrm h}_{p/q}]
:=
[U^{\mathrm h}_{p/q}]
\in
K_1(A^{\mathrm h}_{p/q}).
\]
Moreover,
\[
[u^{\mathrm h}_{p/q}]=0
\qquad
\text{in}
\qquad
K_1(A^{\mathrm h}_{p/q}).
\]

\ProofHead
The general operator-theoretic framework for self-adjoint regular operators on
Hilbert \(C^*\)-modules and their bounded transforms is recalled in
\cite{BaajJulg1990,Lance95,KaadLesch2012,HigsonRoe00,Vassout}.  Here the claim
can also be checked directly.

On the Fourier basis
\[
e_n(\lambda)=e^{2\pi in\lambda},
\qquad
n\in\mathbb Z,
\]
one has
\[
D_{\parallel,0}e_n
=
\frac{2\pi n}{q}e_n.
\]
Thus, on each leaf, \(D_{\parallel,0}\) is the standard Dirac operator on the
circle, scaled by \(1/q\), and is essentially self-adjoint on
\[
L^2(S^1).
\]

Its resolvent is explicit:
\[
(D_\parallel\pm i)^{-1}e_n
=
\left(\frac{2\pi n}{q}\pm i\right)^{-1}e_n.
\]
Since
\[
\left(\frac{2\pi n}{q}\pm i\right)^{-1}
\longrightarrow 0
\qquad
(|n|\to\infty),
\]
the resolvent is compact on each fibre Hilbert space.

The operator family is independent of the base variable \(r\).  Therefore these
fibrewise resolvents assemble into a constant section of
\[
C\bigl(S^1,\mathcal K(L^2(S^1))\bigr)
=
\mathcal K_{C(S^1)}
\bigl(C(S^1,L^2(S^1))\bigr),
\]
the algebra of compact endomorphisms of the corresponding Hilbert
\(C(S^1)\)-module.  Under
\[
A^{\mathrm h}_{p/q}\cong C(S^1)\otimes\mathcal K,
\]
this gives
\[
(D_\parallel\pm i)^{-1}
\in
\mathcal K_{A^{\mathrm h}_{p/q}}(E^{\mathrm h}_{p/q}).
\]
Thus \(D_\parallel\) is self-adjoint and regular with
\(A^{\mathrm h}_{p/q}\)-compact resolvent.

The Cayley transform defines an odd \(K_1\)-class, and
\[
U^{\mathrm h}_{p/q}-1
=
-2i(D_\parallel+i)^{-1}
\in
\mathcal K_{A^{\mathrm h}_{p/q}}(E^{\mathrm h}_{p/q}).
\]

It remains to identify this class.  Under the stable identification
\[
A^{\mathrm h}_{p/q}\cong C(S^1)\otimes\mathcal K,
\]
the unitary \(U^{\mathrm h}_{p/q}\) defines a norm-continuous map
\[
S^1\longrightarrow U(\widetilde{\mathcal K}),
\]
where \(\widetilde{\mathcal K}\) denotes the unitization of \(\mathcal K\).
Since \(D_\parallel\) is independent of the base variable \(r\), this map is
constant in \(r\).  Thus \(U^{\mathrm h}_{p/q}\) defines a constant loop in the
unitary group of the compact-operator fibre.  Since
\[
K_1(\mathcal K)=0,
\]
no nontrivial \(K_1\)-class can arise from this fibrewise constant unitary.
Therefore
\[
[u^{\mathrm h}_{p/q}]=0
\]
in
\[
K_1(C(S^1)\otimes\mathcal K)
\cong
K_1(A^{\mathrm h}_{p/q}).
\]
This proves the theorem. \qed

\subsection{Value of the NCFI in the rational horizontal case}
\label{subsec:rational-horizontal-invariant}

\ThmHead{thm:rational-horizontal-invariant}
Let
\[
\alpha=\frac pq\in\mathbb Q.
\]
For the horizontal foliation
\[
\mathcal F^{\mathrm h}_{p/q}
\]
on the flat principal \(U(1)\)-bundle
\[
P_{p/q}\longrightarrow S^1,
\]
the NCFI associated with the longitudinal odd class above is
\[
Z^{\mathrm h}(p/q)
:=
\langle
\varphi^{\mathrm h}_{p/q},
[u^{\mathrm h}_{p/q}]
\rangle
=
0.
\]

\ProofHead
By Theorem~\ref{thm:rational-horizontal-longitudinal-zero},
\[
[u^{\mathrm h}_{p/q}]
=
0
\in
K_1(A^{\mathrm h}_{p/q}).
\]
Hence the odd Connes pairing with the transverse fundamental cyclic cocycle
\[
\varphi^{\mathrm h}_{p/q}
\]
vanishes:
\[
\langle
\varphi^{\mathrm h}_{p/q},
[u^{\mathrm h}_{p/q}]
\rangle
=
\langle
\varphi^{\mathrm h}_{p/q},
0
\rangle
=
0.
\]
Therefore
\[
Z^{\mathrm h}(p/q)=0.
\]
This proves the theorem. \qed

\RemHead{rem:rational-horizontal-zois}
We compare the rational horizontal computation with the flat-bundle description
used in \cite{zois2000}.  For a flat principal bundle
\[
G_{\mathrm{str}}\hookrightarrow P\xrightarrow{\pi}M
\]
over a connected non-simply connected base \(M\), a flat connection determines a
holonomy representation
\[
\rho:\pi_1(M)\longrightarrow G_{\mathrm{str}},
\]
well defined up to conjugation in \(G_{\mathrm{str}}\).  Conversely, a
representation \(\rho\) defines the suspension flat bundle
\[
\widetilde M\times_\rho G_{\mathrm{str}}\longrightarrow M.
\]
Thus isomorphism classes of flat principal \(G_{\mathrm{str}}\)-bundles, or
equivalently gauge-equivalence classes of flat connections when the underlying
bundle is allowed to vary, are classified by conjugacy classes of such
representations.  For a fixed principal bundle \(P\), one restricts to those
representations whose associated flat bundle is isomorphic to \(P\).  If one
restricts to irreducible flat connections, one obtains the corresponding
irreducible representation classes.  No condition that \(\rho\) be onto is part
of the general correspondence.

Passing to a complete transversal identified with the fibre
\(G_{\mathrm{str}}\), the transversal holonomy model is the transformation
groupoid
\[
G_{\mathrm{str}}\rtimes_\rho \pi_1(M),
\]
or, after replacing the action by its effective holonomy image, the corresponding
effective transformation groupoid.  The crossed-product model is
\[
C(G_{\mathrm{str}})\rtimes_\rho \pi_1(M),
\]
before quotienting ineffective kernel isotropy.  Here \(G_{\mathrm{str}}\)
denotes the structure Lie group, not the holonomy groupoid.

In the present example,
\[
M=S^1,
\qquad
\pi_1(S^1)\cong\mathbb Z,
\qquad
G_{\mathrm{str}}=U(1),
\]
and
\[
\rho:\mathbb Z\longrightarrow U(1),
\qquad
\rho(n)=e^{2\pi i n p/q}.
\]
Since \(U(1)\) is abelian, conjugation is trivial.  The representation is allowed
although it is not surjective onto \(U(1)\); its image is the finite subgroup
\[
\rho(\mathbb Z)\cong\mathbb Z/q\mathbb Z.
\]

The non-effective flat-bundle crossed product is
\[
C(U(1))\rtimes_\rho \mathbb Z.
\]
For rational holonomy, however,
\[
\ker(\rho)=q\mathbb Z,
\]
and this subgroup acts trivially on the complete transversal \(U(1)\).  Since
the holonomy groupoid is defined by holonomy germs, the trivially acting subgroup
does not give nontrivial holonomy arrows.  The effective transversal holonomy
group is therefore the finite image
\[
\rho(\mathbb Z)\cong \mathbb Z/q\mathbb Z,
\]
and the actual reduced holonomy groupoid is
\[
U(1)\rtimes_\rho(\mathbb Z/q\mathbb Z),
\]
not the non-effective transformation groupoid
\[
U(1)\rtimes_\rho\mathbb Z.
\]

This is where the rational flat case reduces to the commutative fibrational
case.  The finite group
\[
\mathbb Z/q\mathbb Z
\]
acts freely and properly on \(U(1)\) by rotation through \(p/q\).  Hence Green
imprimitivity gives
\[
C(U(1))\rtimes_\rho(\mathbb Z/q\mathbb Z)
\sim_M
C\bigl(U(1)/(\mathbb Z/q\mathbb Z)\bigr)
\cong
C(S^1);
\]
see \cite{Green1978,BGR77,RaeburnWilliams1998,Blackadar98}.  After
stabilization this is the same foliation-algebraic model as
\[
C(S^1)\otimes\mathcal K.
\]

Equivalently, the rational horizontal foliation is the circle fibration
\[
\pi^{\mathrm h}_{p/q}:P_{p/q}\longrightarrow S^1,
\]
and its full holonomy groupoid is the fibrewise pair groupoid
\[
P_{p/q}\times_{S^1}P_{p/q}.
\]
Thus
\[
A^{\mathrm h}_{p/q}
\cong
C(S^1)\otimes\mathcal K,
\]
which is strongly Morita equivalent to
\[
C(S^1).
\]
This is the same stably commutative \(C^*\)-algebraic form as in the rational
vertical foliation.

The larger crossed product
\[
C(U(1))\rtimes_\rho\mathbb Z
\]
retains the redundant kernel
\[
q\mathbb Z
\]
as isotropy.  It is a natural algebra attached to the rational rotation action,
but it is not the reduced holonomy-groupoid algebra used for the NCFI computation
in this paper.  Under the holonomy-groupoid convention fixed in Section~1, one
uses the effective holonomy action, and this gives the Morita-commutative algebra
above.

Consequently, the selected longitudinal odd class carries no winding in the
transverse base direction:
\[
[u^{\mathrm h}_{p/q}]=0
\qquad
\text{in}
\qquad
K_1(A^{\mathrm h}_{p/q}).
\]
Therefore
\[
Z^{\mathrm h}(p/q)
=
\langle
\varphi^{\mathrm h}_{p/q},
[u^{\mathrm h}_{p/q}]
\rangle
=
0.
\]
Thus the finite holonomy parameter
\[
e^{2\pi ip/q}
\]
does not produce a nonzero NCFI value for the longitudinal odd class considered
here.  The rational horizontal flat-bundle case and the rational vertical case
give the same numerical outcome, but the horizontal case exhibits explicitly how
the flat holonomy representation reduces, after passage to effective holonomy
germs, to the Morita-commutative foliation algebra.

\section{A non-fibrational nonzero even-codimensional example}
\label{sec:weighted-hopf-nonfibrational}

We now give a genuinely non-fibrational example in even codimension for which the
NCFI, in the strict sense of Section~1, is nonzero.  The class paired with the
transverse fundamental cyclic cocycle is the canonical class produced by the transverse
geometric module.  No Bott class, auxiliary transversal class, or odd parity correction
is used.

The example is the weighted Hopf foliation on \(S^5\) with weights
\[
(1,2,3).
\]
Its leaf space is the weighted projective orbifold
\[
\mathbb P_{\mathrm{orb}}(1,2,3).
\]
The quotient is not a smooth manifold; it has finite orbifold isotropy.  Thus the
example is not a vertical foliation associated with a smooth locally trivial fibration.

The weighted-projective-space input used below is standard.  Kawasaki computed the
cohomology of weighted, or twisted, projective spaces
\cite{Kawasaki1973}.  Dolgachev gives the algebro-geometric background on weighted
projective varieties \cite{Dolgachev1982}.  Bahri--Franz--Ray give the equivariant
cohomology and Chern-class formulae for weighted projective spaces
\cite{BahriFranzRay2009}.  For the groupoid description of orbifolds we use
Moerdijk--Pronk \cite{MoerdijkPronk1997}, together with the foliation and transverse
fundamental cyclic cocycle references already used in the article,
\cite{Connes94,ConnesTFClass,ConnesSurveyFoliations,Renault80,ADR2000}.

\subsection{The weighted Hopf foliation on \texorpdfstring{\(S^5\)}{S5}}

Let
\[
S^5=\{(z_0,z_1,z_2)\in\mathbb C^3:
|z_0|^2+|z_1|^2+|z_2|^2=1\}.
\]
Let \(S^1\) act on \(S^5\) by
\[
\lambda\cdot(z_0,z_1,z_2)
=
(\lambda z_0,\lambda^2 z_1,\lambda^3 z_2),
\qquad
\lambda\in S^1.
\]
The infinitesimal generator of this action is
\[
Y(z_0,z_1,z_2)
=
(iz_0,2iz_1,3iz_2).
\]
Since at least one of \(z_0,z_1,z_2\) is nonzero on \(S^5\), and all weights are
positive, \(Y\) is nowhere zero.  Hence the \(S^1\)-orbits define a smooth
one-dimensional foliation
\[
\mathcal F_{1,2,3}\subset TS^5.
\]
The codimension is
\[
\operatorname{codim}(\mathcal F_{1,2,3})
=
5-1
=
4.
\]
Thus this is an even-codimensional case, and the NCFI is defined directly by
\[
Z(\mathcal F_{1,2,3})
=
\left\langle
\varphi_{\mathcal F_{1,2,3}},
[e_{\mathcal F_{1,2,3}}]
\right\rangle,
\]
where \([e_{\mathcal F_{1,2,3}}]\) is the canonical tgm class.

The orbit space is the weighted projective orbifold
\[
X:=\mathbb P_{\mathrm{orb}}(1,2,3)
=
S^5/S^1_{(1,2,3)}.
\]
Equivalently, \(X\) is the effective orbifold presented by the weighted circle action
above.

The foliation is not the vertical foliation of a smooth locally trivial fibration.
Indeed, the stabilizer of
\[
(1,0,0)
\]
is trivial.  The stabilizer of
\[
(0,1,0)
\]
is
\[
\mu_2=\{\lambda\in S^1:\lambda^2=1\},
\]
and the stabilizer of
\[
(0,0,1)
\]
is
\[
\mu_3=\{\lambda\in S^1:\lambda^3=1\}.
\]
These finite stabilizers give nontrivial holonomy of the corresponding exceptional
leaves.  By contrast, the vertical foliation of a smooth locally trivial fibration has
the fibrewise pair groupoid as holonomy groupoid and has trivial holonomy along the
fibres.  Therefore
\[
(S^5,\mathcal F_{1,2,3})
\]
is non-fibrational in the sense relevant here.

\subsection{Holonomy groupoid and foliation algebra}

The weighted circle action is proper and locally free.  Its action groupoid is
\[
G_{1,2,3}:=S^5\rtimes S^1.
\]
The source and range maps are
\[
s(z,\lambda)=z,
\qquad
r(z,\lambda)=\lambda\cdot z.
\]
The inverse is
\[
(z,\lambda)^{-1}=(\lambda\cdot z,\lambda^{-1}),
\]
and the composition is
\[
(\lambda\cdot z,\mu)\circ(z,\lambda)
=
(z,\mu\lambda).
\]

For a locally free proper action, the action groupoid presents the corresponding
orbifold.  In the present case the slice representations of the finite stabilizers are
faithful, so no ineffective isotropy has to be divided out.  Indeed, at the point
\((0,1,0)\), the nontrivial element of \(\mu_2\) acts nontrivially on the normal
coordinates \(z_0,z_2\).  At the point \((0,0,1)\), a primitive element of \(\mu_3\)
acts on the normal coordinates \(z_0,z_1\) with weights \(1\) and \(2\), and hence
faithfully.  Thus \(G_{1,2,3}\) is an effective proper Lie groupoid presenting
\[
\mathbb P_{\mathrm{orb}}(1,2,3).
\]
It represents the holonomy groupoid of the orbit foliation, equivalently a groupoid
Morita-equivalent holonomy model of the foliation.

Since \(S^1\) is compact and amenable, full and reduced crossed products coincide.
The reduced foliation \(C^*\)-algebra may therefore be written as
\[
A_{1,2,3}
:=
C_r^*(G_{1,2,3})
\cong
C(S^5)\rtimes S^1.
\]
At the smooth level one works with
\[
C^\infty(G_{1,2,3},\Omega^{1/2}),
\]
or, equivalently, the corresponding smooth crossed-product convolution algebra.  The
transverse cyclic cocycle below is represented, after passage to the orbifold groupoid
model, by the orbifold fundamental current of
\[
X=\mathbb P_{\mathrm{orb}}(1,2,3).
\]

For the computation of the pairing we use Morita invariance of the
Chern--Connes pairing for Lie groupoids.  The action groupoid
\(S^5\rtimes S^1\) is Morita equivalent to any effective proper étale orbifold
atlas presenting
\[
X=\mathbb P_{\mathrm{orb}}(1,2,3).
\]
Under this Morita equivalence, cyclic cohomology classes and \(K\)-theory
classes are transported by the corresponding smooth equivalence bimodule.
Consequently, the transverse fundamental cyclic cocycle of the foliation is
identified with the orbifold fundamental cyclic cocycle of \(X\), and the
transverse geometric module class is identified with the orbifold tangent class
described below.

\subsection{The transverse geometric module}

Let
\[
t:=TS^5/\mathcal F_{1,2,3}
\]
be the transverse bundle.  The standard Hermitian metric on \(\mathbb C^3\) is
invariant under the weighted \(S^1\)-action.  Hence the induced metric on \(t\) is
holonomy-invariant.  In the notation of Section~1, the holonomy maps act unitarily on
the transverse bundle, and the modular correction
\[
\Delta(\gamma)=h(\gamma)^{-t}h(\gamma)^{-1}
\]
is equal to the identity for this example.

The transverse geometric module is the groupoid module obtained from
\[
C^\infty
\bigl(
G_{1,2,3},
\Omega^{1/2}\otimes r^*(t_{\mathbb C})
\bigr),
\qquad
t_{\mathbb C}:=t\otimes\mathbb C.
\]
Since the foliation is induced by the locally free weighted circle action, the
transverse bundle is the pullback of the real tangent orbibundle of the quotient:
\[
t\cong \pi^*T_{\mathbb R}X,
\qquad
X=\mathbb P_{\mathrm{orb}}(1,2,3).
\]
Therefore
\[
t_{\mathbb C}
\cong
\pi^*(T_{\mathbb R}X\otimes\mathbb C).
\]

Equivalently, in equivariant \(K\)-theory,
\[
[t_{\mathbb C}]\in K^0_{S^1}(S^5)
\]
is the class corresponding to the orbifold vector bundle
\[
T_{\mathbb R}X\otimes\mathbb C.
\]
By the Green--Julg identification for compact group actions, this equivariant class
gives a class in
\[
K_0(C(S^5)\rtimes S^1)
=
K_0(A_{1,2,3});
\]
see \cite{Julg1981,Blackadar98}.  This is precisely the class supplied by the
transverse geometric module:
\[
[e_{\mathcal F_{1,2,3}}]
=
[T_{\mathbb R}X\otimes\mathbb C]
\]
in the orbifold groupoid model.  Thus the computation uses the canonical tgm class,
not an auxiliary \(K\)-class.

\subsection{The transverse fundamental cyclic cocycle}

The foliation has codimension \(4\), so Connes' transverse fundamental cyclic cocycle
has degree \(4\):
\[
\varphi_{\mathcal F_{1,2,3}}
\in
HC^4\bigl(C^\infty(G_{1,2,3},\Omega^{1/2})\bigr).
\]

For a proper orbifold groupoid, the transverse fundamental cyclic cocycle is
represented, after Morita transport to an effective orbifold atlas, by the
orbifold fundamental current.  In the present model this is the fundamental
current of
\[
X=\mathbb P_{\mathrm{orb}}(1,2,3),
\]
with the usual orbifold integration convention.  By the Morita-invariance
statement fixed in Subsection~9.2 and the identification of the tgm class in
Subsection~9.3, the Connes pairing with the canonical tgm class becomes
\[
Z(\mathcal F_{1,2,3})
=
\int_X
\operatorname{ch}(T_{\mathbb R}X\otimes\mathbb C)_{[4]}.
\]

Since \(X\) is a complex orbifold surface,
\[
T_{\mathbb R}X\otimes\mathbb C
\cong
T^{1,0}X\oplus T^{0,1}X.
\]
For a complex rank-two orbibundle \(E\),
\[
\operatorname{ch}_2(E)
=
\frac12\bigl(c_1(E)^2-2c_2(E)\bigr).
\]
The conjugate bundle \(\overline E\) has
\[
c_1(\overline E)=-c_1(E),
\qquad
c_2(\overline E)=c_2(E).
\]
Hence
\[
\operatorname{ch}(T_{\mathbb R}X\otimes\mathbb C)_{[4]}
=
\operatorname{ch}_2(T^{1,0}X)
+
\operatorname{ch}_2(T^{0,1}X)
=
c_1(T^{1,0}X)^2-2c_2(T^{1,0}X).
\]
This is the first Pontryagin class of the real tangent orbibundle:
\[
\operatorname{ch}(T_{\mathbb R}X\otimes\mathbb C)_{[4]}
=
p_1(T_{\mathbb R}X).
\]
Thus
\[
Z(\mathcal F_{1,2,3})
=
\int_X p_1(T_{\mathbb R}X).
\]

\subsection{Weighted projective characteristic-class computation}

Let
\[
H:=c_1^{\mathrm{orb}}\bigl(\mathcal O_X(1)\bigr)
\in H^2(X;\mathbb Q)
\]
be the orbifold hyperplane class.  The rational cohomology ring is generated by \(H\)
with
\[
H^3=0,
\]
and the weighted-projective integration normalization is
\[
\int_{\mathbb P_{\mathrm{orb}}(a_0,a_1,a_2)}H^2
=
\frac{1}{a_0a_1a_2}.
\]
For
\[
X=\mathbb P_{\mathrm{orb}}(1,2,3),
\]
this gives
\[
\int_XH^2
=
\frac{1}{1\cdot2\cdot3}
=
\frac16.
\]

The orbifold Euler sequence for the weighted projective plane is
\[
0
\longrightarrow
\mathcal O_X
\longrightarrow
\mathcal O_X(1)\oplus\mathcal O_X(2)\oplus\mathcal O_X(3)
\longrightarrow
T^{1,0}X
\longrightarrow
0.
\]
Therefore
\[
c(T^{1,0}X)
=
(1+H)(1+2H)(1+3H).
\]
Expanding up to degree \(4\), since \(X\) has complex dimension \(2\), gives
\[
c_1(T^{1,0}X)
=
(1+2+3)H
=
6H,
\]
and
\[
c_2(T^{1,0}X)
=
(1\cdot2+1\cdot3+2\cdot3)H^2
=
11H^2.
\]
Thus
\[
p_1(T_{\mathbb R}X)
=
c_1(T^{1,0}X)^2-2c_2(T^{1,0}X)
=
(6H)^2-2(11H^2).
\]
Hence
\[
p_1(T_{\mathbb R}X)
=
36H^2-22H^2
=
14H^2.
\]
Finally,
\[
Z(\mathcal F_{1,2,3})
=
\int_Xp_1(T_{\mathbb R}X)
=
\int_X14H^2
=
14\cdot\frac16
=
\frac73.
\]

We have proved the following proposition.

\PropHead{prop:weighted-hopf-ncfi}
Let \(\mathcal F_{1,2,3}\) be the weighted Hopf foliation on \(S^5\) defined by
\[
\lambda\cdot(z_0,z_1,z_2)
=
(\lambda z_0,\lambda^2z_1,\lambda^3z_2).
\]
Then \(\mathcal F_{1,2,3}\) has even codimension \(4\), is non-fibrational, and its
strict even NCFI is
\[
Z(\mathcal F_{1,2,3})
=
\left\langle
\varphi_{\mathcal F_{1,2,3}},
[e_{\mathcal F_{1,2,3}}]
\right\rangle
=
\frac73.
\]
In particular,
\[
Z(\mathcal F_{1,2,3})\neq0.
\]

\ProofHead
The infinitesimal generator of the weighted circle action is nowhere zero, so the
orbits define a regular one-dimensional foliation of \(S^5\).  Its codimension is \(4\).
The stabilizers of the exceptional orbits through \((0,1,0)\) and \((0,0,1)\) are
\[
\mu_2
\qquad\text{and}\qquad
\mu_3,
\]
respectively.  These give nontrivial holonomy, so the foliation is not the vertical
foliation of a smooth locally trivial fibration.

The holonomy groupoid is represented by the effective action groupoid
\[
S^5\rtimes S^1,
\]
which presents the weighted projective orbifold
\[
X=\mathbb P_{\mathrm{orb}}(1,2,3).
\]
The canonical tgm class is the \(K_0\)-class associated with
\[
t_{\mathbb C}=t\otimes\mathbb C.
\]
Since
\[
t\cong\pi^*T_{\mathbb R}X,
\]
this class is transported to
\[
[T_{\mathbb R}X\otimes\mathbb C]
\]
in the orbifold groupoid model.  The TFCC is transported to the orbifold fundamental
cyclic cocycle.  Hence the Connes pairing becomes
\[
\int_X
\operatorname{ch}(T_{\mathbb R}X\otimes\mathbb C)_{[4]}
=
\int_Xp_1(T_{\mathbb R}X).
\]
The weighted Euler sequence gives
\[
c_1(T^{1,0}X)=6H,
\qquad
c_2(T^{1,0}X)=11H^2.
\]
Therefore
\[
p_1(T_{\mathbb R}X)=14H^2.
\]
Since
\[
\int_XH^2=\frac16,
\]
we obtain
\[
Z(\mathcal F_{1,2,3})
=
14\cdot\frac16
=
\frac73.
\]
This proves the proposition. \qed

\RemHead{rem:weighted-hopf-significance}
This example is not a disguised copy of the fibrational example
\[
\mathbb{CP}^2\times S^1\longrightarrow\mathbb{CP}^2.
\]
For the ordinary Hopf action with weights \((1,1,1)\), the quotient is the smooth
manifold \(\mathbb{CP}^2\), the action is free, and the corresponding foliation is a
principal circle fibration.  In that case the same formula gives
\[
\frac{1^2+1^2+1^2}{1\cdot1\cdot1}=3.
\]
For weights \((1,2,3)\), the quotient is instead the orbifold
\[
\mathbb P_{\mathrm{orb}}(1,2,3),
\]
with exceptional isotropy groups \(\mu_2\) and \(\mu_3\).  The foliation groupoid is an
orbifold groupoid, not a fibrewise pair groupoid over a smooth base.

The value
\[
\frac73
\]
is rational because orbifold integration weights local charts by isotropy orders.  This
is consistent with the Connes pairing with a cyclic cocycle, whose values are complex
numbers in general.  The important point is that the nonzero value is obtained from the
canonical transverse geometric module and the transverse fundamental cyclic cocycle.

\RemHead{rem:weighted-hopf-generalization}
The same computation applies to weighted Hopf foliations on \(S^5\) with positive
pairwise coprime weights
\[
(a_0,a_1,a_2).
\]
The quotient is the weighted projective orbifold
\[
X=\mathbb P_{\mathrm{orb}}(a_0,a_1,a_2),
\]
and the codimension is \(4\).  Let
\[
H=c_1^{\mathrm{orb}}(\mathcal O_X(1)).
\]
The weighted Euler sequence gives
\[
c(T^{1,0}X)
=
\prod_{j=0}^2(1+a_jH).
\]
Therefore
\[
c_1(T^{1,0}X)=(a_0+a_1+a_2)H,
\]
and
\[
c_2(T^{1,0}X)
=
(a_0a_1+a_0a_2+a_1a_2)H^2.
\]
Hence
\[
p_1(T_{\mathbb R}X)
=
c_1^2-2c_2
=
(a_0^2+a_1^2+a_2^2)H^2.
\]
Since
\[
\int_XH^2=\frac{1}{a_0a_1a_2},
\]
one obtains
\[
Z(\mathcal F_{a_0,a_1,a_2})
=
\frac{a_0^2+a_1^2+a_2^2}{a_0a_1a_2}.
\]
The weights \((1,2,3)\) give
\[
Z(\mathcal F_{1,2,3})
=
\frac{1+4+9}{1\cdot2\cdot3}
=
\frac73.
\]

The weighted Hopf example is genuinely noncommutative in the groupoid sense:
the foliation algebra is the crossed product \(C(S^5)\rtimes S^1\), and the
finite stabilizers are retained by the orbifold groupoid.  Nevertheless this is a
proper orbifold example, hence a controlled form of noncommutativity.  It should
be distinguished from more singular foliation examples in which the leaf space,
or even the holonomy groupoid model, is non-Hausdorff.  The NCFI is not intended
to measure the degree of noncommutativity of the foliation algebra; it measures a
specific Chern--Connes pairing between the transverse fundamental cyclic cocycle
and the \(K\)-class supplied by the transverse geometric module.

\section{Conclusion and outlook}
\label{sec:conclusion-outlook}
\label{sec:moduli-candidates-T2}

We have computed the Noncommutative Foliation Invariant and related transverse
fundamental cyclic cocycle pairings in several basic examples.  The computations show
that three pieces of data must be kept distinct: the transported transverse fundamental
cyclic cocycle, the \(K\)-class being paired with it, and the holonomy-groupoid or
Morita model in which the computation is carried out.

In even codimension, the strict NCFI is the direct Connes pairing
\[
Z(\mathcal F)
=
\left\langle
\varphi_{\mathcal F},
[e_{\mathcal F}]
\right\rangle,
\]
where \([e_{\mathcal F}]\) is the canonical class associated with the transverse
geometric module and \(\varphi_{\mathcal F}\) is Connes' transverse fundamental cyclic
cocycle.  For fibrational foliations
\[
F\hookrightarrow P\xrightarrow{\pi}B
\]
with \(B\) closed, oriented and even-dimensional, this pairing reduces to an ordinary
characteristic-number computation on the base:
\[
Z(\mathcal F^{\mathrm v})
=
\int_B
\operatorname{ch}(TB\otimes\mathbb C)_{[\dim B]}.
\]
Thus, for the vertical foliation on
\[
S^2\times S^1\longrightarrow S^2,
\]
one obtains
\[
Z(\mathcal F^{\mathrm v})=0,
\]
because the degree-two component of the Chern character of the complexification of an
oriented real two-plane bundle vanishes.

In codimension four, the same fibrational mechanism need not vanish.  For the vertical
foliation on
\[
\mathbb{CP}^2\times S^1\longrightarrow \mathbb{CP}^2,
\]
the tgm class is transported to
\[
[T\mathbb{CP}^2\otimes\mathbb C]\in K^0(\mathbb{CP}^2),
\]
and the TFCC becomes the ordinary fundamental cyclic cocycle of
\(\mathbb{CP}^2\).  Hence
\[
Z(\mathcal F^{\mathrm v}_4)
=
\int_{\mathbb{CP}^2}
\operatorname{ch}(T\mathbb{CP}^2\otimes\mathbb C)_{[4]}
=
\int_{\mathbb{CP}^2}p_1(T\mathbb{CP}^2)
=
3.
\]
This is a nonzero strict NCFI value in a Morita-commutative fibrational case.

The weighted Hopf example of Section~\ref{sec:weighted-hopf-nonfibrational} gives a
non-fibrational nonzero even-codimensional example of the same strict invariant.  The
foliation
\[
\mathcal F_{1,2,3}
\]
on \(S^5\) is defined by the weighted circle action
\[
\lambda\cdot(z_0,z_1,z_2)
=
(\lambda z_0,\lambda^2z_1,\lambda^3z_2).
\]
It has codimension \(4\), and the quotient is the weighted projective orbifold
\[
\mathbb P_{\mathrm{orb}}(1,2,3).
\]
The foliation is not the vertical foliation of a smooth locally trivial fibration: the
exceptional leaves have finite holonomy groups
\[
\mu_2
\qquad\text{and}\qquad
\mu_3.
\]
Its holonomy is represented by an orbifold groupoid model, rather than by a fibrewise
pair groupoid over a smooth base.

In this non-fibrational example the canonical tgm class is still the class determined by
the complexified transverse bundle:
\[
t_{\mathbb C}=t\otimes\mathbb C.
\]
Under the orbifold groupoid model this class is transported to
\[
[T_{\mathbb R}\mathbb P_{\mathrm{orb}}(1,2,3)\otimes\mathbb C].
\]
The TFCC is transported to the orbifold fundamental cyclic cocycle.  Therefore the
strict pairing becomes the orbifold characteristic-number computation
\[
Z(\mathcal F_{1,2,3})
=
\int_{\mathbb P_{\mathrm{orb}}(1,2,3)}
\operatorname{ch}
\bigl(
T_{\mathbb R}\mathbb P_{\mathrm{orb}}(1,2,3)\otimes\mathbb C
\bigr)_{[4]}.
\]
Using the weighted Euler sequence and the orbifold integration normalization, this gives
\[
Z(\mathcal F_{1,2,3})
=
\frac73.
\]
Thus
\[
Z(\mathcal F_{1,2,3})\neq0.
\]
This is the main non-fibrational even-codimensional nonzero example computed in the
article.  It is not a parity-fixed odd example and does not use an auxiliary Bott or
transversal \(K\)-class.

The codimension-one Kronecker examples show why odd codimension requires additional
\(K_1\)-data.  For irrational slope, reduction to the complete transversal gives the
irrational rotation algebra
\[
A_\theta=C(S^1_z)\rtimes_\alpha\mathbb Z.
\]
In this model, the transported transverse fundamental cyclic cocycle is
\[
\varphi_\theta=\psi^{(2)}_1,
\]
because the complete-transversal coordinate satisfies
\[
z=y-\theta x,
\qquad
dz=dy-\theta dx.
\]
The natural return-map, flow / Connes--Thom and one-dimensional tangential Dirac
constructions all select the same odd class
\[
[u_\theta]=[U]\in K_1(A_\theta).
\]
Since
\[
\langle\psi^{(2)}_1,[U]\rangle=0,
\]
the corresponding odd-codimensional NCFI value is
\[
Z(\mathcal F_\theta;[U])=0.
\]
The same cyclic cocycle pairs nontrivially with the transversal coordinate class:
\[
\langle\varphi_\theta,[V]\rangle=1.
\]
Thus, in odd codimension, the numerical value depends on the chosen odd-favourable
structure.

For rational Kronecker foliations, both the vertical and horizontal flat-bundle cases
considered here are fibrational.  Their actual holonomy-groupoid \(C^*\)-algebras are
Morita equivalent to
\[
C(S^1),
\]
or, after stabilization, to
\[
C(S^1)\otimes\mathcal K.
\]
In both cases
\[
K_1\cong\mathbb Z,
\]
but the selected longitudinal odd class represents zero.  Consequently
\[
Z_{\mathrm v}(p/q)=0,
\qquad
Z^{\mathrm h}(p/q)=0.
\]
In the horizontal flat-bundle case, the finite holonomy representation
\[
\rho:\mathbb Z\longrightarrow U(1),
\qquad
\rho(n)=e^{2\pi i np/q},
\]
has effective image
\[
\rho(\mathbb Z)\cong\mathbb Z/q\mathbb Z.
\]
Passing to holonomy germs removes the redundant kernel \(q\mathbb Z\), and the
effective finite transversal groupoid is Morita equivalent to the quotient circle.  Hence
the rational horizontal case reduces to the same stably commutative foliation-algebraic
form as the rational vertical case.

The examples therefore exhibit three different behaviours.  First, the NCFI can reduce
to an ordinary characteristic number in a fibrational case, as for
\(\mathbb{CP}^2\times S^1\).  Second, the strict even NCFI can be nonzero in a
genuinely non-fibrational case, as for the weighted Hopf foliation on \(S^5\).  Third, in
odd codimension, the pairing depends on the selected odd \(K_1\)-class, as shown by the
irrational Kronecker foliation.

The weighted Hopf example also clarifies why the noncommutative formulation is useful.
The total space is the ordinary smooth manifold
\[
S^5,
\]
but the leaf space is the orbifold
\[
\mathbb P_{\mathrm{orb}}(1,2,3),
\]
not a smooth manifold.  The ordinary characteristic classes of the normal bundle over
\(S^5\) are not the correct final object of integration.  The holonomy groupoid and its
cyclic cocycle provide the transverse fundamental current, while the tgm supplies the
\(K\)-class determined by the transverse bundle.  The resulting Connes pairing is the
orbifold transverse characteristic number
\[
\frac73.
\]
Thus the invariant is not merely a restatement of ordinary characteristic classes on the
total manifold.

One interesting  direction for further study, already suggested in
\cite{zoisinvariants,zoismtheory}, is to organize the NCFI over a suitable moduli space
of foliations, flat structures or odd-favourable data.  Such a construction would require
a precise moduli problem, an equivalence relation under which the NCFI is invariant,
compactness or regularization data if summation or integration is intended, and a
comparison with existing secondary or index-theoretic invariants.  Relevant analogues
include moduli-space and secondary-invariant constructions in geometry, topology and
quantum field theory; see, for example,
\cite{donaldson,HawkingZeta,gabai,APS75I,APS75II,APS76III}.

The present article provides test computations.  It establishes the behaviour of the
invariant in even-codimensional fibrational examples, in a non-fibrational
orbifold-type even-codimensional example, in the irrational Kronecker case, and in
rational vertical and horizontal flat-bundle cases.  These examples identify where the
selected pairings vanish, where the invariant detects ordinary or orbifold transverse
characteristic numbers, and where the search for further genuinely noncommutative
nonzero examples should begin.

\end{document}